\documentclass[reqno,10pt, centertags,draft]{amsart}
\usepackage{amsmath,amsthm,amscd,amssymb,latexsym,upref}
\usepackage{latexsym}

\makeatletter
\def\theequation{\@arabic\c@equation}

\newcommand{\bbN}{{\mathbb{N}}}
\newcommand{\bbR}{{\mathbb{R}}}

\newcommand{\bbC}{{\mathbb{C}}}

\newcommand{\cB}{{\mathcal B}}

\newcommand{\cH}{{\mathcal H}}

\newcommand{\cJ}{{\mathcal J}}

\newcommand{\cM}{{\mathcal M}}

\newcommand{\cU}{{\mathcal U}}

\newcommand{\cX}{{\mathcal X}}
\newcommand{\cY}{{\mathcal Y}}

\newcommand{\no}{\nonumber}
\newcommand{\lb}{\label}

\newcommand{\ol}{\overline}

\newcommand{\wti}{\widetilde  }

\newcommand{\DT}{{\det}_2}

\newcommand{\spec}{\sigma}

\newcommand{\ran}{\text{\rm{ran}}}

\newcommand{\bi}{\bibitem}

\numberwithin{equation}{section}

\newcommand{\loc}{\operatorname{loc}}
\newcommand{\dom}{\operatorname{dom}}
\newcommand{\im}{\operatorname{ran}}

\newcommand{\tr}{\operatorname{tr}}
\newcommand{\re}{\operatorname{Re}}

\renewcommand{\Im}{\text{\rm Im}}

\newcommand{\diag}{\operatorname{diag}}

\theoremstyle{plain}
\newtheorem{theorem}{Theorem}[section]
\newtheorem{hypothesis}[theorem]{Hypothesis}
\newtheorem{lemma}[theorem]{Lemma}
\newtheorem{corollary}[theorem]{Corollary}

\theoremstyle{definition}
\newtheorem{definition}[theorem]{Definition}

\newtheorem{example}[theorem]{Example}

\newtheorem{remark}[theorem]{Remark}

\begin{document}

\allowdisplaybreaks

\title[Evans Functions, Jost Functions, and Fredholm
Determinants]{Evans Functions, Jost Functions, \\ and Fredholm
Determinants}

\author[F.\ Gesztesy, Y.\ Latushkin, and K.\ A.\ Makarov]{Fritz Gesztesy,
Yuri Latushkin, and Konstantin A.\ Makarov}

\address{Department of Mathematics,
University of Missouri, Columbia, MO 65211, USA}
\email{fritz@math.missouri.edu}
\urladdr{http://www.math.missouri.edu/personnel/faculty/gesztesyf.html}
\address{Department of Mathematics, University of
Missouri, Columbia, MO 65211, USA}
\email{yuri@math.missouri.edu}
\urladdr{http://www.math.missouri.edu/personnel/faculty/latushkiny.html}
\address{Department of Mathematics, University of
Missouri, Columbia, MO 65211, USA}
\email{makarov@math.missouri.edu}
\urladdr{http://www.math.missouri.edu/personnel/faculty/makarovk.html}
\date{\today}
\thanks{Based upon work supported by the US National Science
Foundation under Grant Nos.\ DMS-0405526, DMS-0338743, and DMS-0354339,
and by the CRDF grant UP1-2567-OD-03.}

\begin{abstract}
The principal results of this paper consist of an intrinsic definition
of the Evans function in terms of newly introduced generalized
matrix-valued Jost solutions for general first-order matrix-valued
differential equations on the real line, and a proof of the fact that the
Evans function, a finite-dimensional determinant by construction,
coincides with a modified Fredholm determinant associated with a
Birman--Schwinger-type integral operator up to a nonvanishing 
factor.
\end{abstract}

\maketitle

 \section{Introduction} \lb{s1}

In this paper we study connections between the asymptotic behavior
of solutions of the following first-order $d\times d$
matrix-valued systems of linear differential equations on the real line
$\bbR$,
\begin{equation}
y'(x) = A(x)y(x), \quad x\in\bbR,  \label{unperturbed}
\end{equation}
\begin{equation}
y'(x) = (A(x)+R(x))y(x),  \quad x\in\bbR,  \label{perteqn}
\end{equation}
which will often be referred to as the {\em unperturbed} and {\em
perturbed} equation in the following.

Three types of results are obtained: 

$\bullet$ First, we establish the
existence of special matrix-valued solutions of the perturbed
equation \eqref{perteqn}, called the generalized matrix-valued Jost
solutions, that are asymptotic 
to some reference solutions of the unperturbed equation
\eqref{unperturbed}. In spirit, these results resemble the celebrated
Levinson Theorem \cite[Theorem\ 1.3.1]{E}; however, we are using a
different approach involving  Bohl and Lyapunov exponents. 

$\bullet$ Second, we use the determinant of a finite-dimensional matrix
composed of initial data of the generalized matrix-valued Jost solutions
to calculate an infinite-dimensional modified Fredholm determinant of a
Birman--Schwinger-type integral operator associated with equations
\eqref{unperturbed} and \eqref{perteqn}. These results generalize a
classical relation identifying the Jost function and a Fredholm
determinant in the case of one-dimensional half-line  Schr\"odinger
operators first derived by Jost and Pais \cite{JP51}. 

$\bullet$ Third, using the generalized matrix-valued Jost
solutions, we give a coordinate-free definition of the 
Evans function, and relate the Evans function to the above-mentioned
finite-dimensional determinant. Simultaneously, this identifies the Evans
function with an infinite-dimensional modified Fredholm determinant
associated with a Birman--Schwinger-type integral operator associated
with \eqref{unperturbed} and \eqref{perteqn} up to a nonvanishing
factor. As a result, the Evans function is expressed via an
infinite-dimensional modified Fredholm determinant, and, for the special
case of the Schr\"odinger equation, is proved to be equal to the classical
Jost function familiar in scattering theory.

The Evans function, $E$, is a Wronskian-type analytic function
which is widely used to trace the spectrum of ordinary
differential operators, $\mathcal{L}$. Most frequently, the
operators $\mathcal{L}$ appear as linearizations of partial
differential equations along special solutions such as travelling
waves and steady states. For instance, $\mathcal{L}$ can be the
one-dimensional Schr\"{o}dinger operator, $d^2/dx^2+V$, obtained by
linearizing  the reaction-diffusion equation, $u_t=u_{xx}+r(u)$, about a
steady state $u_0$, so that $V(x)=r'(u_0(x))$, $x\in\bbR$.

Recently, the Evans function became one of the
most important tools in stability analysis. It was originally
introduced by J. W. Evans \cite{Ev72}--\cite{Ev75} to treat a particular
model of nerve impulses, however, quite soon this object was generalized,
and numerous connections of the Evans function to many fields of
mathematical physics have been discovered (from Chern numbers in
topology \cite{AGJ} to scattering data in quantum mechanics
\cite{KS}). To give the reader just a small sample of the work done
in this rapidly developing area, we cite, for instance,
\cite{AGJ}, \cite{AGJS97}, 
\cite{BSZ01}, \cite{BD}, \cite{DN03}, \cite{FS}, \cite{Ga93}, \cite{Ga96},
\cite{Ga97}, \cite{GJ89}, \cite{GZ}, \cite{GMWZ}, \cite{HZ04}, \cite{JK}, 
\cite{Ka94}, \cite{KS98}, \cite{KS}, \cite{KS1}, \cite{LZ}, \cite{MW},
\cite{OZ}, \cite{PW1}, \cite{PW93}, \cite{PW2}, \cite{PZ04}, \cite{SS},
\cite{SS1}, \cite{Sa76}, \cite{Sa77}, \cite{SZ01}, \cite{Zu01},
\cite{Zu03}, \cite{ZH98}, and \cite{ZS99}. In addition, excellent reviews
of this subject and further references can be found in \cite{JK} and 
\cite{S}. 

Re-writing the {\it eigenvalue problem} for
$\mathcal{L}$ as a first-order system of differential equations,
\begin{equation}
y'(x)=B(z,x)y(x), \quad x\in\bbR, \lb{ep}
\end{equation}
where $z\in\bbC$ is the spectral
parameter, one is interested in conditions under which this system
has $\bbC^d$-valued solutions $y$ on $\bbR$ exponentially
decaying at $\pm\infty$. The Evans function, $E=E(z)$, is defined in
such a way that $E(z)=0$ if and only if decaying $\bbC^d$-valued
solutions $y$ exist (thus detecting the eigenvalues of $\mathcal{L}$).
Since in this paper we are not concerned with function theoretic
properties of the Evans function such as its analyticity, etc., deferring
this topic to a forthcoming publication, we simply fix a value $z_0$ of
$z$ and study equations \eqref{unperturbed} and \eqref{perteqn} (and hence
\eqref{ep} with $B(z_0,x)=A(x)+R(x)$, $x\in\bbR$). We note that,
according to all definitions available in the literature, the Evans
function for \eqref{unperturbed} and \eqref{perteqn} is {\em not} uniquely
defined, and one of the central issues in the current paper is
to motivate a ``canonical'' choice of the definition of $E$.

To give an informal outline\footnote{See the glossary of notation at
the end of this introduction.} of our results, let us temporarily assume
for simplicity of exposition, that $A$ and $R$  belong to
$L^\infty(\bbR)^{d\times d}$, and, in addition, that $R\in
L^1(\bbR)^{d\times d}$. On the Hilbert space $L^2(\bbR)^d$ consider
operators, $G_A$ and $G_{A+R}$, with domains $H^{1} (\bbR)^d$,
the usual Sobolev space, defined by  
\begin{align}
\begin{split} 
(G_Au)(x)&=-u'(x)+A(x)u(x), \quad x\in\bbR, \\
(G_{A+R}u)(x)&=-u'(x)+(A(x)+R(x))u(x),\quad x\in\bbR.
\end{split}
\end{align}
Next, we represent the perturbation $R$ as a product 
\begin{equation}
R(x)=R_\ell(x)R_r(x), \quad x\in\bbR. 
\end{equation}
Assuming the invertibility of $G_A$ and $G_{A+R}$, one can write
\begin{align}
G_{A+R}^{-1} = G_A^{-1} - G_A^{-1} R_\ell R_r G_{A+R}^{-1},
\lb{1.7}
\end{align}
and, in addition, one can show that the operator $I+ R_r G_A^{-1}
R_\ell$ is invertible. Multiplying \eqref{1.7} from the left with
$R_r$ and solving for $R_r G_{A+R}^{-1}$ then yields
\begin{equation}\label{INVFORMGA}
G_{A+R}^{-1} = G_A^{-1} - G_A^{-1} R_\ell [I  + R_r G_A^{-1}
R_\ell]^{-1} R_r G_A^{-1}.
\end{equation}
Conversely, assuming the invertibility of $G_A$ and $I+ R_r
G_A^{-1} R_\ell$, we conclude that $G_{A+R}$ is invertible and
\eqref{INVFORMGA} holds. Next, still assuming that $G_A$ is invertible,
one computes
\begin{align}
\begin{split}
G_A^{-1/2} G_{A+R} G_A^{-1/2}  &= G_A^{-1/2} (G_A + R_\ell R_r)
G_A^{-1/2}    \\
& = I + G_A^{-1/2} R_\ell R_r G_A^{-1/2},
\end{split}
\end{align}
and assuming that $R_\ell G_A^{-1/2}$ and $R_r G_A^{-1/2}$ are
operators in $\cB_4(L^2(\bbR)^d)$
 (for this it suffices that $R_\ell, R_r
\in L^2(\bbR)^{d\times d}$, cf. Lemma \ref{tech2}; also, see
\cite{GGKr,GK69,Si79} for general information about the
Schatten--von Neumann ideals $\mathcal{B}_p(\cdot)$), one obtains for the
symmetrized perturbation determinant
\begin{align} 
\begin{split}
{\det}_{2} (G_A^{-1/2} G_{A+R} G_A^{-1/2})  
& = {\det}_{2}(I + G_A^{-1/2} R_\ell R_r G_A^{-1/2})   \\ 
& = {\det}_{2}(I + R_r G_A^{-1} R_\ell).  \lb{1.10}
\end{split}
\end{align}
Here ${\det}_{2}(I+T)$ denotes the modified Fredholm determinant of a
Hilbert--Schmidt operator $T$ on $L^2(\bbR)^d$ (and we used the fact that
${\det}_{2}(I+T_1T_2) = {\det}_{2}(I+T_2T_1)$ for $T_1, T_2$ bounded
operators such that $T_1T_2$ and $T_2T_1$ are Hilbert--Schmidt). In
removing the restriction
$A, R \in L^\infty(\bbR)^{d\times d}$ and permitting $A$ to be locally
integrable on $\bbR$ one needs to supply proper operator closures in
formulas \eqref{1.7}--\eqref{1.10} (see, e.g., \cite{GLMZ05} and the
extensive literature therein for a more general discussion of
factorizations of perturbations). 

To make the connection with
\eqref{unperturbed} and \eqref{perteqn} we recall that the operator
$G_A$ is invertible on $L^2(\bbR)^d$ if and only if \eqref{unperturbed} has
an exponential dichotomy $Q$ on $\bbR$ (see, e.g., \cite[Ch.\ 3]{CL}). By
a classical result in the theory of dichotomic differential equations
(see, e.g., \cite[Proposition\ 4.1]{Cop}, \cite[Ch.\ IV]{DK}), since
\eqref{unperturbed} has an exponential dichotomy on $\bbR$ and the
perturbation satisfies $\|R\|_{\bbC^{d\times d}}\in L^1(\bbR)$,
one infers that the perturbed equation \eqref{perteqn} has exponential
dichotomies on both half-lines $\bbR_+$ and $\bbR_-$, and moreover (by yet
another well-known Dichotomy Theorem, see \cite{BAG},
\cite{Pa}, \cite{Palm}, \cite{Sack}), the operator $G_{A+R}$ is then a
Fredholm operator with zero Fredholm index. Thus, $G_{A+R}$ is either
invertible on
$L^2(\bbR)^d$, or has a nontrivial null space (consisting exactly of the 
solutions of \eqref{perteqn} that belong to $H^1(\mathbb{R})^d$).

Due to the assumption $\|R\|_{\bbC^{d\times d}} \in L^1(\bbR)$, the
Birman--Schwinger-type operator 
\begin{equation}
K=R_r G_A^{-1} R_\ell
\end{equation}
is a Hilbert--Schmidt operator on $L^2(\bbR)^d$ (see Lemma \ref{tech2}).
Thus, employing the celebrated  Birman-Schwinger argument we conclude from
\eqref{INVFORMGA} that the operator $G_{A+R}$ is not invertible (that is,
\eqref{perteqn} has solutions in $H^1(\mathbb{R})^d$) if and only if $-1$
is an eigenvalue of $K$, or, in other words, if and only if the
2-modified perturbation determinant $\det_2(I+K)$ is equal to zero. 

Without going into further details at this point, we mention that
the Birman-Schwinger-type argument (cf. \cite[Sect.\ X.III.3]{RS78} for
this terminology) in our present context provides a bijection between the
geometric eigenspace (nullspace) of
$G_{A+R}$ and the geometric eigenspace of $R_r G_A^{-1} R_\ell$
corresponding to the eigenvalue $-1$   (see, e.g., \cite{GH87},
\cite{GLMZ05}, \cite[Sect.\ XIII.3]{RS78},
\cite[Ch.\ III]{Si71}, \cite{Si76} and the recent discussion in
\cite{GLMZ05} which particularly targets non-self-adjoint operators).

Equation \eqref{1.10} finally illustrates the sense in which
$\det_2(I+K)$ can be viewed as a properly symmetrized modified
perturbation determinant (cf.\ \cite[Sect.\ IV.3]{GK69}) associated with
the pair $(G_A, G_{A+R})$. We note at this point that $\det_2(I+K)$ will
be one of the central objects in this paper.

Our next step is to reduce the calculation of the {\em
infinite-dimensional} determinant $\det_2(I+K)$ to
a {\em finite-dimensional} determinant. This is possible because the
integral kernel of the operator $K$ happens to be of a special type: it
belongs to the class of so-called semi-separable integral kernels (cf.\
\cite{GM}, \cite{GGK97}, \cite[Ch.\ XIII]{GGKr}). For the integral operator
with semi-separable integral  kernel we construct a finite-dimensional
determinant, $D$, such that the following representation holds:
\begin{equation}\label{MainRepr}
{\det}_2(I+K)=e^{\Theta}  D,
\end{equation}
where $\Theta$ is a certain constant.
It is known from the work in \cite{GM} that $D$ can be expressed in terms
of solutions of certain Volterra integral equations associated with
\eqref{unperturbed} and \eqref{perteqn}. We make yet another step,
and express $D$ in terms of solutions of \eqref{perteqn} that are
asymptotic to some reference solutions of \eqref{unperturbed}. We
called these solutions of \eqref{perteqn} the {\em generalized
matrix-valued Jost solutions}, since they appear as a
generalization of the classical Jost solutions of the
Schr\"odinger equation \cite[Ch.\ XVII]{CS}, \cite[Sect.\ 12.1]{Ne02}. In
particular, we prove that the determinant $D$ in \eqref{MainRepr} is equal
to the classical Jost function in the case of the perturbed equation 
\eqref{perteqn} corresponding to the eigenvalue problem for the
Schr\"odinger operator $\mathcal{L}$. The reference solutions are chosen
using the Lyapunov exponents while the rate of approximation of the
reference solutions by the generalized matrix-valued Jost
solutions is controlled by the Bohl exponents for
\eqref{unperturbed} (see the terminology in \cite[Ch.\ III]{DK}). The
existence of the generalized matrix-valued Jost solutions is
proved by a contraction mapping argument applied to a ``mixed'' system of
integral equations of Volterra- and Fredholm-type (which appears
to be a new element in the literature). In addition, we identify the
abstract asymptotic properties of the generalized matrix-valued
Jost solutions leading to formula \eqref{MainRepr}. Finally, we
complete the picture by proving that the finite-dimensional
determinant $D$ in \eqref{MainRepr} is equal to the Evans function
when $E$ is defined by means of the initial data of the
generalized matrix-valued Jost solutions, thus making this choice
of $E$ canonical.

We emphasize that our results on the existence of the generalized
matrix-valued Jost solutions can be viewed (cf.\ Remarks
\ref{remLT1} and \ref{remLT2}) as a further refinement of the
celebrated Levinson theorem (see, for instance,
\cite[Theorem\ 1.3.1]{E}).

Some results of the current paper have been announced in
\cite{GLM}. A version of this theory for difference equations, based on the
material developed in this paper, has been derived in \cite{CrL}.

 The paper is organized as follows. In Section
\ref{NAP} we recall some known facts from \cite{Cop,DK} regarding
exponential dichotomies and the Bohl and Lyapunov exponents (which
are our means to measure the asymptotic behavior of solutions of
differential equations), and also show that $K$ is a
Hilbert--Schmidt operator on $L^2(\bbR)^d$ (Lemma \ref{tech2}). In
addition, we formulate Lemmas \ref{stub_dich} and \ref{stab_split}
(proved in Appendix \ref{AppendB}) to the effect that the exponential
dichotomy on half-lines and, more generally, exponential splittings and the
corresponding Bohl exponents for \eqref{unperturbed} persist
under $L^1$-perturbations. Note that the robustness of the Bohl 
exponents can also be proved under assumptions different from
$\|R\|_{\bbC^{d\times d}} \in L^1(\bbR)$ (e.g., assuming that the
perturbation is continuous and decays at $\pm\infty$ to zero in norm), and
thus our subsequent results can be developed in this different setting,
but to keep this manuscript at a reasonable length we do not pursue this
in the current paper.

In Section \ref{secjost} we introduce in Definition \ref{jost} the
matrix-valued Jost solutions, $Y_\pm$ on $\bbR_\pm$, which
are proved to exist (and to be unique) under the assumptions
that \eqref{unperturbed} has an exponential dichotomy on $\bbR$
and $\|R\|_{\bbC^{d\times d}}$ decays exponentially at $\pm\infty$ at
a rate controlled by the width of the two Bohl segments corresponding to
the dichotomy projections (cf.\ Theorems \ref{+} and \ref{-}). Although
this result is obtained under strong exponential decay
assumptions, it is shown to be optimal if one does not
allow further exponential splitting of the dichotomy projections (cf.\ 
Example \ref{ex}).

In Section \ref{PDSEC} (cf.\  Theorem \ref{detcomp}) we establish
formula \eqref{MainRepr} with
\begin{equation}
D={\det}_{\bbC^d}(Y_+(0)+Y_-(0)) 
\end{equation}
under the assumptions of Section
\ref{secjost}. We note that these assumptions
are satisfied for compactly supported perturbations, and our
method of proof of \eqref{MainRepr} is to establish this formula
for a truncated perturbation first, and then pass to the limit as
the support of the truncation expands. Formula \eqref{MainRepr} in
the compactly supported perturbation case (which is our main intermediate 
calculation) is proved using results from \cite{GM} collected in
Appendix \ref{AP1}. We call  $D$ in \eqref{MainRepr} the {\em
Evans determinant}.

Section \ref{SES} deals with the class of perturbed equations
\eqref{perteqn}, where each of the Bohl segments for the dichotomy
projections degenerates to a single point (cf.\ Theorem \ref{detcomt}). In
particular, this assumption covers the case of \eqref{perteqn}
corresponding to the Schr\"odinger equation, and in this section
we also show that our matrix-valued Jost solutions indeed
generalize the classical Jost solutions, and, moreover, that the
Evans determinant $D$ in \eqref{MainRepr} coincides with the Jost
function (see Theorem \ref{SEMAIN}). We refer to \cite{KS,KS1} for
results relating the Evans function and a Fredholm determinant in
the Schr\"odinger equation case.

In Sections \ref{GMVJS} and \ref{ESAO} we introduce and study the
{\em generalized} matrix-valued Jost solutions, assuming that the
exponential dichotomy for \eqref{unperturbed} admits further (and
finer) exponential splitting (of order two in Section \ref{GMVJS}
and of an arbitrary order in Section \ref{ESAO}). The generalized
matrix-valued Jost solutions are defined abstractly as solutions
of \eqref{perteqn} that approximate decaying
reference solutions of \eqref{unperturbed} at $\pm\infty$ at an exponential
rate which is better than the rate of decay of the reference solutions
themselves (cf.\ Definitions \ref{jostsplit} and \ref{jostsplitn}).
Unlike the matrix-valued Jost solutions from Section
\ref{secjost}, the generalized matrix-valued Jost solutions are
not unique due to the possible appearance of ``lower-order modes'', 
that is, solutions of the perturbed equation \eqref{perteqn} having
worse exponential rate of decay than the corresponding reference
solutions. However, they are unique up to terms that decay
exponentially faster than the respective reference solutions
(see Corollary \ref{Cor6.6} and Theorem \ref{+splitn}). Nevertheless, we
prove in Lemma \ref{mix} that the Evans determinant $D$ defined by
means of the generalized matrix-valued Jost solutions is
determined uniquely. The existence of the generalized matrix-valued
Jost solutions is proved using a ``mixed'' system of integral
equations of Volterra- and Fredholm-type in Theorems \ref{+split}
and \ref{+splitn}. The use of the finer exponential splitting
allows us to prove the existence of the generalized matrix-valued
Jost solutions and formula \eqref{MainRepr} under much weaker assumptions
on the exponential decay of the perturbation than those
in Section \ref{secjost}: In Sections \ref{GMVJS} and
\ref{ESAO} the decay of the perturbation is controlled by the width of {\em
each} of the Bohl segments corresponding to the fine exponential splitting
(cf.\ Theorems \ref{+split} and \ref{posl}). This constitutes one of the
main new  effects observed in this paper: Passing from the case of two- to
the case of many Bohl segments in the unperturbed equation
\eqref{unperturbed} requires replacing the system of two Volterra integral
equations by a cascade-type mixed system of Volterra- and 
Fredholm-type integral equations so that one can take care of the
appearance of the above mentioned ``lower-order modes''. As a
generalization in yet another direction, we express in Corollary
\ref{RATIO} the infinite-dimensional determinant $\det_{2}(I+K)$ in
\eqref{MainRepr} as a ratio of two Evans determinants corresponding to the
perturbed and the unperturbed equation \eqref{unperturbed} thus making the
connection with the interpretation of $\det_{2}(I+K)$ as a (modified)
perturbation determinant (cf.\ \cite[Sect.\ IV.3]{GK69}).

In Section \ref{APERTE} we treat the special but important case of the
constant coefficient unperturbed equation \eqref{unperturbed}, $A(x)=A$,
$x\in\bbR$. We introduce the generalized matrix-valued Jost solutions as
solutions of the ``mixed'' system of Volterra- and Fredholm-type
integral equations in Definition \ref{AUTjsp}. Using the
Jordan normal form of $A$, the existence of the solutions and formula
\eqref{MainRepr} are proved under the following mild assumptions in Theorem
\ref{MAINAUT}: Either we suppose $\|R\|_{\bbC^{d\times d}}\in
L^1(\bbR)$ if the spectrum of $A$ is semi-simple, or we
assume a polynomial decay of $\|R\|_{\bbC^{d\times d}}$. These
hypotheses are shown to be optimal.

Finally, in Section \ref{EVFUN} we show that the Evans determinant
$D$ equals the Evans function for
\eqref{unperturbed} and \eqref{perteqn} defined via the generalized
matrix-valued Jost solutions.

{\bf A glossary of notation:}\, For $d,d_1,d_2\in\bbN$, let $\bbC^{d_1\times
d_2}$ be the set of $d_1\times d_2$ matrices with complex entries; by
${\tr}_{\bbC^d}(M)$ and ${\det}_{\bbC^d} (M)$ we denote the trace and
determinant of a $d\times d$ matrix $M\in \bbC^{d\times d}$.
$\|\cdot\|_{\bbC^{d}}$ denotes a vector norm in $\bbC^d$;
$\|\cdot\|_{\bbC^{d\times d}}$ a matrix norm in $\bbC^{d\times d}$;
$\|\cdot\|_{\mathcal{X}}$ a norm in a Banach space $\mathcal{X}$. For
$p\geq 1$ and $J\subseteq\bbR$, $L^p(J)$,
$L^p(J)^d$, and $L^p(J)^{d\times d}$ are the
usual Lebesgue spaces on $J$ with values in
$\bbC$, $\bbC^d$, and $\bbC^{d\times d}$, associated with Lebesgue measure
$dx$ on $J$. Similarly,
$L^p(J;w(x)dx)$, $L^p(J;w(x)dx)^d$ and
$L^p(J;w(x)dx)^{d\times d}$ are the weighted spaces
with a weight $w\geq 0$; the spaces of {\em
bounded continuous} functions on $J$ are denoted by
$C_{\rm b}(J)$, $C_{\rm b}(J)^d$ and
$C_{\rm b}(J)^{d\times d}$,
respectively. The identity matrix in $\bbC^{d\times d}$ is denoted by
$I_d$ and the identity operator on a Banach space $\cX$ is denoted by
$I$ (or by $I_{\cX}$ if its dependence on $\cX$ needs to be stressed). For
a projection $Q$ we write $\im (Q)$ and
$\ker(Q)$ to denote the range and the null space (kernel) of $Q$. If $T$ is
a linear operator on a Banach space $\mathcal{X}$, then $T^{-1}\in\cB(\cX)$
denotes the (bounded) inverse operator of $T$, with $\cB(\cX)$ being the
Banach space of bounded linear operators on $\cX$. Moreover,
$\spec(T)$ denotes the spectrum of $T$, and
$T|_{\mathcal{Y}}$ denotes the restriction of $T$ to a subspace
$\mathcal{Y}$ of $\mathcal{X}$. If $\mathcal{X}_1$ and
$\mathcal{X}_2$ are two subspaces of $\mathcal{X}$, then
$\mathcal{X}_1\dot{+}\mathcal{X}_2$ denotes their direct (but not
necessarily orthogonal) sum. $\cB_p(\cH)$ denotes the Schatten--von
Neumann ideals of compact operators on a Hilbert space $\cH$ with
singular values in $\ell^p$, $p\geq 1$. A generic constant is denoted by
$c$ or $C$, and we also use the abbreviations $\bbR_+=[0,\infty)$,
$\bbR_-=(-\infty,0]$.

\section{Preliminaries}\label{NAP}

Throughout this paper we make the following assumptions on $A$ and $R$ in
\eqref{unperturbed} and \eqref{perteqn}: 
\begin{hypothesis} \lb{h2.1}
Suppose
\begin{equation}
A, R \in L^1_{\loc}(\bbR)^{d\times d}.
\end{equation}
\end{hypothesis}

Let $\Phi$ denote the fundamental matrix solution of the unperturbed
equation \eqref{unperturbed} so that $\Phi'(x)=A(x)\Phi(x)$, $x\in\bbR$
and $\Phi(0)=I_d$. Throughout, we assume that the propagator
$\Phi(x)\Phi(x')^{-1}$ of \eqref{unperturbed} is {\em
exponentially bounded} on $\bbR$, that is, we suppose that
\begin{equation}
\sup_{x,x'\in\bbR, \, |x-x'|\le
1}\|\Phi(x)\Phi(x')^{-1}\|_{\bbC^{d\times d}} < \infty,
\end{equation}
or equivalently, there exist constants $C \in [1,\infty)$ and $\alpha \in
\bbR$, such that
\begin{equation}
\|\Phi(x)\Phi(x')^{-1}\|_{\bbC^{d\times d}} \leq C e^{\alpha |x-x'|}, 
\quad x, x' \in \bbR.
\end{equation}
First, we recall basic definitions and some standard facts
regarding Bohl and Lyapunov exponents and dichotomies (see, e.g.,
\cite[Sec.1]{Cop}, \cite[Ch.\ III]{DK}).

A projection $Q$ in $\bbC^d$ is said to be {\em
uniformly conjugated} by $\Phi$ on $\bbR$ if 
\begin{equation}
\sup_{x\in\bbR}
\|\Phi(x)Q\Phi(x)^{-1}\|_{\bbC^{d\times d}} < \infty.
\end{equation}
The geometrical
meaning of this condition is that the angle between the range and the
null space of the projection $\Phi(x)Q\Phi(x)^{-1}$ in $\bbC^d$ is
uniformly separated from zero for all $x\in\bbR$ (see, e.g.,
\cite[Corollary IV.1.1]{DK}).

Let $Q$ be a projection in $\bbC^d$ which is uniformly conjugated
by $\Phi$ on $\bbR$. The {\em upper Bohl exponent on
$\bbR$} associated with the projection $Q$, denoted by
$\varkappa(Q)$, is defined as the infimum of all 
$\varkappa\in\bbR$, such that for some constant $C(\varkappa)\in
[1,\infty)$, the following inequality holds for all $x,x'\in\bbR$ such that
$x\ge x'$,
\begin{equation}\label{DEFVAR}
\|\Phi(x)Q\Phi(x')^{-1}\|_{\bbC^{d\times d}}\le
C(\varkappa)e^{\varkappa(x-x')}.
\end{equation}

The {\em lower Bohl exponent on $\bbR$}, denoted by
$\varkappa'(Q)$, is defined as the supremum of all $\varkappa\in\bbR$, such
that \eqref{DEFVAR} holds for all $x\le x'$, where $x,x'\in\bbR$.

Similarly, given a projection $Q$ which is uniformly conjugated by
$\Phi$ on $\bbR_+$ (resp., on $\bbR_-$), the upper
Bohl exponent, $\varkappa_+(Q)$, and the lower Bohl exponent,
$\varkappa'_+(Q)$, on $\bbR_+$ associated with the projection $Q$
(resp., the Bohl exponents $\varkappa_-(Q)$ and
$\varkappa'_-(Q)$ on $\bbR_-$), are defined in the same way except that
one takes $x,x'\in\bbR_+$ (respectively,
$x,x'\in\bbR_-$) in \eqref{DEFVAR} . The upper and lower Bohl exponents on
$\bbR_+$ can also be expressed as follows (see 
\cite[Theorem III.4.4]{DK}) 
\begin{align}
\varkappa_+(Q)&=\limsup_{(x-x')\to \infty, \, x'\to  \infty}
\frac{\log \|\Phi(x)Q\Phi(x')^{-1}\|_{\bbC^{d\times d}}}{x-x'},\\
\varkappa'_+(Q)&=-\limsup_{(x-x')\to \infty, \, x'\to 
\infty} \frac{\log \|\Phi(x')Q\Phi(x)^{-1}\|_{\bbC^{d\times d}}}{x-x'};
\end{align}
similar formulas hold for $\varkappa_-(Q)$ and $\varkappa'_-(Q)$.

We introduce  the upper, $\lambda_+(Q)$, and the lower,
$\lambda'_+(Q)$, {\em Lyapunov exponents on $\bbR_+$} associated
with the projection $Q$ by the formulas
\begin{equation}
\lambda_+(Q)=\limsup_{x\to  \infty} \frac{\log
\|\Phi(x)Q\|_{\bbC^{d\times d}}}{x},\quad \lambda'_+(Q)=-\limsup_{x\to 
\infty}
\frac{\log \|Q\Phi(x)^{-1}\|_{\bbC^{d\times d}}}{x},
\end{equation}
and the Lyapunov exponents on $\bbR_-$ by the formulas
\begin{equation}
\lambda_-(Q)=\limsup_{x\to  -\infty} \frac{\log
\|Q\Phi(x)^{-1}\|_{\bbC^{d\times d}}}{|x|}, \quad
\lambda'_-(Q)=-\limsup_{x\to  -\infty} \frac{\log
\|\Phi(x)Q\|_{\bbC^{d\times d}}}{|x|}.
\end{equation}

We remark that the upper (lower) Bohl exponent measures the best (worst)
exponential growth of the propagator $\Phi(x)\Phi(x')^{-1}$ and
the upper (lower) Lyapunov  exponent measures the best (worst)
exponential growth of the fundamental solution $\Phi$ relative to the
projection $Q$. Under the assumption that the propagator of
\eqref{unperturbed} is exponentially bounded, the Bohl exponents are finite
\cite[Theorem III.4.2]{DK}. We note the inequalities
 \begin{equation}\label{INEQbohls}
 \varkappa'(Q)\le\varkappa_{\pm}'(Q)\le \lambda_{\pm}'(Q)\le
 \lambda_{\pm}(Q)\le\varkappa_{\pm}(Q)\le\varkappa(Q),\end{equation}
  and stress that they are {\em strict} in general (in particular, the
possible inequality of the Bohl and Lyapunov exponents is exhibited in
 the classical Perron example, see, e.g., \cite[Sect.\ III.4.4]{DK}).

\begin{remark}\label{RESC}
 The set of Bohl and Lyapunov exponents for the {\em rescaled
 equation},
 \begin{equation}\label{RESCeqn}
 y'(x)=(A(x)-\mu I_d) y(x),\quad x\in\bbR, \; \mu\in\bbC, 
 \end{equation}
with the associated propagator
 $e^{-\mu(x-x')}\Phi(x)\Phi(x')^{-1}$, $x,x'\in\bbR$, is obtained from the
set of the exponents of equation
 \eqref{unperturbed} by shifting the latter by the amount $-\re (\mu)$.
 \hfill$\Diamond$
\end{remark}

\begin{example}\label{CC} If the unperturbed
equation \eqref{unperturbed} is autonomous then the Bohl and Lyapunov
exponents are the same and equal to the real parts of the eigenvalues of
$A$. To give more details and fix notations, suppose that
$A(x)=A$, $x\in\bbR$, for some $A\in\bbC^{d\times d}$, and let
$\nu_k\in\bbC$, $k=1,\dots,d$, be the eigenvalues of $A$. We
split the set of the eigenvalues as
$\spec(A)=\cup_{j=1}^{d'}\Sigma_j$, where $1\le d'\le d$ and for
each $j$ all eigenvalues that belong to $\Sigma_j$ have the same
real part denoted by $\varkappa_j$, $j=1,\dots,d'$. Let $Q_j$
denote the spectral projection for $A$ such that
$\spec(A|_{\im(Q_j)})=\Sigma_j$. Then the Bohl and the Lyapunov exponents on
$\bbR$ (and on $\bbR_\pm$) associated with $Q_j$, are equal to
$\varkappa_j$ (we note that we will always enumerate $\varkappa_j$ in
increasing order of magnitude, $\varkappa_j<\varkappa_{j'}$ for $j<j'$,
$j,j'=1,\dots,d'$).
\hfill$\Diamond$
\end{example}

Returning to the general nonautonomous case, equation
\eqref{unperturbed} is said to have a {\em bounded} dichotomy $Q$ on $\bbR$
if $Q$ is a projection in $\bbC^d$ (called the {\em dichotomy
projection}) so that for some constant $M\in [1,\infty)$, 
 the following inequalities hold for all $x,x'\in\bbR$:
\begin{equation}  \label{bouddich}
\begin{split}\| \Phi (x)Q\Phi (x')^{-1}\|_{\bbC^{d\times d}} & \leq
M, \quad x\geq x', \\
\| \Phi (x)(I_d-Q)\Phi (x')^{-1}\|_{\bbC^{d\times d}} &\le M, \quad x\leq
x'.
\end{split}
\end{equation}
Equation \eqref{unperturbed} is said to have an {\em exponential}
dichotomy $Q$ on $\bbR$ if $Q$ is a projection in $\bbC^d$
 so that for some positive constants $\varkappa$ and $\varkappa'$, and
 some constants $C(\varkappa),C(\varkappa')\in [1,\infty)$,
 the following inequalities hold for all $x,x'\in\bbR$:
\begin{eqnarray}\| \Phi (x)Q\Phi (x')^{-1}\|_{\bbC^{d\times d}} &
\leq&
C(\varkappa)e^{-\varkappa(x-x')}, \quad x\geq x', \label{expdich1}\\
\| \Phi (x)(I_d-Q)\Phi (x')^{-1}\|_{\bbC^{d\times d}} &\le&
C(\varkappa')e^{\varkappa'(x-x')}, \quad x\leq x'.\label{expdich2}
\end{eqnarray}
We note that $Q$ is an exponential dichotomy for
\eqref{unperturbed} on $\bbR$ if and only if the following
inequalities for the Bohl exponents hold:
 \begin{equation}\label{EDBE}
 \varkappa(Q)<0<\varkappa'(I_d-Q).
 \end{equation}
Dichotomies on $\bbR_+$, respectively on $\bbR_-$, are defined in
the same way except in \eqref{bouddich} or
\eqref{expdich1}, \eqref{expdich2} one takes only
$x,x'\in\bbR_+$, respectively, $x,x'\in\bbR_-$.

\begin{example}\label{CC1} Assume in Example \ref{CC} that
$\spec(A)\cap i\,\bbR=\emptyset$ and $d'\ge 2$. Consider
$k_0\in\{1,\dots,d'-1\}$ such
that $\varkappa_{k_0}<0<\varkappa_{k_0+1}$. Then $Q=\sum_{j=1}^{k_0} Q_j$
is the exponential dichotomy projection for the equation $y'(x)=Ay(x)$ on
$\bbR$, $\bbR_+$, and
$\bbR_-$. \hfill$\Diamond$
\end{example}

\begin{definition}\label{DEFexpsplit} (\cite[Sect.\ IV.4]{DK}.)
Assume $d\ge 2$, and let
$1< d'\le d$. A  system $\{Q_j\}_{j=1}^{d'}$ of disjoint
projections in $\bbC^d$ is called an {\em exponential splitting}
for \eqref{unperturbed} on $\bbR$ if the following conditions hold: \\
$(i)$ $\sum_{j=1}^{d'}Q_j=I_d$. \\
$(ii)$ Each projection $Q_j$ is uniformly
conjugated by $\Phi$ on $\bbR$. \\
$(iii)$ The segments $[\varkappa'(Q_j), \varkappa(Q_j)]$ are disjoint.
\end{definition}

 The segments $[\varkappa'(Q_j), \varkappa(Q_j)]$ of
the real axis are called the {\em Bohl segments} associated with
$Q_j$; they are determined by the lower and upper Bohl exponents
on $\bbR$ associated with the projections $Q_j$, and in what
follows they will always be numbered so that
$\varkappa(Q_j)<\varkappa'(Q_{j+1})$, $j=1,\dots,d'-1$. If $d'\ge
2$ and \eqref{unperturbed} has an exponential dichotomy $Q$ on
$\bbR$, then for some $k_0$, $1\le k_0<d'-1$, we have the
following splitting:
\begin{equation}
 Q=\sum_{j=1}^{k_0}Q_j \, \text{ and } \,
I_d-Q=\sum_{j=k_0+1}^{d'}Q_j,
\end{equation}
so that
$\varkappa(Q_{k_0})<0<\varkappa'(Q_{k_0+1})$.

Replacing $\bbR$ by $\bbR_+$ or $\bbR_-$ in Definition
\ref{DEFexpsplit}, one can also consider exponential splittings on
$\bbR_+$ or $\bbR_-$. We note that the exponential dichotomy projection $Q$
on $\bbR$ is uniquely defined, while exponential dichotomy
projections on $\bbR_+$ or $\bbR_-$ are {\em not}. Indeed, for example, in
the case of $\bbR_+$,  the subspace $N_s$ of initial data of
bounded solutions of \eqref{unperturbed} on $\bbR_+$ is
determined uniquely; to obtain an exponential dichotomy projection
$Q$ on $\bbR_+$ one can choose an arbitrary direct complement $N_u$
of the subspace $N_s$ in $\mathbb{C}^d$ so that 
$\mathbb{C}^d=N_s\dot{+}N_u$, and define the exponential dichotomy
projection $Q$ such that $\ran(Q)=N_s$ and $\ker(Q)=N_u$ 
(see \cite[Remark IV.3.4]{DK}). We recall that if $Q$ is an exponential
dichotomy for \eqref{unperturbed} on $\bbR$, then $Q$ is
also  an exponential dichotomy on both semi-axes $\bbR_+$ and
$\bbR_-$. Similarly, an exponential splitting $\{Q_j\}_{j=1}^{d'}$
on $\bbR$ induces exponential splittings on $\bbR_+$ and $\bbR_-$;
we stress that in view of \eqref{INEQbohls}, the Bohl segments for
$\bbR$ are, generally, wider than the segments for $\bbR_+$ and
$\bbR_-$. This may happen since an exponential dichotomy on $\bbR_+$ does
not necessarily imply an exponential dichotomy on $\bbR$.

\begin{example}\label{CC2} In Examples \eqref{CC} and \eqref{CC1}
each Bohl segment $[\varkappa'(Q_j),\varkappa(Q_j)]$ associated
with $Q_j$ degenerates into a single point $\varkappa_j$,
$j=1,\dots,d'$, while the Bohl segments associated with the
(exponential dichotomy) splitting $\{Q,I_d-Q\}$ are given by
$[\varkappa_1,\varkappa_{k_0}]$ and
$[\varkappa_{k_0+1},\varkappa_{d'}]$; we recall that in this case 
\begin{align}
&\varkappa_1=\inf\{\re(\nu)\in\bbR \,|\,\nu\in\spec(A)\}, \quad 
\varkappa_{k_0}=\sup\{\re(\nu)\in\bbR \,|\, \nu\in\spec(A), \,
\re(\nu)<0\}, \no \\
& \varkappa_{k_0+1}=\inf\{\re(\nu)\in\bbR \,|\,
\nu\in\spec(A), \, \re(\nu)>0\}, \quad 
\varkappa_{d'}=\sup\{\re(\nu)\in\bbR \,|\, \nu\in\spec(A)\}.
\end{align}
\hfill$\Diamond$
\end{example}

\begin{remark}\label{VVSOLS} The Lyapunov exponents for nonzero
$\bbC^d$-valued solutions $y=y(x)$ of \eqref{unperturbed} are defined as
follows:
\begin{equation}
\lambda_\pm(y)=\limsup_{x\to\pm\infty}
\frac{\log(\|y(x)\|_{\bbC^{d}})}{x},\quad
\lambda'_\pm(y)=\liminf_{x\to\pm\infty}
\frac{\log(\|y(x)\|_{\bbC^{d}})}{x}.
\end{equation}
For any nonzero $\bbC^d$-valued solution $y$ of \eqref{unperturbed}, and
for any given exponential splitting $\{Q_k\}_{k=1}^{d'}$, there exists a
$j \in \{1,\dots,d'\}$ such that all four Lyapunov exponents
$\lambda_{\pm}(y), $
$\lambda_{\pm}'(y)$ belong to the $j$th Bohl segment. One verifies
that the actual $j$th Bohl segment depends on the initial data
$y_0=y(0)\neq 0$ of the solution
$y$ of \eqref{unperturbed}, say, on $\bbR_+$, as follows: 
\begin{align}
& \text{If } \, y_0\in\im\bigg(\sum_{k=1}^{j}Q_k\bigg) \, \text{ then } \, 
\lambda_+(y)\le\lambda_+(Q_j), \no \\ 
& \text{if } \, y_0\in\im\bigg(\sum_{k=1}^{j}Q_k\bigg)
\bigg\backslash\im\bigg(\sum_{\ell=1}^{j-1}Q_{\ell}\big) \, \text{ then }
\, 
\lambda'_+(y)\ge\lambda'_+(Q_j), \no \\
& \text{if } \, \lambda_+(y)<\lambda'_+(Q_j) \, \text{ then } \,
y_0\in\im\bigg(\sum_{k=1}^{j-1}Q_k\bigg),  \label{IMPINT} \\
& \text{if } \, \lambda'_+(y)>\lambda_+(Q_j) \, \text{ then } \, 
y(0)\in\bbC^d\big\backslash\im\bigg(\sum_{k=1}^{j}Q_k\bigg), \no 
\end{align}
where $j\in\{1,\dots,d'\}$ and we set $Q_0=0$. Similar assertions hold for
the Laypunov exponents on $\bbR_-$. \hfill$\Diamond$
\end{remark}

For future references, we record the assumptions on the coefficient
$A$ of the unperturbed equation \eqref{unperturbed}.

\begin{hypothesis}\label{diffeq} Assume that
$ A \in L^1_{\loc}(\bbR)^{d\times d}$. In addition, we suppose that
the propagator
$\Phi (x)\Phi (x')^{-1}$, $x,x'\in\bbR$, is exponentially bounded on $\bbR$,
that is,
\begin{equation}
\sup_{x,x'\in\bbR, \, |x-x'|\le 1} \|\Phi(x)\Phi(x')^{-1}\|_{\bbC^{d\times
d}}<\infty,
\end{equation}
and that \eqref{unperturbed} has an exponential dichotomy $Q$ on $\bbR$.
\end{hypothesis}

Turning to the perturbed equation \eqref{perteqn}, we assume that
the perturbation $R$ in \eqref{perteqn} satisfies the condition
\begin{equation}\label{L1R}
\|R\|_{\bbC^{d\times d}}\in L^1(\bbR).
\end{equation}
Let $U$  and $|R|$ denote the $d\times d$ matrices in the
polar decomposition of $R$:
\begin{equation}
R(x)=U(x)|R(x)|, \quad |R(x)|=(R(x)^*R(x))^{1/2},\quad x\in\bbR.
\end{equation}
Throughout this paper we will use the notation
\begin{equation}\label{FACTOR}
R_\ell(x)=U(x)|R(x)|^{\frac{1}{2}}, \quad
R_r(x)=|R(x)|^{\frac{1}{2}}, \quad x\in\bbR,
\end{equation}
so that
\begin{equation}
R(x)=R_\ell(x)R_r(x), \quad x\in \bbR.
\end{equation}

\begin{lemma}\label{tech2}
 Assume Hypothesis \ref{diffeq} and condition \eqref{L1R}.
Then the integral operator  $K$ with integral kernel given  by
\begin{equation}\label{KK}
K(x,x')=\begin{cases}
-R_r(x) \Phi(x)Q\Phi(x')^{-1} R_\ell(x'),& x\ge x',\\
R_r(x)\Phi(x)(I_d-Q)\Phi(x')^{-1} R_\ell(x'),&x< x',
\end{cases}  \quad x,x'\in\bbR,
\end{equation}
is a Hilbert--Schmidt operator on $L^2(\bbR)^d$. Moreover, the map
\begin{equation} \label{contt}
\begin{cases} L^1(\bbR)^{d\times d} \to \cB_2(L^2(\bbR)^d) \\
\hspace*{1.17cm} R\mapsto K  \end{cases}
\end{equation}
is continuous.
\end{lemma}
\begin{proof} Since by hypothesis  equation \eqref{unperturbed} has an
 exponential dichotomy $Q$ on $\bbR$, the projection $Q$ is also a bounded
dichotomy on $\bbR$, and thus \eqref{bouddich} holds.
 Using \eqref{bouddich},  the
 Hilbert--Schmidt norm of the integral operator
  $K$ can be estimated as follows:
\begin{align}
\|K\|_{{\cB}_2(L^2(\bbR)^d)}&=\int\int_{\bbR\times
\bbR} dx\, dx'\, \|K(x,x')\|_{\bbC^{d\times d}}^2   \no \\
&=\int_\bbR dx \int_{-\infty}^x
dx'\|R_r(x)\Phi(x)Q\Phi(x')^{-1}R_\ell(x')\|_{\bbC^{d\times d}}^2
\no \\
& \quad +\int_\bbR dx\int^{\infty}_x
dx'\|R_r(x)\Phi(x)(I_d-Q)\Phi(x')^{-1}R_\ell(x')\|_{\bbC^{d\times d}}^2
\no \\
&\le M^2\int_\bbR dx\int_\bbR dx' \|R(x)\|_{\bbC^{d\times d}}
\|R(x')\|_{\bbC^{d\times d}} \no \\
& = M^2\bigg (\int_\bbR dx \, \|R(x)\|_{\bbC^{d\times d}}\bigg )^2.
\end{align}
In the last estimate we used the inequalities
\begin{equation} \label{15.4}
\| R_r(x)\|_{\bbC^{d\times d}}
\leq \| R(x)\|_{\bbC^{d\times d}}^{\frac{1}{2}}, \quad
\| R_\ell (x)\|_{\bbC^{d\times d}}
\leq \|R(x)\|_{\bbC^{d\times d}}^{\frac{1}{2}}, \quad x\in\bbR,
\end{equation}
 which follow from the spectral theorem for self-adjoint matrices.
The continuity of the mapping in \eqref{contt} holds because the
mappings $R\mapsto R_r$ and $R\mapsto
R_\ell$ are continuous from the Banach space
$L^1(\bbR)^{d\times d}$ to the Hilbert space $L^2(\bbR)^{d\times
d}$.
\end{proof}

\begin{remark}\label{invert}
Since the propagator of \eqref{unperturbed} is exponentially
bounded, the formula
\begin{equation}\label{semi}
(E^t_\Phi u)(x)=\Phi(x-t)u(x-t),\quad u\in L^2(\bbR)^d, \quad x\in
\mathbb{R},\; t\geq 0,
\end{equation}
defines a strongly continuous semigroup, $\{E^t_\Phi\}_{t\geq 0}$, on
$L^2(\bbR)^d$, called the {\em evolution semigroup} (see \cite[Sect.\
3.2]{CL} for a detailed discussion and extensive bibliography). It is
well-known (see, e.g., \cite[Theorem\ 3.17]{CL}) that the generator, $G$,
of the semigroup in \eqref{semi} has a bounded inverse if and only if
\eqref{unperturbed} has an exponential dichotomy $Q$ on $\bbR$ and, if
this is the case, the inverse operator
$G^{-1}$ is an integral operator on $L^2(\bbR)^d$ with
integral kernel given by
\begin{equation}\label{KERNF}
F(x,x')=\begin{cases}
- \Phi(x)Q\Phi(x')^{-1} ,& x\ge x',\\
\Phi(x)(I_d-Q)\Phi(x')^{-1} ,&x< x',
\end{cases}  \quad x, x' \in\bbR,
\end{equation}
where $Q$ is the dichotomy projection. Thus, under the  hypotheses
of  Lemma \ref{tech2}, the integral operator $K$ admits the
representation $K=\ol{{\cM}_{R_r} G^{-1} {\cM}_{R_\ell}}$, where
${\cM}_{R_r}$ and ${\cM}_{R_\ell}$ denote the maximally defined operators
of multiplication by the matrix-valued functions $R_r, R_\ell\in
L^2(\bbR)^{d\times d}$, and the operator
$S={\cM}_{R_r} G^{-1} {\cM}_{R_\ell}$ (with the bounded closure $\ol{S}$) 
is defined on the set of functions $u\in\dom({\cM}_{R_\ell})$ such that 
$Su\in L^2(\mathbb{R})^{d\times d}$. We refer to \cite{GLMZ05}
for a detailed discussion of Birman--Schwinger-type operators and the
Birman--Schwinger principle.  
\hfill$\Diamond$
\end{remark}

\begin{remark} \label{invertNEW}
Using the semigroup in \eqref{semi}
one can describe the exponential dichotomy projection $Q$ for
\eqref{unperturbed} on $\bbR$ by means of the Riesz spectral projection
for the operator $E_\Phi^1$, with $t=1$ in \eqref{semi}, as follows
(see \cite[Theorems\ 3.13, 3.17]{CL}): Assume that
the generator $G$ of the semigroup in \eqref{semi}
is invertible on $L^2(\bbR)^d$. By \cite[Theorem\ 3.13]{CL},
it follows that $E_\Phi^1$ has no spectrum on the unit circle.
Let $\mathcal{Q}$ denote the spectral projection
for $E_\Phi^1$ on $L^2(\bbR)^d$ such that
$\spec(E_\Phi^1|_{\im(\mathcal{Q})})$
is the part of the spectrum of $E_\Phi^1$ contained inside the unit disc.
Theorem 3.17 in \cite{CL} states that $\mathcal{Q}$ is an operator
of multiplication in $L^2(\bbR)^d$ by a projection-valued
function $Q(\cdot)\in C_{\rm b}(\bbR)^{d\times d}$
and, in addition, $Q=Q(0)$ is the exponential dichotomy projection for
\eqref{unperturbed} on $\bbR$.
\hfill$\Diamond$
\end{remark}

\begin{remark} \label{look}
The proof of Lemma \ref{tech2} shows that the
assumption of an {\em exponential} dichotomy on $\bbR$ for
\eqref{unperturbed} can be relaxed to require only a {\em bounded}
dichotomy. Under this weaker assumption the generator $G$ of the
semigroup in \eqref{semi} may not be invertible (cf.\ Remark
\ref{invert}), or may not even be a Fredholm operator on
$L^2(\bbR)^d$. Indeed, by a well-known Dichotomy Theorem
(sometimes called Palmer's theorem), the operator $G$  is Fredholm
on $L^2(\bbR)^d$ if and only if \eqref{unperturbed} has  {\it exponential}
dichotomies $Q_-$  on $\bbR_-$ and $Q_+$ on $\bbR_+$; see \cite{BAG},
\cite{Pa}, \cite{Palm}, \cite{Sack} or \cite[Theorem\ 3.2]{S}, and also
\cite{LT}, \cite{SS2} for more recent versions of the dichotomy theorem.   
\hfill$\Diamond$
\end{remark}

Next, we will discuss the Bohl exponents and exponential
splittings for the perturbed equation \eqref{perteqn}. Since the propagator
of \eqref{unperturbed} is exponentially
bounded, and condition \eqref{L1R} holds, it follows from the
variation of constants formula and Gronwall's inequality that the
propagator of the perturbed equation \eqref{perteqn} is also
exponentially bounded (see, e.g., \cite[Lemma\ IV.4.1]{H}). Thus, the
Bohl exponents for the perturbed equation \eqref{perteqn} are finite.

We conclude this section with two technical results to the effect that
first, the exponential dichotomy and exponential splitting of
\eqref{unperturbed} on $\bbR_+$ persist under
$L^1(\bbR_+)^{d\times d}$-perturbations, and second, that the
corresponding Bohl exponents do not change under these
perturbations (similar facts hold for $\bbR_-$). A proof of the
first part of these results can be found in \cite[Proposition 4.1]{Cop}
or \cite[Theorem IV.5.1]{DK}. However, we were not able to find the second
part in the literature, and so we will briefly sketch a proof in Appendix
\ref{AppendB} by modifying some of the arguments in \cite[Ch.\ IV]{DK}.
Lemma \ref{stub_dich} is used in Section \ref{GMVJS}, while
its direct generalization, Lemma \ref{stab_split}, is used in Section 
\ref{ESAO}.

\begin{lemma}\label{stub_dich} Assume that the propagator
of the unperturbed equation \eqref{unperturbed} on $\bbR_+$ is
exponentially bounded, that the unperturbed equation \eqref{unperturbed}
has an exponential dichotomy $Q$ on $\bbR_+$, and the perturbation $R$
satisfies $\|R\|_{\bbC^{d\times d}} \in L^1(\bbR_+)$. Then
the perturbed equation \eqref{perteqn} has an exponential dichotomy $P$ on
$\bbR_+$ such that
\begin{equation}\label{bohl12}
\varkappa_+(P)=\varkappa_+(Q) \, \text{ and } \,
\varkappa'_+(I_d-P)=\varkappa'_+(I_d-Q),
\end{equation}
where
$\varkappa_+(P)$ and $\varkappa'_+(I_d-P)$ are the Bohl exponents
for the perturbed equation \eqref{perteqn} on $\bbR_+$. Here, $P$ is the
projection in $\bbC^d$ parallel to $\ker (Q)$ onto the subspace $N$
consisting of the values $y(0)$ at zero
 of all {\em bounded} solutions $y$ of \eqref{perteqn} on $\bbR_+$.
In addition,
\begin{equation}
\dim(\im (P))=\dim(\im (Q))\, \text{ and } \, \dim(\ker (P))=\dim(\ker(Q)).
\end{equation}
\end{lemma}

\begin{lemma}\label{stab_split}
Assume that the propagator of the unperturbed equation \eqref{unperturbed}
on $\bbR_+$ is exponentially bounded, the unperturbed equation
\eqref{unperturbed} has an exponential
splitting $\{Q_j\}_{j=1}^{d'}$, $1<d'\le d$, on $\bbR_+$, and the
perturbation $R$ satisfies $\|R\|_{\bbC^{d\times d}}\in
L^1(\bbR_+)$. Then, the perturbed equation \eqref{perteqn} has an
exponential splitting $\{P_j\}_{j=1}^{d'}$ on $\bbR_+$ such that
\begin{equation}
\dim(\im (P_j))=\dim(\im (Q_j)), \quad j=1,\dots,d',
\end{equation}
and
\begin{equation}\label{bohl12split}
\varkappa_+(P_j)=\varkappa_+(Q_j), \quad j=1,\dots,d',
\end{equation}
where $\varkappa_+(P_j)$ are the Bohl exponents for the perturbed
equation \eqref{perteqn} on $\bbR_+$.
\end{lemma}

\section{Matrix-Valued Jost Solutions}\label{secjost}

In this section we discuss optimal conditions on the perturbation $R$
under which, assuming an exponential dichotomy of
\eqref{unperturbed}, one can establish existence and uniqueness of
bounded matrix-valued solutions of \eqref{perteqn} (we will call them {\em
matrix-valued Jost solutions}) that are asymptotically close 
with respect to an exponential weight factor to bounded solutions of the
unperturbed equation \eqref{unperturbed}. These results are motivated by
the study of the one-dimensional Schr\"odinger equation on the real axis,
\begin{equation}\label{schr}
-u''(x)+V(x)u(x)=k^2u(x),\quad x\in\bbR,\; k\in\bbC,
\end{equation}
with an integrable potential $V\in L^1(\bbR;dx)$, where the Jost solutions,
$u_\pm(k,\cdot)$, are introduced as the solutions of \eqref{schr}
that are asymptotically close to the free plane waves:
$u_{\pm}(k,x)\underset{x\to \pm\infty}{\sim}
e^{\pm ikx}$ (see, e.g., \cite[Ch.\ XVII]{CS}).

\begin{definition}\label{jost}
Under Hypothesis \ref{diffeq}, $d \times d$ matrix-valued
solutions $Y_\pm$ of the differential equation
\begin{equation}
Y_\pm(x)=(A(x)+R(x))Y_\pm(x), \quad x\in \bbR_\pm,
\end{equation}
with a locally integrable $d\times d$ matrix-valued function
$R$, are called  {\em matrix-valued Jost solutions} if
\begin{equation}\label{lapun+}
\limsup_{x\to  \infty } \frac{\log\|Y_+(x)-\Phi(x)Q
\|_{\bbC^{d\times d}}}{x}<\varkappa_+'(Q)
\end{equation}
and
\begin{equation}\label{lapun-}
\liminf_{x\to  -\infty } \frac{\log\|Y_-(x)-\Phi(x)(I_d-Q)
\|_{\bbC^{d\times d}}}{x}>\varkappa_-(I_d-Q).
\end{equation}
\end{definition}

To motivate this definition, we first
note that the Lyapunov exponents $\lambda_+(Q)$ and ${\lambda'}_+(Q)$ of
the bounded matrix-valued solution $\Phi \, Q$ of \eqref{unperturbed} on
$\bbR_+$ belong to the Bohl segment $[\varkappa'_+(Q),\varkappa_+(Q)]$. 
If $Y_+$ is a matrix-valued solution of \eqref{perteqn} satisfying
\eqref{lapun+}, then $Y_+$ is bounded and therefore Lemma \ref{stub_dich}
implies that its Lyapunov exponents belong to the same segment. Thus,
the significance of \eqref{lapun+} is that the solution $Y_+$ approximates
the solution $\Phi \, Q$ exponentially better than either of
these solutions decays as $x\to \infty$.

Our first result shows the existence and uniqueness of the
matrix-valued Jost solution $Y_+$ on $\bbR_+$ using a
rather strong assumption on the exponential fall-off of the
perturbation $R$.

\begin{theorem}\label{+}
Assume Hypothesis \ref{diffeq} and the condition
\begin{equation}\label{falloff_1}
\|R\|_{\bbC^{d\times d}} \in L^1(\bbR_+;e^{\beta x }dx)
\end{equation}
for some constant $\beta>0$ satisfying the inequality
\begin{equation}\label{fallooof_1}\beta>\lambda_+(Q)-\varkappa_+'(Q).
\end{equation}
 Then the perturbed differential equation \eqref{perteqn}
has a unique matrix-valued Jost  solution $Y_+(x)$, $x\ge 0$.
Moreover, the map
\begin{equation}\label{cont_1}
\begin{cases}
L^1(\bbR_+;e^{\beta x }dx )^{d\times d} \to\bbC^{d\times d} \\
\hspace*{2.47cm} R \mapsto Y_+(0) \end{cases}
\end{equation}
is continuous.
\end{theorem}
\begin{proof}
We split the proof of Theorem \ref{+} into three steps.

{\it Step 1. Uniqueness.} Suppose  that $Y_{+,1}$ and
$Y_{+,2}$ are two different $d\times d$ matrix-valued
Jost solutions of \eqref{perteqn}, and introduce the nonzero
solution $Y(x)=Y_{+,1}(x)-Y_{+,2}(x)$, $x\ge 0$. Using
\eqref{lapun+}, we infer:
\begin{equation}
 \limsup_{x\to  \infty } \frac{\log\|Y(x)\|_{\bbC^{d\times d}}}{x}
\le\max_{k=1,2}\bigg(
\limsup_{x\to  \infty }
\frac{\log\|Y_{+,k}(x)-\Phi(x)Q\|_{\bbC^{d\times d}}}{x}\bigg)
<\varkappa_+'(Q).
\end{equation}
This shows, in particular, that $Y$ is a bounded solution
of \eqref{perteqn}. Since 
$\|R\|_{\bbC^{d\times d}}\in L^1(\bbR_+)$ by hypothesis,
 Lemma \ref{stub_dich} yields an exponential dichotomy
$P$ for \eqref{perteqn} such that
$\varkappa_+'(P)=\varkappa_+'(Q)$. Since $Y$ is a bounded
solution, $Y(0)\in \im (P)$ and hence we arrive at the contradiction
\begin{equation}
 \limsup_{x\to  \infty }
\frac{\log\|Y(x)\|_{\bbC^{d\times d}}}{x}
 \ge \varkappa_+'(P)=\varkappa'_+(Q),
\end{equation}
proving the uniqueness part of the theorem.

{\it Step 2. Existence.} By hypothesis, \eqref{unperturbed} has an
exponential dichotomy  $Q$ on $\bbR_+$. Therefore, for any
$\varepsilon>0$ there exists a constant $C(\varepsilon)\ge 1$ such
that applying \eqref{DEFVAR} in the definitions of
$\varkappa_+'(Q)$ and  $\varkappa_+'(I_d-Q)$, we obtain the estimate
 \begin{align}
\begin{split}
 & \|\Phi(x)\Phi(x')^{-1}\|_{\bbC^{d\times d}}\le
 \|\Phi(x)Q\Phi(x')^{-1}\|_{\bbC^{d\times d}}
 +\|\Phi(x)(I_d-Q)\Phi(x')^{-1}\|_{\bbC^{d\times d}}\\
 & \quad \le
\big(e^{(\varkappa_+'(Q))(x-x')}+
e^{(\varkappa_+'(I_d-Q))(x-x')} \big)f_\varepsilon (|x-x'|),\quad
 0\le x\le x',
\end{split}
\end{align}
where we denoted
\begin{equation}\label{expscaling}
f_\varepsilon(x)=C(\varepsilon)e^{\varepsilon x}, \quad x\ge 0.
\end{equation}
Since $\varkappa_+'(Q)<\varkappa_+'(I_d-Q)$, we get the estimate
\begin{equation}\label{estPhiPhi}
 \|\Phi(x)\Phi(x')^{-1}\|_{\bbC^{d\times d}} \le 2
e^{\varkappa_+'(Q)(x-x')}f_\varepsilon
(|x-x'|), \quad 0\le x\le x'.
\end{equation}
Next, fix $\varepsilon$ such that
 $0<2\varepsilon< \beta-(\lambda_+(Q)-\varkappa_+'(Q))$,
 introduce the function
\begin{equation}
p(x)=e^{(\lambda_+(Q)-\varkappa_+'(Q))x}f_\varepsilon^2(x)
 \|R(x)\|_{\bbC^{d\times d}},\quad x\ge0,
\end{equation}
and observe that $0 \leq p\in L^1(\bbR_+)$ due to
\eqref{falloff_1} and \eqref{fallooof_1}.

On the space $C_{\rm b}(\bbR_+)^{d\times d}$
we now define a Volterra integral operator, $V$, by
\begin{equation}\label{volt_1} (VZ)(x)=\int_x^\infty dx'\,
 e^{-\lambda_+(Q)(x-x')}
 f_\varepsilon(x)^{-1}f_\varepsilon(x')\Phi(x)\Phi(x')^{-1}R(x')Z(x'),
\end{equation}
and consider the corresponding Volterra integral equation
\begin{equation}
 Z(x)=Z^{(0)}(x)-(VZ)(x),  \quad x\geq 0,  \label{volt} 
\end{equation}
where $ Z^{(0)}$ is defined by
\begin{equation}
 Z^{(0)}(x)=e^{-\lambda_+(Q) (x)}f_\varepsilon(x)^{-1}\Phi(x) Q, \quad
x\ge0.
\end{equation}
Using \eqref{estPhiPhi} and the inequality
$\varkappa_+'(Q)-\lambda_+(Q)\le 0$ (cf.\ \eqref{INEQbohls}), one
derives the following estimate for the integral kernel of the
integral operator $V$:
 \begin{align}
  \|V(&x,x') \|_{\bbC^{d\times d}}\le
2e^{(\varkappa_+'(Q)-\lambda_+(Q))(x-x')}
 f_\varepsilon(x)^{-1}f_\varepsilon(x')
 f_\varepsilon(|x-x'|)\|R(x')\|_{\bbC^{d\times d}}\no \\&
 \le2e^{(\varkappa_+'(Q)-\lambda_+(Q))x}f_\varepsilon(x)^{-1}
 e^{-(\varkappa_+'(Q)-\lambda_+(Q))x'}
f^2_\varepsilon(x') \|R(x')\|_{\bbC^{d\times d}} \nonumber \\
&= 2e^{(\varkappa_+'(Q)-\lambda_+(Q))x}f_\varepsilon(x)^{-1} p(x')
\label{amper_1}\\&
 \le2f_\varepsilon(x)^{-1}p(x'), \quad\quad 0\le x\le
 x'.\label{amper}
 \end{align}
Since $f_\varepsilon(x)\ge 1$, $x\in \bbR_+$, and $p\in
L^1(\bbR_+)$, the estimate \eqref{amper} shows that $V$ is a
compact operator on $C_{\rm b}(\bbR_+)^{d\times d}$ with spectral
radius equal to zero. Since the exponential decay of
$\|\Phi \, Q\|_{\bbC^{d\times d}}$ is controlled by $\lambda_+(Q)$, we
see that $Z^{(0)}\in C_{\rm b}(\bbR_+)^{d\times d}$. Thus, the
integral equation \eqref{volt} has a unique solution $Z\in
C_{\rm b}(\bbR_+)^{d\times d}$ that can be
obtained by the iteration process
\begin{equation}
Z(x)=\sum_{j=1}^\infty Z^{(j)}(x), \quad
Z^{(j)}(x)=(VZ^{(j-1)})(x),\quad j\in\bbN, \; x\ge 0.
\end{equation}
Moreover, from \eqref{volt} and \eqref{amper_1} it also follows
that
\begin{align}
\begin{split}
& f_\varepsilon(x)\|Z(x)-Z^{(0)}(x)\|_{\bbC^{d\times d}}   \label{zzz_1} \\
& \quad \le
2 \|Z\|_{C_{\rm b}(\bbR_+)^{d\times d}}
e^{(\varkappa_+'(Q)-\lambda_+(Q))x} \int_x^\infty dx' \, p(x').
\end{split}
\end{align}
A straightforward computation shows that the function
\begin{equation}\label{via}
Y_+(x)=e^{\lambda_+(Q) x}f_\varepsilon(x)Z(x),\quad x\ge0,
\end{equation}
yields a solution of the perturbed differential equation
\eqref{perteqn}, and hence \eqref{zzz_1} and $p\in L^1(\bbR_+)$ imply
\begin{equation}\label{vesa1}
e^{-\varkappa'_+(Q)x}\|Y_+(x)-\Phi(x)Q \|_{\bbC^{d\times d}}
\underset{x\to  \infty}{=}o(1).
\end{equation}
In addition, using \eqref{falloff_1} and \eqref{fallooof_1} again,
we may choose $\delta>0$ sufficiently small such that the function
$p_\delta(x)=e^{\delta x}p(x)$ is integrable on $\bbR_+$. If $x\le
x'$, then $p(x')\le e^{-\delta x}p_\delta(x')$, and
thus \eqref{zzz_1} implies
\begin{align}
\begin{split}
& f_\varepsilon(x)\|Z(x)-Z^{(0)}(x)\|_{\bbC^{d\times d}}  \\
& \quad \le 2 \|Z\|_{C_{\rm b}(\bbR_+)^{d\times d}}
e^{(\varkappa_+'(Q)-\lambda_+(Q))x} e^{-\delta x}
\int_x^\infty dx' \, p_\delta(x').
\end{split}
\end{align}
This gives a better asymptotic relation than \eqref{vesa1},
\begin{equation}
e^{-\varkappa'_+(Q)x}\|Y_+(x)-\Phi(x)Q \|_{\bbC^{d\times d}}
\underset{x\to  \infty}{=}o(e^{-\delta x}),
\end{equation}
leading to
 \begin{equation}
\limsup_{x\to  \infty } \frac{\log\|Y_+(x)-\Phi(x)Q
\|_{\bbC^{d\times d}}}{x}\le\varkappa'_+(Q)-\delta<\varkappa'_+(Q),
\end{equation}
proving the existence part of the theorem.

\noindent {\it Step 3. Continuity}. Let
 $\{R_n\}_{n=1}^\infty$ be a sequence of $d\times d$ matrix-valued
 functions converging to $R$ in
 $L^1(\bbR_+;e^{\beta x }dx)^{d\times d}$ as $n\to\infty$.
On the Banach space $C_{\rm b}(\bbR_+)^{d\times d}$
introduce Volterra integral operators $V_n$ by
\begin{equation}
(V_nZ)(x)=\int_x^\infty dx' \,
 e^{-\lambda_+(Q)(x-x')}f_\varepsilon(x)^{-1}
 f_\varepsilon(x')\Phi(x)\Phi(x')^{-1}R_n(x')Z(x'), \quad n\in\bbN.
\end{equation}
An estimate similar to \eqref{amper} shows that the sequence of
the operators $V_n$ converges in operator norm to the integral operator $V$,
 which in turn,  yields  convergence in $C_{\rm b}(\bbR_+)^{d\times d}$ of
the unique solutions $Z_n$ of the Volterra integral equations
\begin{equation}
Z_n(x)=Z^{(0)}(x)-(V_nZ_n)(x), \quad x\ge 0,
\end{equation}
 to the  unique solution $Z$ of the integral equation \eqref{volt}.
 In particular,
\begin{equation}
\lim_{n\to  \infty}Z_n(0)=Z(0) \,\text{ in } \,
 \bbC^{d\times d}\,\text{ and hence }\,
 \lim_{n\to  \infty}Y_{+,n}(0)=Y_+(0),
\end{equation}
where, similarly to  \eqref{via}, we denote $
Y_{+,n}(x)=e^{\lambda_+(Q) x}f_\varepsilon(x)Z_n(x)$. This proves 
continuity of the mapping \eqref{cont_1}.
\end{proof}

\begin{remark}\label{netdich+}
This theorem holds under a weaker assumption
than the exponential dichotomy hypothesis (cf.\eqref{EDBE}). In fact, it is
sufficient to require that
$\varkappa_+(Q)<\varkappa_+'(I_d-Q)$ and $ \varkappa_+(Q)<0$ only.
\hfill$\Diamond$
\end{remark}

\begin{remark}\label{inteq}
Replacing assumptions \eqref{falloff_1} and \eqref{fallooof_1} by the
hypothesis
\begin{equation}
\|R\|_{\bbC^{d\times d}}\in L^1(\bbR_+;e^{\beta x }dx) \,
\text{ with } \, \beta> -\lambda_+'(Q),  \lb{3.29}
\end{equation}
the matrix-valued Jost solution $Y_+$ can be uniquely
determined by solving the Volterra integral equation
\begin{equation}\label{VENM}
Y_+(x)=\Phi(x)Q-\int_x^\infty dx' \, \Phi(x) \Phi(x')^{-1}R(x')Y_+(x'),
\quad x\ge 0; 
\end{equation}
that is, in this case, the ``modified"   Volterra integral equation
\eqref{volt_1}, \eqref{volt} on
$C_{\rm b}(\bbR_+)^{d\times d} $ is   not
needed. Indeed, assume $\beta> -\lambda_+'(Q)$, choose
$\varepsilon$ such that 
$0<\varepsilon<\min\{\beta+\lambda'_+(Q),\varkappa'_+(I-Q)\}$, and
use the inequalities $x'\ge x\ge0$ and
$\lambda'_+(Q)<0<\varkappa_+'(I-Q)$ to estimate the integral kernel of the integral
operator in \eqref{VENM} as follows:
\begin{equation}\label{NEW3.31}
\begin{split}
& \|\Phi(x)\Phi(x')^{-1}R(x')\|_{\bbC^{d\times d}} \le 
\Big(\|\Phi(x)Q\|_{\bbC^{d\times d}}\|Q\Phi(x')^{-1}\|_{\bbC^{d\times d}}\\
&\hspace*{4.3cm}    +\|\Phi(x)(I-Q)\Phi(x')^{-1}\|_{\bbC^{d\times d}}\Big)
\|R(x')\|_{\bbC^{d\times d}}\\
&\quad \le  \Big(C(\varepsilon)e^{(-\lambda'_+(Q)+\varepsilon)x'}
+C(\varepsilon)e^{(\varkappa'_+(I-Q)-\varepsilon)(x-x')}\Big)
\|R(x')\|_{\bbC^{d\times d}}\\
&\quad \le  2C(\varepsilon)e^{-(\lambda'_+(Q)
+\beta-\varepsilon)x'}e^{\beta
x'}\|R(x')\|_{\bbC^{d\times d}}
\le 2C(\varepsilon)e^{\beta x'}\|R(x')\|_{\bbC^{d\times d}}.
\end{split}
\end{equation}
We recall that $\varkappa'_+(Q)\le\lambda'_+(Q)\le\lambda_+(Q)<0$.
Condition \eqref{fallooof_1} is formulated in terms of the difference
$\lambda_+(Q)-\varkappa'_+(Q)$ while condition \eqref{3.29} is 
formulated in terms of the quantity $-\lambda'_+(Q)$. These two numbers,
$\lambda_+(Q)-\varkappa'_+(Q)$ and $-\lambda'_+(Q)$ are fairly independent,
and thus either of the two conditions, \eqref{fallooof_1} and 
\eqref{3.29}, can be better than the other. Thus, the existence of 
matrix-valued Jost solutions follows under the assumption
$\beta>\min\{\lambda_+(Q)-\varkappa'_+(Q), -\lambda'_+(Q)\}$.
\hfill$\Diamond$
\end{remark}

\begin{remark}\label{lala}
Example \ref{ex} below shows that the statement of  Theorem
\ref{+} is optimal in the  sense that, under the exponential
dichotomy hypothesis, the assumption in \eqref{fallooof_1} on the
exponential fall-off of the perturbation $R$ in general
cannot be relaxed: Indeed, if condition \eqref{fallooof_1} is
violated, then the matrix-valued Jost  solution may not exist.
Example \ref{ex} also shows that if Definition \ref{jost} is
modified to replace the strict inequality in \eqref{lapun+} by the
condition
\begin{equation}\label{flex}
\limsup_{x\to  \infty } \frac{\log\|Y_+(x)-\Phi(x)Q
\|_{\bbC^{d\times d}}}{x}\le\varkappa_+'(Q),
\end{equation}
then the matrix-valued Jost solution in this modified
sense is generally not unique. \hfill$\Diamond$
\end{remark}

\begin{example}\label{ex}
Consider equations \eqref{unperturbed} and \eqref{perteqn} with
\begin{equation}
A(x)=
 \begin{pmatrix}
-2&0&0\\
0& -1&0\\
0&0&1
 \end{pmatrix} \, \text{ and } \,
R(x)=  \begin{pmatrix}
0&\chi_{\bbR_+}(x)e^{-x}\cos x&0\\
0& 0&0\\
0&0&0
 \end{pmatrix}, \quad x\in \bbR,
\end{equation}
where $\chi_{\bbR_+}(x)$ denotes the characteristic function of
$\bbR_+=[0,\infty)$. The fundamental solution $\Phi$, the exponential
dichotomy $Q$ on $\bbR$, and the bounded solution $\Phi \, Q$
of \eqref{unperturbed} on $\bbR_+$ are given as follows:
\begin{align}
& \Phi(x)=\begin{pmatrix}
e^{-2x}&0&0\\
0& e^{-x}&0\\
0&0&e^{x}
\end{pmatrix}, \quad
Q=\begin{pmatrix}
1&0&0\\
0& 1&0\\
0&0&0
\end{pmatrix}, \quad \Phi(x)Q=\begin{pmatrix}
e^{-2x}&0&0\\
0& e^{-x}&0\\
0&0&0
\end{pmatrix}, \no \\
& \hspace*{10cm}  x\in\bbR.
\end{align}
The Bohl and Lyapunov exponents on  $\bbR_+$ associated with  the
dichotomy projection $Q$ are $ \varkappa_+'(Q)=\lambda_+'(Q)=-2$
and $\varkappa_+(Q)=\lambda_+(Q)=-1$. Clearly, $\|R\|_{\bbC^{d\times
d}}\in L^1(\bbR_+;e^{\beta x}dx)$ with
$\beta=\lambda_+(Q)-\varkappa_+'(Q)=1$, and hence condition
\eqref{fallooof_1} is violated. A direct computation shows that
the perturbed equation \eqref{perteqn} has two linearly
independent bounded $\bbC^d$-valued solutions on $\bbR_+$
\begin{equation}
y_1(x)= \begin{pmatrix}
e^{-2x}\\
0\\
0
 \end{pmatrix},
\quad  y_2(x)= \begin{pmatrix}
e^{-2x}\sin x\\
e^{-x}\\
0
 \end{pmatrix}, \quad x\ge 0,
\end{equation}
and thus any bounded $\bbC^d$-valued solution of
\eqref{perteqn} on $\bbR_+$ is a linear combination of these two. Hence, any
bounded matrix-valued  solution $Y$ of
\eqref{perteqn} on $\bbR_+$ is necessarily of the form
\begin{equation}\label{would}
Y(x)= \begin{pmatrix}
e^{-2x}&e^{-2x}(\sin x+C)&0\\
0&e^{-x}&0\\
0&0&0
 \end{pmatrix}, \quad x\ge 0,
\end{equation}
where $C\in\bbC$ is an arbitrary constant. Thus, it satisfies
\begin{equation}\label{ploxo}
\limsup_{x\to  \infty } \frac{\log\|Y(x)-\Phi(x)Q
\|_{\bbC^{d\times d}}}{x}=-2=\varkappa_+'(Q).
\end{equation}
This observation shows that the $d\times d$ matrix-valued Jost
solution as introduced in Definition \ref{jost} does not exist if
\eqref{fallooof_1} fails, while the solutions satisfying
\eqref{flex} are not unique. \hfill$\Diamond$
\end{example}

\begin{remark}\label{ravno2}
In the proof of Theorem \ref{+}  we suggested a method of
introducing the Jost solution by $Y_+(x)=e^{\mu x}Z(x)$, $x\geq 0$, where
$Z\in C_{\rm b}(\bbR_+)^{d\times d}$ is a
solution of the ``$\mu$-modified'' Volterra integral equation
\begin{equation}\label{volt2}
 Z(x)=e^{-\mu x}\Phi(x)Q-\int_x^\infty dx' \,
 e^{-\mu(x-x')}\Phi(x)\Phi(x')^{-1}R(x')Z(x') , \quad x\geq 0,
 \end{equation}
for an appropriate choice of $\mu\in \bbR$. We stress that this
method fails if the exponential fall-off hypothesis
\eqref{fallooof_1} is violated. Indeed, for the equations in
Example \ref{ex},
 if $\mu< \lambda_+(Q)=-1$, the term $e^{-\mu
x}\Phi(x)Q$ in the integral equation \eqref{volt2} is not a
bounded function, while if $\mu\ge  \varkappa_+(Q)=-1$, the
Volterra integral operator is unbounded, both in the Banach space
$C_{\rm b}(\bbR_+)^{d\times d}$ and in the
Hilbert space $L^2(\bbR_+)^{d\times d}$. \hfill$\Diamond$
\end{remark}

Similar to the proof of Theorem \ref{+} (replacing
$\varkappa'_+(Q)$ by $\varkappa_-(I_d-Q)$ and $\lambda_+(Q)$ by
$\lambda_-'(I_d-Q)$) one proves the following result for the
negative half-line $\bbR_-$.

\begin{theorem}\label{-}
Assume Hypothesis \ref{diffeq} and the condition
\begin{equation}
\|R\|_{\bbC^{d\times d}}\in L^1(\bbR_-;e^{\beta |x| }dx)
\, \text{ for some } \, \beta>\varkappa_-(I_d-Q)-\lambda_-'(I_d-Q).  
\lb{new3.38}
\end{equation}
 Then the perturbed differential equation
in \eqref{perteqn} has a unique matrix-valued Jost  solution
$Y_-(x)$, $x\le 0$. Moreover, the map
\begin{equation}\label{cont}
\begin{cases}
L^1(\bbR_-;e^{\beta |x| })^{d\times d} \to \bbC^{d\times d}  \\
\hspace*{2.27cm}  R \mapsto Y_-(0)  \end{cases}
\end{equation}
is continuous.
\end{theorem}

\section{The Perturbation Determinant}\label{PDSEC}

In this section we employ results from \cite{GM}, succinctly summarized in
Appendix \ref{AP1}, to express the infinite-dimensional (modified) Fredholm
determinant of the integral operator $I+K$, given by \eqref{KK},
via a finite-dimensional determinant obtained by means of the
matrix-valued Jost solutions introduced in Definition \ref{jost}.
This result is a generalization of the classical relation between
a Fredholm determinant  and the Jost function due to Jost and Pais
\cite{JP51} (see, also
\cite[Sect.\ 12.1.2]{Ne02}). As we will see in Section \ref{EVFUN}, the
finite-dimensional determinant is related to the Evans function associated
with
\eqref{unperturbed} and \eqref{perteqn}. We recall Lemma \ref{tech2}, and
also the general fact (see, e.g., \cite[Ch.\ XIII]{GGKr}, 
\cite[Sect.\ IV.2]{GK69}, \cite{Si77}, \cite[Ch.\ 9]{Si79}) that if $K$ is
any Hilbert--Schmidt operator, then the 2-modified Fredholm perturbation
determinant $\DT(I+K)$ is given by the formulas
\begin{equation}
\DT(I+K)=\det\big((I+K)e^{-K}\big)=
\prod_{\lambda\in\spec(K)}\big((1+\lambda)e^{-\lambda}\big),
\end{equation}
counting algebraic multiplicites of the eigenvalues of $K$. Moreover, if 
$K$ is a trace class operator, then
\begin{equation}\label{NEWtr1}
{\det}_2(I+K)=e^{-\tr(K)}\det(I+K),\, \text{ where }\,
\det(I+K)=\prod_{\lambda\in\sigma(K)}(1+\lambda).
\end{equation}

Throughout, we will use the following notation:
\begin{equation}\label{theta}
\Theta=
\int^\infty_0 dx\, {\tr}_{\bbC^d} \big(\Phi(x)Q\Phi(x)^{-1}R(x)\big)
 -\int^0_{-\infty} dx\, {\tr}_{\bbC^d}
\big(\Phi(x)(I_d-Q)\Phi(x)^{-1}R(x)\big).
\end{equation}

\begin{theorem}\label{detcomp}
Assume Hypothesis \ref{diffeq}  and the condition
\begin{equation}
\|R\|_{\bbC^{d\times d}}\in L^1(\bbR;e^{\beta |x| }dx)   \lb{4.3}
\end{equation}
for some
\begin{equation}\label{fall}
\beta>\max\Big\{\lambda_+(Q)-\varkappa'_+(Q),
\varkappa_-(I_d-Q)-\lambda_-'(I_d-Q)\Big\}.
\end{equation}
Let $K$ be the integral
 operator on $L^2(\bbR)^d$ whose integral kernel is given by \eqref{KK}.
Then the  2-modified perturbation determinant $\DT(I+K)$ admits
the representation
\begin{equation}  \label{detdetcomp}
\DT(I+K)=e^{\Theta} {\det}_{\bbC^d}(Y_+(0)+Y_-(0)), 
\end{equation}
where $\Theta$ is defined in \eqref{theta},
and $Y_\pm$ are the matrix-valued Jost solutions
 on $\bbR_\pm$, respectively.
\end{theorem}
\begin{proof} Assume temporarily that
\begin{equation}\label{supR}
\text{supp}\,(R) \, \text{ is compact}.
\end{equation}
Then $\|R\|_{\bbC^{d\times d}}\in L^1(\bbR;e^{\beta |x| }dx)$ for
any real
$\beta$. In particular, by Remark \ref{inteq} and its obvious
$\bbR_-$-analog, the $d\times d$ matrix-valued Jost solutions
 $Y_\pm$ are the unique bounded solutions of
 the following Volterra integral equations on $\bbR_\pm$,
\begin{equation}\label{3.9.5}
\begin{split}  Y_+(x)&= \Phi (x)Q-\int^\infty_x dx' \,
\Phi (x)\Phi (x')^{-1}R(x')Y_+(x'),
\quad x\ge 0,\\
Y_-(x)&=
 \Phi (x)(I_d-Q)+\int^x_{-\infty} dx' \,
\Phi (x)\Phi (x')^{-1}R(x')Y_-(x'),  \quad x\le 0.
\end{split}\end{equation}
We note that $Y_+(x)=Y_+(x)Q$ and $Y_-(x)=Y_-(x)(I_d-Q)$ by 
uniqueness of the solutions. Moreover, by \eqref{supR} we have
$Y_+(x)=\Phi(x)Q$ for $x\ge n$ and $Y_-(x)=\Phi(x)(I_d-Q)$ for $x\le
-n$ for sufficiently large $n\in\bbN$.

Treating matrices as operators on respective spaces,
we introduce the following notations for $x\in\bbR$,
\begin{equation}\label{14.2} \begin{split}
f_1(x)&=R_r(x) \Phi (x)Q:\im (Q)\to \bbC^d,\\  f_2(x)& =R_r(x)
 \Phi (x)(I_d-Q) :\ker (Q)\to \bbC^d,\\
g_1(x)&= Q\Phi (x)^{-1}  R_\ell(x):\bbC^d\to \im (Q), \\ g_2(x)& =
 -(I_d-Q)\Phi (x)^{-1}R_\ell(x):\bbC^d\to \ker (Q).
\end{split}
\end{equation}
In addition, let $d_1=\dim(\im (Q))$ and $d_2=\dim(\ker (Q))$, and denote
\begin{equation*}\begin{split} H(x,x')=f_1(x)g_1(x')-f_2(x)g_2(x')
= R_r(x)\Phi (x)\Phi (x')^{-1}R_\ell (x'),\quad
x,x'\in\bbR.\end{split}
\end{equation*}
We note that the functions $f_j$ and $g_j$, $j=1,2$,
are compactly supported on $\bbR$ due to \eqref{supR}, and thus
assumption \eqref{2.1} holds. Hence, the results recorded in Appendix
\ref{AP1} are at our disposal. Next, we introduce the
Volterra integral equations
\begin{equation}\label{some }
\begin{split} \hat{f}_1 (x) &=f_1(x)-\int^\infty_x
dx' \, H(x,x')\hat{f}_1(x'), \quad x\in\bbR, \\
\hat{f}_2(x)&=f_2(x)+\int^x_{-\infty}  dx' \, H(x,x')\hat{f}_2(x'),
\quad x\in\bbR,
\end{split}
\end{equation}
and let $\hat f_1(x):\im (Q)\to \bbC^d$ and $\hat f_2(x):\ker (Q)\to
\bbC^d$ be the unique solutions of \eqref{some } that satisfy
$\hat f_j\in L^2(\bbR)^{d\times d_j}$, $j=1,2$ (cf.\
Appendix \ref{AP1}). The $d\times d$ matrix $U(x)$, $x\in\bbR$,
defined in \eqref{2.37}, written as a block-operator with
respect to the decomposition
 $\bbC^d=\im (Q)\dot+\ker (Q)$, then reads as  follows:
\begin{equation}\label{19.3}
U(x)= \begin{pmatrix}
I_{d_1}-\int^\infty_x dx'\, g_1(x')\hat{f}_1(x') & \int^x_{-\infty} dx' \,
g_1(x')\hat{f}_2(x')  \\
\int^\infty_x dx' g_2(x')\hat{f}_1(x') & I_{d_2} -\int^x_{-\infty} dx' \,
g_2(x')\hat{f}_2 (x') \end{pmatrix}.
\end{equation}
Writing
$R(x)=R_\ell(x)R_r(x)$, and multiplying \eqref{3.9.5} by $R_r(x)$
from the left, we arrive at \eqref{some } and observe that
\begin{equation}\label{14.3new} \begin{split}
\hat f_1(x)&= R_r(x)Y_+(x)|_{\im (Q)},\quad  x\ge0, \\
\hat f_2(x) &=R_r(x)Y_-(x)|_{\ker (Q)},\quad
x\le0.
\end{split}
\end{equation} Setting $x=0$ in
\eqref{3.9.5}, using \eqref{14.2} and \eqref{14.3new}, and writing
the matrix $Y_+(0)+Y_-(0)$ as a block-operator with respect to the
decomposition $\bbC^d=\im (Q)\dot+\ker (Q)$, we conclude from
\eqref{19.3} at $x=0$ that $U(0)=Y_+(0)+Y_-(0)$. Since  the integral kernel
$K(x,x')$  in \eqref{KK} can be represented as
\begin{equation}
K(x,x')  = \begin{cases} -f_1(x)g_1(x'),& x\ge x',\\
-f_2(x)g_2(x'), & x<x',\end{cases} \quad x, x' \in \bbR,
\end{equation}
 Theorem \ref{t2.6}\,$(ii)$  proves \eqref{detdetcomp} for compactly
supported perturbations. Specifically, one can use \eqref{3.30} with
$x_0=0$ and then apply some elementary properties of matrix traces.

To remove the compact support assumption on the matrix-valued
function $R$  one proceeds as follows. Given $n\in \bbN$, we
introduce the truncations
\begin{equation}\label{TRUNC}
R_n(x)=\begin{cases} R(x), & |x|\le n,\\
0, & |x|>n.
\end{cases}
\end{equation}
Let $R_n(x)=V_n(x)|R_n(x)|$, $x\in \bbR$,
 be the  polar decompositions of $R_n(x)$, denote
 $ R_{\ell,n}(x)=V_n(x)|R_n(x)|^{{1}/{2}}$
and $R_{r,n}(x)=|R_n(x)|^{{1}/{2}}$, and introduce the integral kernels of
Hilbert--Schmidt integral operators $K_n$ on $L^2(\bbR)^d$ by
\begin{equation}\label{KKn}
K_n(x,x')=\begin{cases}
-R_{r,n}(x) \Phi(x)Q\Phi(x')^{-1} R_{\ell,n}(x'),& x\ge x',\\
R_{r,n}(x)\Phi(x)(I_d-Q)\Phi(x')^{-1} R_{\ell,n}(x'),&x< x',
\end{cases}  \quad x, x' \in\bbR.
\end{equation}
Since the support of $R_n$ is compact, one infers by the first step of
the proof that
\begin{equation}
\DT(I+K_n)=e^{\Theta_n} {\det}_{\bbC^d}(Y_{+,n}(0)+Y_{-,n}(0)),\quad n\in
\bbN,
\end{equation}
where
\begin{align}
\begin{split}
\Theta_n&=
\int^\infty_0 dx \, {\tr}_{\bbC^d} \big(\Phi(x)Q\Phi(x)^{-1}R_n(x)\big) \\
& \hspace*{.4cm} -\int^0_{-\infty} dx \, {\tr}_{\bbC^d}
\big(\Phi(x)(I_d-Q)\Phi(x)^{-1}R_n(x)\big)
\end{split}
\end{align}
and $Y_{\pm,n}$ are the matrix-valued Jost solutions of the truncated
perturbed equation $y'=(A+R_n)y$ on $\bbR_\pm$. Since $R_n$ converges
to $R$ in $L^1(\bbR)^{d\times d}$ as $n\to \infty$,
 one concludes that
\begin{equation}\label{raz}
\lim_{n\to  \infty}\Theta_n=\Theta,
\end{equation}
where $\Theta$ is given by \eqref{theta}. Using the estimate in
the proof of Lemma \ref{tech2}, one checks that the operators
$K_n$ converge to the operator $K$ in Hilbert--Schmidt norm,
and hence (see, e.g., \cite[Ch.\ IX]{GGKr}, \cite[Ch.\ 9]{Si79})
\begin{equation}\label{dva}
\lim_{n\to  \infty}\DT(I+K_n)=\DT(I+K).
\end{equation}
Finally, since the sequence  $R_n$ converges to $R$
in $L^1(\bbR; e^{\beta|x |}dx)^{d\times d}$ as $n\to \infty$, applying
assertions \eqref{cont_1} and \eqref{cont} of  Theorems \ref{+}
and \ref{-}, we obtain
 \begin{equation}\label{tri}
 \lim_{n\to  \infty}Y_{\pm,n}(0)=Y_\pm(0).
 \end{equation}
 Combining \eqref{raz}--\eqref{tri} then completes the proof.
\end{proof}

\begin{remark}\label{remNEWTRACE}
The following heuristic argument may be helpful in understanding
the role of the factor $e^\Theta$ in \eqref{detdetcomp} and
shows that the appearance of this factor is quite natural. First,
we observe that the integral kernel \eqref{KK} of the operator $K$
is generally discontinuous on the diagonal $x=x'$. Accordingly, the
definition of $\Theta$ in \eqref{theta} can be re-written as
follows:
\begin{equation}\label{NEWtr2}
-\Theta=\int_0^\infty \,dx \,
\tr_{\mathbb{C}^d}(K(x+0,x))+\int_{-\infty}^0 \, dx\,
\tr_{\mathbb{C}^d}(K(x-0,x)).
\end{equation}
Heuristically, the right-hand side of \eqref{NEWtr2} can be
viewed as a ``regularized integral trace'' of the operator $K$.
Next, we assume in addition that $K$ is a trace class operator 
(this requires quite a stretch of imagination!), and that $\tr (K)$ equals 
the integral trace $-\Theta$, 
\begin{equation}
\tr (K)=-\Theta. 
\end{equation}
Then {\it formally} applying \eqref{NEWtr1} and \eqref{detdetcomp} yields
\begin{equation}
\DT(I+K)=e^{-\tr(K)}\det(I+K)=e^{\Theta} {\det}_{\bbC^d}(Y_+(0)+Y_-(0))
\end{equation}
and hence {\it formally},
\begin{equation}
\det(I+K) = {\det}_{\bbC^d}(Y_+(0)+Y_-(0)).
\end{equation}
As we will see in
Section \ref{EVFUN}, the determinant ${\det}_{\bbC^d}(Y_+(0)+Y_-(0))$ is
in fact the Evans function for the equations \eqref{unperturbed} and
\eqref{perteqn}. \hfill$\Diamond$
\end{remark}

\begin{remark}\label{look_1}
The proof of Theorem \ref{detcomp} shows that the $2$-modified
determinant of the operator $I+K$ can be computed by the formula
\begin{equation}\label{formula}
\DT(I+K)=e^{\Theta} \lim_{n\to 
\infty}{\det}_{\bbC^d}(Y_{+,n}(0)+Y_{-,n}(0))\\
\end{equation}
under the weaker assumption $\|R\|_{\bbC^{d\times d}}\in
L^1(\bbR)$ as opposed to \eqref{4.3} and \eqref{fall}. Here 
$Y_{\pm,n}$ are the unique matrix-valued Jost solutions of the
Volterra integral equations
\begin{align}
\begin{split}
Y_{+,n}(x)&=\Phi(x)Q-\int_x^\infty dx' \,
\Phi(x)\Phi(x')^{-1}R_n(x')Y_{+,n}(x')
,\quad x\ge 0,\label{volttrunc}\\
Y_{-,n}(x)&=\Phi(x)(I_d-Q)+\int_{-\infty}^x dx' \,
\Phi(x)\Phi(x')^{-1}R_n(x')Y_{-,n}(x'), \quad x\le 0,
\end{split}
\end{align}
and $R_n$ is the truncated perturbation \eqref{TRUNC}.
We also have the following identities:
\begin{equation}\label{NEW4.16}
Y_{+,n}(0)=Y_{+,n}(0)Q, \quad
Y_{-,n}(0)=Y_{-,n}(0)(I_d-Q), \end{equation} since the solutions
$Y_{\pm,n}$ are unique.\hfill$\Diamond$
\end{remark}

 \begin{remark}
 The exponential fall-off assumption \eqref{fall} has been
 imposed in Theorem \ref{detcomp} to make
sure that the matrix-valued Jost solutions $Y_\pm$ exist
and are unique (cf.\ Remark \ref{lala}). Moreover, under this assumption,
one uses continuity of the mapping \eqref{cont_1} to conclude that
\begin{align}
\begin{split}
\lim_{n\to  \infty}{\det}_{\bbC^d}(Y_{+,n}(0)+Y_{-,n}(0))&
={\det}_{\bbC^d}\Big(\lim_{n\to  \infty}(Y_{+,n}(0)+Y_{-,n}(0))\Big)
\label{inter_1}\\
&={\det}_{\bbC^d}(Y_+(0)+Y_-(0)).
\end{split}
\end{align}
Example \ref{ex} shows that
 the exponential fall-off assumptions \eqref{fall} are indeed needed
to interchange the limit $n\to \infty$ and the determinant  in
\eqref{inter_1}. In fact, in this example one immediately verifies that the
Jost solutions $Y_{\pm,n}$ associated with the truncated perturbation
$R_n$ are given by
\begin{align}
\begin{split}
Y_{+,n}(x)&=\begin{cases}
 \begin{pmatrix}
e^{-2x}&e^{-2x}[\sin (x)-\sin (n)])&0\\
0&e^{-x}&0\\
0&0&0
 \end{pmatrix}, & \, 0\le x \le n, \\[1mm]
 \begin{pmatrix}
e^{-2x}&0&0\\
0&e^{-x}&0\\
0&0&0
 \end{pmatrix}, & \,  x > n,
\end{cases} \\
Y_{-,n}(x)&= \begin{pmatrix} 0&0&0\\0&0&0\\0&0&e^x \end{pmatrix},
\quad x\le0.
\end{split}
\end{align}
 Clearly, the sequence of the $3\times 3$ matrices
 $Y_{+,n}(0)$
 does not converge as $n\to  \infty$, while $\lim_{n\to 
\infty} {\det}_{\bbC^d}( Y_{+,n}(0)+Y_{-,n}(0))$ does exist
 since  ${\det}_{\bbC^d}( Y_{+,n}(0)+Y_{-,n}(0))=1$ for all $n\in \bbN$.
Consequently,  \eqref{inter_1} fails.
 \hfill$\Diamond$
\end{remark}

\section{Sub-exponential Weights and the Jost Function}\label{SES}

As we have seen in Section \ref{secjost}, the exponential fall-off
hypotheses \eqref{fallooof_1} and \eqref{new3.38} on the perturbation
$R$ cannot be relaxed in general because the unperturbed equation
might have at least one of the following two properties: First,
the upper Lyapunov exponent associated with the dichotomy
projection $Q$ may not coincide with the lower Bohl exponent and
then, necessarily, the corresponding Bohl interval is of  positive
width; and, second, the estimate \eqref{DEFVAR} used in the
definition of the lower Bohl exponent $\varkappa'(Q)$ may not hold
for $\varkappa=\varkappa'(Q)$ but only for
$\varkappa=\varkappa'(Q)-\varepsilon$ for any $\varepsilon>0$. 
This results in the {\em de facto} presence of an exponential weight factor
 $f_\varepsilon(|x-x'|)=C(\varepsilon)e^{\varepsilon|x-x'|}$ on the
right-hand side of \eqref{DEFVAR}; cf.\ 
 also \eqref{expscaling} and the effect of this on the proof of Theorem
 \ref{+}. In this section we restrict our attention to the class of 
unperturbed equations \eqref{unperturbed} that do not have either of these
properties. That is, we consider the case when each of the Bohl
segments associated with the dichotomy projections $Q$ and $(I_d-Q)$
degenerates into a single point and, in addition, we assume that the
exponential weight factors in the estimates such as \eqref{DEFVAR}
are replaced by {\em sub-exponential} weight factors induced by a given
monotone weight function $f$.

We introduce a {\em weight function}, $f$, which, by definition, is
 a nondecreasing function $f: \bbR_+\to \bbR_+$,
satisfying the following conditions:
\begin{equation}
\limsup_{x\to  \infty}\frac{\log (f(x))}{x}=0, \quad f(0)\ge 1.
\end{equation}
\begin{hypothesis}\label{furt}
Assume Hypothesis \ref{diffeq}. Suppose that the Bohl segments on
$\bbR$ associated with the projections $Q$ and $(I_d-Q)$ have zero
width, that is,
\begin{equation}
\varkappa'(Q)=\varkappa(Q)<0<\varkappa'(I_d-Q)=\varkappa(I_d-Q).
\end{equation}
Given a weight function $f$, we asssume, in addition, that the
following estimates hold for all $x, x' \in\bbR$:
\begin{equation}
\begin{split}\| \Phi (x)Q\Phi (x')^{-1}\|_{\bbC^{d\times d}} & \leq
e^{\varkappa(Q) (x-x')}f(|x-x'|), \quad 0\le x\le x', \\
\label{5.1}\| \Phi (x)(I_d-Q)\Phi (x')^{-1}\|_{\bbC^{d\times d}}  &\leq
e^{\varkappa(I_d-Q)(x-x')} f(|x-x'|),\quad 0\ge x\geq x'.
\end{split}\end{equation}
\end{hypothesis}

Given Hypothesis \ref{furt} we now
introduce the matrix-valued Jost solutions as solutions  of
certain ``modified'' Volterra integral equations (cf.\
\eqref{volt_1}--\eqref{volt} and \eqref{3.9.5}).

\begin{definition}\label{FGjost}
Assume Hypothesis \ref{furt} and suppose
\begin{equation}
\|R\|_{\bbC^{d\times d}}\in L^1(\bbR; f^2(|x|)dx).
\end{equation}
Then matrix-valued solutions $Y_\pm$ of
\eqref{perteqn} on $\bbR_\pm$ are called {\em matrix-valued
Jost solutions} of \eqref{perteqn} if
\begin{align}
\begin{split}
Y_+(x)&=f(|x|)e^{\varkappa(Q)x}Z_+(x), \quad x\ge 0,\\
Y_-(x)& =f(|x|)e^{\varkappa(I_d-Q)x}Z_-(x), \quad x\le 0,
\end{split}
\end{align}
where $Z_+(x)$, $x\ge0$, and $Z_-(x)$, $x\le0$, are the unique
bounded solutions of the Volterra
integral equations
\begin{align}
Z_+(x)&=e^{-\varkappa(Q)x}f(x)^{-1}\Phi(x) Q  \no \\
& \hspace*{-.5cm} -\int_x^\infty dx' \,
e^{-\varkappa(Q)(x-x')}f(x)^{-1}f(x')
\Phi(x)\Phi(x')^{-1}R(x')Z_+(x'), \quad x\geq 0, \no \\
Z_-(x)&=e^{-\varkappa(I_d-Q)x}f(|x|)^{-1}\Phi(x)(I- Q) \label{FTSYS} \\
& \hspace*{-.5cm} +\int_{-\infty}^x dx' \,
e^{-\varkappa(I_d-Q)(x-x')}f(|x|)^{-1}f(|x'|)
\Phi(x)\Phi(x')^{-1}R(x')Z_-(x'), \quad  x\leq 0, \no
\end{align}
on $\bbR_+$ and $\bbR_-$.
\end{definition}

Almost literally repeating the proof of Theorems \ref{+} and
\ref{-} (replacing $f_\varepsilon$ by $f$),
one concludes that the matrix-valued Jost solutions $Y_\pm$
in the sense of Definition \ref{FGjost} satisfy the following
asymptotic relations (cf.\ \eqref{vesa1}):
\begin{align}
\begin{split}
& e^{-\varkappa(Q)x}\|Y_+(x)-\Phi(x)Q\|_{\bbC^{d\times d}}
\underset{x\to  \infty}{=}o(1),\\
& e^{-\varkappa(I_d-Q)x}\|Y_-(x)-\Phi(x)(I_d-Q)\|_{\bbC^{d\times d}}
\underset{x\to  -\infty}{=}o(1).
\end{split}
\end{align}

Moreover, one obtains the following version of  Theorem
\ref{detcomp}.

\begin{theorem}\label{detcomt} Assume Hypothesis \ref{furt}, suppose that
\begin{equation}
\|R\|_{\bbC^{d\times d}}\in L^1(\bbR;
f^2(|x|)dx),
\end{equation}
and let the integral kernel of the operator $K$ be given by
\eqref{KK}. Then the 2-modified perturbation determinant $\DT(I+K)$ admits
the representation
\begin{equation}\label{ES5.2}
\DT(I+K)=e^{\Theta} {\det}_{\bbC^d}(Y_+(0)+Y_-(0)),\\
\end{equation}
where $\Theta$ is defined in \eqref{theta} and $Y_\pm$ are
the matrix-valued Jost solutions on $\bbR_\pm$ introduced in
Definition \ref{FGjost}.
\end{theorem}

In the remaining part of this section we will apply Theorem
\ref{detcomt} to study the Schr\"odinger equation
\eqref{schr} with an integrable potential,
\begin{equation}
V\in L^1(\bbR).
\end{equation}
We consider equations \eqref{unperturbed} and \eqref{perteqn} with
\begin{equation}\label{SE0}
A(k)= \begin{pmatrix} 0 & 1\\ -k^2 & 0 \end{pmatrix},
\quad  R(x)= \begin{pmatrix} 0 & 0\\ V(x) & 0 \end{pmatrix}, \quad
x\in\bbR.
\end{equation}
To avoid confusion we emphasize that $A$ is an $x$-independent function
of the parameter $k$ in \eqref{schr}. Clearly,
\eqref{perteqn}, with $A=A(k)$ and $R$ as in \eqref{SE0}, is the
first-order system corresponding to the Schr\"odinger equation
\eqref{schr}. We note that $\spec(A(k))=\{ik,-ik\}$, where, without loss
of generality, we choose $k$ so that $\Im(k)\ge0$. Since we intend to
apply Theorem
\ref{detcomt}, we have to make sure that Hypothesis \ref{diffeq}
is satisfied. In particular, the unperturbed equation with the
$x$-independent coefficient $A(k)$ must have an exponential
dichotomy on $\bbR$. This is equivalent to the requirement
$\spec(A(k))\cap i\,\bbR=\emptyset$ or, taking into account the
choice of $k$, to the inequality $\Im(k)>0$, which we will assume
to hold in the remaining part of this section. 

Let $Q(k)$ be the spectral projection for $A(k)$ so that
$\spec(A(k)|_{\im (Q(k))})=\{ik\}$. We remark that
\begin{equation}
\varkappa'(Q(k))=\varkappa(Q(k))=\re(ik)
\, \text{ and } \,
\varkappa'(I_d-Q(k))=\varkappa(I_d-Q(k))=\re(-ik).
\end{equation}
Next, we introduce the
sub-exponential weight function by
\begin{equation}
f(x)=c, \quad x\in\bbR_+,
\end{equation}
for an appropriate constant $c\geq 1$. Since the eigenvalues $\pm ik$ of
$A(k)$ are simple, the estimate \eqref{5.1} with $\Phi(k,x)=e^{xA(k)}$,
$x\in\bbR$, holds, and thus Hypothesis \ref{furt} is satisfied. The
matrix-valued Jost solutions $Y_\pm(k,\cdot)$ in the sense of
Definition \ref{FGjost} are the unique solutions of the Volterra
integral equations
\begin{equation}\label{SE2}
\begin{split}
Y_+(k,x)&=e^{xA(k)}Q(k)-\int_x^\infty dx' \,
e^{(x-x')A(k)}R(x')Y_+(k,x'),\quad x\ge0,\\
Y_-(k,x)&=e^{xA(k)}(I_d-Q(k))+\int^x_{-\infty} dx' \,
e^{(x-x')A(k)}R(x')Y_-(k,x'),\quad x\le0,
\end{split}
\end{equation}
such that the matrix-valued functions $Z_\pm(k,x)=e^{\mp
ikx}Y_\pm(k,x)$, $x\in\bbR_\pm$, are bounded. Since $V\in
L^1(\bbR)$, the conclusion of Theorem \ref{detcomt} holds with
$Y_\pm(k,\cdot)$ given by \eqref{SE2}.

Next, we will relate the finite-dimensional determinant in \eqref{ES5.2} and
the classical Jost function (see, e.g., \cite[Ch.\ XVII]{CS}, \cite[Sect.\
12.1]{Ne02} for the latter). First, we recall some well-known notions from
scattering theory
 (see \cite{GM} for a detailed bibliography).
 The {\em Jost solutions}, $u_\pm(k,\cdot)$,
 of the Schr\"odinger equation \eqref{schr} are defined as
 solutions of the  Volterra integral equations
 \begin{align}\label{JEqn}
 \begin{split}
& u_\pm(k,x)=e^{\pm ikx}-\int_x^{\pm\infty} dx' \, k^{-1}\sin (k(x-x'))
V(x')u_\pm(k,x'), \\
& \hspace*{6.75cm} \Im(k) > 0, \; x\in\bbR.
 \end{split}
 \end{align}
The Jost function, $\mathcal{J}=\mathcal{J}(k)$ is defined by
\begin{equation}
\mathcal{J}(k)=\frac{1}{2ik}W(u_-(k,\cdot),u_+(k,\cdot)), \quad \Im(k)>0, 
\end{equation}
where $W(u(x),v(x))=u(x)v'(x)-u'(x)v(x)$ is the Wronskian of $u$ and $v$
with $u, v \in C^1(\bbR)$. For $k> 0$, $\cJ$(k) is the reciprocal
of the transmission coefficient (see, e.g., formula (XVII.1.36) in
\cite{CS}). We note that $\mathcal{J}(k)$ is independent of
$x\in\bbR$ since ${\tr}_{\bbC^2} (A(k))=0$.

Next, we consider the factorization $V=V_\ell V_r$ with
\begin{equation}
V_\ell(x)=|V(x)|^{1/2}, \quad\
V_r(x)=|V(x)|^{1/2}\exp(i\arg(V(x)), \quad x\in\bbR.
\end{equation}
To make the connection with \eqref{FACTOR}, we remark that
\begin{equation}\label{SE2.2}
R_\ell(x)= \begin{pmatrix}0&0\\V_\ell(x)&0 \end{pmatrix}, \quad
R_r(x)= \begin{pmatrix}V_r(x)&0\\0&0 \end{pmatrix},\quad x\in\bbR.
\end{equation}
Finally, we introduce the integral operator, $L(k)$ on $L^2(\bbR)$
with integral kernel
\begin{equation}\label{SE2.1}
L(k,x,x')=\frac{i}{2k}V_r(x)e^{ik|x-x'|}V_\ell(x'),\quad
\Im(k) > 0, \; x,x'\in\bbR.
\end{equation}

The following corollary of Theorem \ref{detcomt} recovers a
well-known relation between the Fredholm determinant of the
operator $I+L(k)$ and the Jost function originally due to Jost and Pais
\cite{JP51} (see also \cite{Ne72}, \cite{Ne80}, \cite[Sect.\ 12.1.1]{Ne02}
and, specifically, \cite[Theorem 4.7]{GM}). More importantly, it shows that
the classical Jost function coincides with the finite-dimensional
determinant
${\det}_{\bbC^d}(Y_+(k,0)+Y_-(k,0))$ in \eqref{ES5.2} obtained by means of
the matrix-valued Jost solutions introduced in Definition
\ref{FGjost}. As we will see in Section \ref{EVFUN}, this
determinant is, in fact, the Evans function associated with
equations \eqref{unperturbed} and \eqref{perteqn}.

\begin{theorem}\label{SEMAIN}
Assume that the potential of the Schr\"odinger equation
\eqref{schr} satisfies $V\in L^1(\bbR)$, and fix $k\in\bbC$ such
that $\Im (k) > 0$.
Let $K(k)$ be the integral operator on $L^2(\bbR)^2$ with the integral
kernel given by \eqref{KK}, where $\Phi(k,x)=e^{xA(k)}$, $x\in\bbR$,
and $A(k)$ is defined in \eqref{SE0}, and $R_\ell$ and
$R_r$ are defined by \eqref{SE2.2}. Let $L(k)$ be the integral
operator on $L^2(\bbR)$ with the integral kernel given by
\eqref{SE2.1} and $\mathcal{J}(k)$ be the Jost
function of \eqref{schr}. Finally, let $Y_\pm(k,\cdot)$
be the matrix-valued Jost solutions \eqref{SE2} on $\bbR_\pm$ for
\eqref{unperturbed} and \eqref{perteqn} with $A(k)$ and $R$ given
by \eqref{SE0}. Then the following assertions hold: \\
$(i)$  $\mathcal{J}(k)={\det}_{\bbC^2}(Y_+(k,0)+Y_-(k,0))$. \\
$(ii)$ The 2-modified Fredholm determinant admits the following
representation:
\begin{equation}\label{sTheta}
{\det}_2(I+K(k))={\det}_2(I+L(k))=e^{\Theta (k)} \mathcal{J}(k),
\end{equation}
where
\begin{equation}
\Theta (k)=\frac{1}{2ik}\int_\bbR dx \, V(x).  \lb{5.21}
\end{equation}
\end{theorem}
\begin{proof}
It is convenient to diagonalize $A(k)$ and $e^{xA(k)}$, $x\in\bbR$. To
this effect we introduce the matrices
\begin{align}
T(k) & = \begin{pmatrix}1&1\\ik&-ik \end{pmatrix}, \quad T(k)^{-1}  =
\frac{1}{2ik}
 \begin{pmatrix} ik & 1\\ ik & -1  \end{pmatrix},  \no \\
\widetilde{A}(k)& =T(k)^{-1} A(k) T(k)= \begin{pmatrix}ik&0\\0&-ik
\end{pmatrix},
\label{DEFATIL} \\
\widetilde{Q}& = \begin{pmatrix}1&0\\0&0 \end{pmatrix}, \quad S(k)
=\frac{1}{2ik} \begin{pmatrix}1&1\\-1&-1 \end{pmatrix},  \no  \\
Q(k)&=T(k)\widetilde{Q}T(k)^{-1}=\frac12 \begin{pmatrix}1&(ik)^{-1}\\ik&1
\end{pmatrix}.
\label{DEFATIL.1}
\end{align}
We note that $\widetilde{Q}$ is the spectral projection for
$\widetilde{A}(k)$ so that
\begin{equation}
\spec(\widetilde{A}(k)|_{\im
(\widetilde{Q})})=\spec(A(k)|_{\im (Q(k))})=\{ik\}.
\end{equation}
Passing
to the matrix-valued functions $\wti{Y}_\pm(k,x)=T(k)^{-1}Y_\pm(k,x)T(k)$,
$x\in\bbR_\pm$, in \eqref{SE2}, we obtain the integral equations
\begin{align}
\wti{Y}_+(k,x)&=e^{x\wti{A}(k)}\wti{Q}-\int_x^\infty dx' \,
e^{(x-x')\wti{A}(k)}V(x')S(k)\wti{Y}_+(k,x'),\quad
x\ge0,  \label{SE2.11}  \\
\wti{Y}_-(k,x)&=e^{x\wti{A}}(I-\wti{Q})+\int^x_{-\infty} dx' \,
e^{(x-x')\wti{A}(k)}V(x')S(k)\wti{Y}_-(k,x'), \quad
x\le0.\label{SE2.12}
\end{align}
The property $\wti{Y}_+(k,\cdot)=\wti{Y}_+(k,\cdot)\wti{Q}$
and \eqref{SE2.11} imply that the second column of the $2\times
2$ matrix $\wti{Y}_+(k,\cdot)$ is equal to zero, while
$\wti{Y}_-(k,\cdot)=\wti{Y}_-(k,\cdot)(I_2 -\wti{Q})$ and \eqref{SE2.12}
imply that the first column of $\wti{Y}_-(k,\cdot)$ is equal to zero. Thus,
in the matrix equations \eqref{SE2.11} and \eqref{SE2.12} we can
separately consider the first and the second column, respectively. Let
$\wti{y}_+(k,\cdot)$ denote  the first column of $\wti{Y}_+(k,\cdot)$ and
$\wti{y}_-(k,\cdot)$ denote the second column of $\wti{Y}_-(k,\cdot)$.
Passing in \eqref{SE2.11} and \eqref{SE2.12} to the $\bbC^2$-valued
functions $y_\pm(k,x)=T(k)\wti{y}_\pm(k,x)$, $x\in\bbR_\pm$, we
observe\footnote{We use $\top$ for transposition so that
$(a\;\; b)^\top$ is a $(2\times 1)$ column
vector.} that they satisfy the following $\bbC^2$-valued integral
equations:
\begin{align}\label{SE2.3}
y_+(k,x)&=e^{ikx} (1\;\;ik)^\top
-\int_x^\infty dx' \, V(x')T(k)e^{(x-x')\wti{A}(k)}S(k)T(k)^{-1}y_+(k,x'),
\no \\
& \hspace*{8.5cm}  x\ge0,    \\
y_-(k,x)&=e^{-ikx} (1\;\;-ik)^\top
+\int^x_{-\infty} dx' \,
V(x')T(k) e^{(x-x')\wti{A}(k)}S(k)T(k)^{-1}y_-(k,x'), \no \\
& \hspace*{10.3cm} x\le0.  \no
\end{align}
A direct calculation using \eqref{DEFATIL} and \eqref{DEFATIL.1} shows that
\begin{align}
\begin{split}
& T(k)e^{(x-x')\wti{A}(k)}S(k)T(k)^{-1}=
\begin{pmatrix}k^{-1}\sin(k(x-x'))&0 \\
\cos(k(x-x'))&0 \end{pmatrix},  \\
& \hspace*{5.35cm} \Im(k) > 0, \; x,x'\in\bbR.
\end{split}
\end{align}
Since
$Y_\pm(k,\cdot)$ are matrix-valued solutions of \eqref{perteqn}, it
follows that $y_\pm(k,\cdot)$ are $\bbC^2$-valued solutions of
\eqref{perteqn}. Next, we denote by $u_\pm(k,\cdot)$ the top entry of the
$2\times 1$ vector $y_\pm(k,\cdot)$. Differentiating the
first components in the $2\times 1$ vector integral equations in
\eqref{SE2.3}, we observe that
$y_\pm(k,x)= (u_\pm(k,x)\;\;u'_\pm(k,x))^\top$,
$x\in\bbR_\pm$, and thus $u_\pm(k,\cdot)$ are (scalar-valued) solutions of
the Schr\"odinger equation \eqref{schr}. Moreover, it follows
from \eqref{SE2.3} that $u_\pm(k,\cdot)$ satisfy \eqref{JEqn},
that is, $u_\pm(k,\cdot)$ are the traditional Jost solutions of
\eqref{schr}. This proves assertion $(i)$ since
\begin{align}
&{\det}_{\bbC^2}(Y_+(k,0)+Y_-(k,0))
={\det}_{\bbC^2}(\wti{Y}_+(k,0)+\wti{Y}_-(k,0))  \no \\
& \quad =
{\det}_{\bbC^2} ((\wti{y}_+(k,0)\;\;\wti{y}_-(k,0)))
={\det}_{\bbC^2}\big(T(k)^{-1} ({y}_+(k,0)\;\;{y}_-(k,0))\big)  \no \\
& \quad = \frac{-1}{2ik}W(u_+(k),u_-(k))=\mathcal{J}(k).
\end{align}
Finally, we turn to the proof of assertion $(ii)$. A
direct computation in \eqref{theta} using \eqref{SE0} and the
formula for $Q(k)$ in \eqref{DEFATIL.1}, verifies the formula for
$\Theta(k)$ in \eqref{5.21}.
Using \eqref{KK} and \eqref{SE2.2}, and with the help of the
diagonalization described in \eqref{DEFATIL} and \eqref{DEFATIL.1},
one computes that the matrix-valued integral kernel $K(k,x,x')$ in
\eqref{KK} and the scalar-valued integral kernel $L(k,x,x')$ in
\eqref{SE2.1} are related by the formula
\begin{equation}
K(k,x,x')=L(k,x,x')\wti{Q}, \quad x,x'\in\bbR,
\end{equation}
where $\wti{Q}$ is the projection in \eqref{DEFATIL.1}. This implies
$\DT(I+K(k))=\DT(I+L(k))$. Thus, relation
\eqref{ES5.2} and assertion $(i)$ yield assertion $(ii)$.
\end{proof}

While we assumed $\Im(k)>0$ throughout this section, we note that
continuity of $K(k)$ and $L(k)$ with respect to $k$, $\Im(k)\geq 0$, $k
\neq 0$, in the Hilbert--Schmidt norm, permits one to extend the results
of Theorem \ref{SEMAIN}\,$(i), (ii)$ to all $\Im(k)\geq 0$, $k\neq 0$, by
continuity.  

\section{Generalized Matrix-Valued Jost Solutions}\label{GMVJS}

In this section we start the discussion of a generalization of the
matrix-valued Jost solutions introduced in Definition \ref{jost},
and prove an extension of Theorem \ref{+} assuming that the
dichotomy projection $Q$ admits further exponential splitting. As
we will see later on, the generalized matrix-valued
Jost solutions will allow us to relax the exponential fall-off
hypothesis on the perturbation imposed in \eqref{fall}. We recall
from Remark \ref{lala} that these hypotheses are optimal, provided
one merely assumes the existence of an exponential dichotomy $Q$
with no further splitting. To simplify the exposition in this
section, we will only consider the case where the dichotomy
projection $Q$ admits a splitting of order two and postpone the
general case of the exponential splitting of arbitrary order until
Section \ref{ESAO}.

\begin{hypothesis}\label{dichsplit}
Assume Hypothesis \ref{diffeq}. In addition, assume that the
dichotomy projection $Q$ can be represented in the form
$Q=Q_1+Q_2$, where $Q_1$ and $Q_2$ are projections which are
uniformly conjugated by $\Phi$ on $\bbR$ such that the
Bohl segment $[\varkappa'(Q_1), \varkappa(Q_1)]$ lies strictly
below  $[\varkappa'(Q_2), \varkappa(Q_2)]$, that is,
\begin{equation}
\varkappa'(Q_1)\le \varkappa(Q_1)<\varkappa'(Q_2)\le
\varkappa(Q_2)<0<\varkappa'(I_d-Q)\le\varkappa(I_d-Q).
\end{equation}
\end{hypothesis}

We will consider the  generalized $d\times d$ matrix-valued Jost
solutions on $\bbR_+$ (cf.\ Definition \ref{jost}), the case of
$\bbR_-$ can be treated similarly.

\begin{definition}\label{jostsplit}
Assume $R\in L^1_{\loc}(\bbR)^{d\times d}$
and  Hypothesis \ref{dichsplit}. Then $d \times d$ matrix-valued
solutions $Y_+^{(1)}$ and $Y_+^{(2)}$ of the
perturbed differential equation \eqref{perteqn} on $\bbR_+$ are called 
{\em generalized matrix-valued Jost solutions} associated with the
projections $Q_1$ and $Q_2$ if
\begin{equation}\label{lapun}
\limsup_{x\to  \infty } \frac{\log\|Y_+^{(j)}(x)-\Phi(x)Q_j
\|_{\bbC^{d\times d}}}{x}<\varkappa_+'(Q_j)
\end{equation}
and
\begin{equation}\label{initial}
Y_+^{(j)}(0)=Y_+^{(j)}(0)Q_j, \quad j=1,2.
\end{equation}
\end{definition}

We emphasize that the additional  technical requirement in \eqref{initial}
can always be satisfied as soon as  \eqref{lapun} holds for some
solutions  $Y_+^{(j)}$, $j=1,2.$ Indeed,  if the solutions
$Y_+^{(j)}$ satisfy \eqref{lapun}, then
$Y_+^{(j)} Q_j$, $j=1,2$, are solutions of the perturbed
differential equation satisfying both \eqref{lapun} and
\eqref{initial}.

\begin{remark}\label{NUN}
We note that the generalized matrix-valued Jost solution
$Y^{(2)}_+$ is not unique. Indeed, if $Y^{(2)}_+$ is
a solution satisfying \eqref{lapun} for $j=2$, and $Y_+$ is
any solution of \eqref{perteqn} whose upper Lyapunov exponent on
$\bbR_+$ belongs to the Bohl segment
$[\varkappa'_+(Q_1),\varkappa_+(Q_1)]$, then the solution
$Y^{(2)}_+ +Y_+$ again satisfies \eqref{lapun} for $j=2$ since
\begin{align}
\begin{split}
& \limsup_{x\to  \infty }
\frac{\log\|Y_+^{(2)}(x)+Y_+(x)-\Phi(x)Q_2\|_{\bbC^{d\times d}}}{x}\\
& \quad \le  \max\Big\{\limsup_{x\to  \infty }
\frac{\log\|Y_+^{(2)}(x)
-\Phi(x)Q_2\|_{\bbC^{d\times
d}}}{x},\varkappa_+(Q_1)\Big\}<\varkappa'_+(Q_2).
\end{split}
\end{align}
\hfill $\Diamond$
\end{remark}

We start with the following elementary fact.

\begin{lemma}\label{nulek} Let $\alpha>0$ and assume that $0\leq p\in
L^1(\bbR_+)$. Then,
\begin{equation}
\lim_{x\to  \infty}e^{-\alpha x}\int_0^x dx' \, e^{\alpha x'}p(x')=0.
\lb{6.5}
\end{equation}
\end{lemma}
\begin{proof} Assuming $\text{mes} \, (\{x\geq 0 \,|\, p(x)>0\})>0$ and
introducing $q(x)=\int_0^x dx' \, p(x')$, $x\geq 0$, an integration by
parts yields
\begin{equation}
e^{-\alpha x}\int_0^x dx' \, e^{\alpha x'}p(x')=q(x)-
\alpha{e^{-\alpha x}}\int_0^x dx' \, e^{\alpha x'}q(x'), \quad x\geq 0.
\end{equation}
Since $p$ is nonnegative  and integrable on $\bbR_+$,
$\lim_{x\to \infty}q(x)$ is positive, and the integral $\int_0^x dx'
\,  e^{\alpha x'}q(x')$ diverges. Using l'H\^opital's rule
one obtains assertion \eqref{6.5}.
\end{proof}

The following result ensures the existence of the {\em generalized}
matrix-valued Jost solutions on $\bbR_+$ introduced in Definition
\ref{jostsplit}
under much weaker exponential fall-off assumptions on the
perturbation than \eqref{falloff_1} and 
\eqref{fallooof_1} (we recall that the requirement \eqref{falloff_1} is
optimal if we  want to deal with the matrix-valued solutions introduced in
Definition \ref{jost} only).

\begin{theorem}\label{+split}
Assume Hypothesis \ref{dichsplit} and the condition
\begin{equation}\label{falloff}
\|R\|_{\bbC^{d\times d}}\in L^1(\bbR_+;e^{\beta x }dx)
\end{equation}
with some
\begin{equation}\label{fallooof}\beta>\max_{j=1,2}\{\lambda_+(Q_j)
-\varkappa_+'(Q_j)\}.
\end{equation}
Then the perturbed differential equation \eqref{perteqn} has
generalized matrix-valued Jost solutions $Y_+^{(j)}$,
$j=1,2$.
\end{theorem}
 \begin{proof} Existence and uniqueness
 of the generalized Jost solution $Y_+^{(1)}$
 on $\bbR_+$ associated with the projection $Q_1$
  follows from Theorem \ref{+} and Remark \ref{netdich+};
  to check condition \eqref{initial} for $j=1$ we note that
  $Y_+^{(1)}(x)=Y_+^{(1)}(x)Q$, $x\ge0$, since the solution
 $Y_+^{(1)}$ is unique.

  By hypothesis and the definition of the lower Bohl exponents
(cf.\ \eqref{DEFVAR}), for any $\varepsilon
 >0$ there exists a positive constant $C(\varepsilon)>0$ such that
\begin{align}
 & \|\Phi(x)(I-Q_1)\Phi(x')^{-1}\|_{\bbC^{d\times d}}  \no \\
& \quad \le
 \|\Phi(x)Q_2\Phi(x')^{-1}\|_{\bbC^{d\times d}}
 +\|\Phi(x)(I_d-Q)\Phi(x')^{-1}\|_{\bbC^{d\times d}} \no \\
 & \quad \le
\big (e^{(\varkappa_+'(Q_2))(x-x')}+
e^{(\varkappa_+'(I_d-Q))(x-x')} \big)f_\varepsilon(|x-x'|),\quad
0\le x\le x',
\end{align}
where $f_\varepsilon(x)=C(\varepsilon)e^{\varepsilon |x|}$. Since
$\varkappa_+'(Q_2)<\varkappa_+'(I_d-Q)$, we infer the estimate
\begin{equation}\label{PhiPhi}
\|\Phi(x)(I-Q_1)\Phi(x')^{-1}\|_{\bbC^{d\times d}} \le
2e^{(\varkappa_+'(Q_2))(x-x')}f_\varepsilon(|x-x'|), \quad
0\le x\le x'.
\end{equation}
By the definition of the upper Bohl exponent, one may also assume that
\begin{equation}\label{PhiPhiQ1}
\|\Phi(x)Q_1\Phi(x')^{-1}\|_{\bbC^{d\times d}}\le
2e^{(\varkappa_+(Q_1))(x-x')}f_\varepsilon(|x-x'|), \quad
0\le x'\le x.
\end{equation}
Taking $\varepsilon>0$ such that
$0<2\varepsilon< \beta-(\lambda_+(Q_2)-\varkappa_+'(Q_2))$, and
introducing the function
\begin{equation}
p(x)=e^{(\lambda_+(Q_2) -\varkappa_+'(Q_2))x}
 f_\varepsilon^2(x)\|R(x)\|_{\bbC^{d\times d}},\quad x\ge 0,
\end{equation}
one observes that $0 \leq p\in L^1(\bbR_+)$ by assumptions
\eqref{falloff} and \eqref{fallooof}. For brevity we denote
$\mu=\lambda_+(Q_2)$,
$\alpha=\varkappa'_+(Q_2)-\varkappa_+(Q_1)$, and, for any
$\tau\ge0$, we introduce the integral kernel
\begin{equation}   \label{yadro}
F_\mu^\tau(x,x')=e^{-\mu(x-x')}
 \frac{f_\varepsilon(x')}{f_\varepsilon(x)}
\begin{cases}
 -\Phi(x)(I-Q_1)\Phi(x')^{-1}R(x'),&\tau\le x<
x', \\\Phi(x)Q_1\Phi(x')^{-1}R(x'),&\tau\le x'\le x.
\end{cases}
\end{equation}
Using \eqref{PhiPhi} and \eqref{PhiPhiQ1} and the inequalities
$\varkappa'_+(Q_2)-\lambda_+(Q_2)\le0$ and $\alpha>0$, one derives
 the following estimate:
\begin{align}
& \| F_\mu^\tau(x,x')\|_{\bbC^{d\times d}} \le
 2f_\varepsilon(x)^{-1}f_\varepsilon^2(x')\|R(x')\|_{\bbC^{d\times d}} \no\\
& \hspace*{2.8cm} \times  \begin{cases}
 e^{(\varkappa'_+(Q_2)-\lambda_+(Q_2))(x-x')},& \tau\le x<x', \\
 e^{(\varkappa_+(Q_1)-\lambda_+(Q_2))(x-x')},& \tau\le x'\le x,
\end{cases}\no\\
& \quad =2f_\varepsilon(x)^{-1}e^{(\varkappa'_+(Q_2)-\lambda_+(Q_2))
x}p(x')
\begin{cases}
 1,& \tau\le x<x', \\
 e^{-\alpha(x-x')},& \tau\le x'\le x,
\end{cases}\label{gutgut}\\
& \quad \le 2f_\varepsilon(x)^{-1}p(x'),\quad
x,x'\ge\tau.\label{gutgut_1}
\end{align}
On the space $C_{\rm b}([\tau,\infty))^{d\times d}$ we
define the integral operator $F_\mu^\tau$ by
\begin{equation}
(F_\mu^\tau Z)(x)=\int_\tau^\infty dx' \,
 F_\mu^\tau(x,x')Z(x'),\quad x\ge \tau,
\end{equation}
and consider the corresponding  Fredholm-type integral equation
 \begin{equation} \label{fred_1}
 Z(x)=Z^{(0)}(x)-(F_\mu^\tau Z)(x), \quad x\geq \tau,  
\end{equation}
where $Z^{(0)}$ is defined by
\begin{equation}
Z^{(0)}(x)=e^{-\mu x}f_\varepsilon^{-1}(x)\Phi(x)Q_2,\quad x\ge\tau.
 \end{equation}
Since $f_\varepsilon\ge 1$ and $p\in L^1(\bbR_+)$, the
estimate in \eqref{gutgut_1} shows that the integral operator
$F_\mu^\tau$ is a contraction on
$C_{\rm b}([\tau,\infty))^{d\times d}$ for $\tau\ge 0$
sufficiently large. Since the exponential decay of
$\|\Phi \, Q_2\|_{\bbC^{d\times d}}$ is controlled by $\mu$, we see that
$Z^{(0)}\in C_{\rm b}([\tau,\infty))^{d\times d}$.
Thus, the integral equation
\eqref{fred_1} has a unique solution
$Z\in C_{\rm b}([\tau,\infty))^{d\times d}$ that  can
be obtained by the iteration process
\begin{equation}\label{it1}
Z(x)=\sum_{j=1}^\infty Z^{(j)}(x), \quad
Z^{(j)}(x)= (F_\mu^\tau Z^{(j-1)})(x), \quad j\in\bbN, \;
x\ge\tau.
\end{equation}
Clearly, by uniqueness of the solution,
\begin{equation}\label{qt}
Z(x)=Z(x)Q_2, \quad x\geq \tau.
\end{equation}
A straightforward computation shows that the function
\begin{equation}\label{viaaa}
Y_+^{(2)}(x)=e^{\mu x}f_\varepsilon(x)Z(x), \quad x\ge \tau,
\end{equation}
yields a solution of the perturbed differential equation in
\eqref{perteqn}. Moreover, from \eqref{fred_1} and \eqref{gutgut}
it also follows that
\begin{align}\label{zzzz_1}
\begin{split}
& f_\varepsilon(x)\|Z(x) -Z^{(0)}(x)\|_{\bbC^{d\times d}} \le
2e^{(\varkappa'_+(Q_2)-\lambda_+(Q_2))x}\bigg(\|p\|_{L^1([x,\infty))}\\
& \quad + e^{-\alpha x}\int_0^x dx'  p(x') e^{\alpha x'} \bigg)
\|Z\|_{C_{\rm b}([\tau,\infty))^{d\times d}},\quad
x\ge\tau.
\end{split}
\end{align}
Using $p\in L^1(\bbR_+)$ and Lemma \ref{nulek},
\eqref{zzzz_1} implies the asymptotic relation
\begin{equation}\label{vano}
e^{-\varkappa_+'(Q_2)x}\|Y_+^{(2)}(x)-\Phi(x)Q_2\|_{\bbC^{d\times d}}
\underset{x\to  \infty}{=}o(1).
\end{equation}
Moreover, using \eqref{falloff} and \eqref{fallooof} again,
choose $\delta\in(0,\alpha)$ sufficiently small such that the function
$p_\delta(x)=e^{\delta x}p(x)$, $x\geq 0$, is integrable on $\bbR_+$. Then
\eqref{zzzz_1} implies
\begin{align}
\begin{split}
& f_\varepsilon(x)\|Z(x)
-Z^{(0)}(x)\|_{\bbC^{d\times d}} \le 2e^{-\delta
x}e^{(\varkappa'_+(Q_2)-\lambda_+(Q_2))x}
\bigg(\|p_\delta\|_{L^1([x,\infty))}   \\
& \quad + e^{-(\alpha-\delta) x}\int_0^x dx'  p_\delta(x')
e^{(\alpha-\delta) x'}\bigg) \|Z\|_{C_{\rm b}([\tau,\infty))^{d\times d}},
\quad  x\ge\tau.
\end{split}
\end{align}
Using $p_\delta\in L^1(\bbR_+)$ and Lemma \ref{nulek}
again, one obtains the asymptotic relation
\begin{equation}\label{vano1}
e^{-\varkappa_+'(Q_2)x}\|Y_+^{(2)}(x)-\Phi(x)Q_2\|_{\bbC^{d\times
d}} \underset{x\to  \infty}{=}o(e^{-\delta x}).
\end{equation}
Finally, the matrix-valued solution $Y_+^{(2)}$ can be
uniquely extended to the interval $[0, \tau)$ by solving the
initial value problem
\begin{equation}
{Y_+^{(2)}}'(x)=(A(x)+B(x))Y_+^{(2)}(x),\quad x\geq \tau, \quad
Y_+^{(2)}(\tau)=e^{\mu \tau }f_\varepsilon(\tau)Z(\tau).
\end{equation}
In addition, from \eqref{vano1} one derives \eqref{lapun}, and from
\eqref{qt} it follows that $Y_+^{(2)}(0)=Y_+^{(2)}(0)Q_2$,
completing the proof.
\end{proof}

In order to discuss uniqueness properties of the generalized
$d\times d$ matrix-valued Jost solution $Y_+^{(2)}$ up to lower-order terms
(cf.\  Remark \ref{NUN}), we recall that if $Q$ is the (unique)
exponential dichotomy projection for the unperturbed equation
\eqref{unperturbed} on $\bbR$, then equation \eqref{unperturbed} also has
an exponential dichotomy on
$\bbR_+$ with the same dichotomy projection $Q$. However, the corresponding
Bohl segment on $\bbR_+$ may be strictly smaller than that on $\bbR$, that
is, $[\varkappa_+'(Q),\varkappa_+(Q)] \subset
[\varkappa'(Q),\varkappa(Q)]$, see \eqref{INEQbohls}. By Lemma
\ref{stub_dich}, the perturbed differential equation with
$\|R\|_{\bbC^{d\times d}}\in L^1(\bbR_+)$ will also have an
exponential dichotomy on $\bbR_+$ with a dichotomy projection $P$ having 
the same Bohl exponents as the unperturbed differential
equation for the dichotomy projection $Q$. We emphasize that the
dichotomy projection $P$ is not unique: Only the subspace $\im (P)$ is
determined uniquely since it consists precisely of those initial  data
$y_0$ such that the $\bbC^d$-valued solutions of the initial value problem
\begin{equation}\label{INPR}
y'(x)=(A(x)+R(x))y(x) ,\quad x\geq 0, \quad y(0)=y_0,
\end{equation}
are bounded on $\bbR_+$. Moreover, by Lemma \ref{stab_split}, if
$Q=Q_1+Q_2$ is an exponential splitting for \eqref{unperturbed} on
$\bbR$ (and therefore on $\bbR_+$), then the dichotomy projection
$P$ for \eqref{perteqn} on $\bbR_+$ also admits an exponential splitting
$P=P_1+P_2$ with the same Bohl segments, that is,
$[\varkappa_+'(P_j),\varkappa_+(P_j)]=
 [\varkappa_+'(Q_j),\varkappa_+(Q_j)]$, $j=1,2$.
 The projection $P_1$ is not uniquely determined but its range is,
since $\im (P_1)$ consists of precisely those 
 initial  data $y_0$ such that the Lyapunov exponent of the
 $\bbC^d$-valued solution of \eqref{INPR}
 satisfies $\lambda_+(y)
 \le\varkappa_+(P_1)<\varkappa_+'(P_2)
 =\varkappa_+'(Q_2)$ (cf.\ Remark \ref{VVSOLS}).
 Thus, in view of Theorem \ref{+split} and Remark \ref{NUN},  we
have proved the following uniqueness result:

\begin{corollary}\label{Cor6.6}
Assume the hypotheses of Theorem  \ref{+split}. Then: \\
$(i)$ The generalized matrix-valued Jost solution
$Y_+^{(1)}$ associated with the projection $Q_1$ is
unique.  \\
$(ii$ If $Y_+^{(2)}$ and $\widetilde Y_+^{(2)}$ are
any two generalized matrix-valued Jost solutions associated with
the projection $Q_2$ then $ \im (Y_+^{(2)}(0)-\widetilde
Y_+^{(2)}(0))\subseteq \im (P _1)$.
\end{corollary}

\section{Exponential Splitting of Arbitrary Order}\label{ESAO}

In this section we consider the generalized matrix-valued Jost
solutions of \eqref{perteqn} in the general case where the
unperturbed equation \eqref{unperturbed} has an exponential splitting of
arbitrary order on $\bbR$.

\begin{hypothesis}\label{mnogo}
Assume Hypothesis \ref{diffeq} with  $d\ge 2$. Suppose, in
addition, that for some $d'$, $2\le d'\le d$ and $k_0$, $1\le
k_0\le d'-1$, the dichotomy projection $Q$ on $\bbR$ for
\eqref{unperturbed} admits an exponential splitting
$Q=\sum_{j=1}^{k_0}Q_j$ of order $k_0$, and the projection $(I_d-Q)$
admits an exponential splitting $ I-Q=\sum_{j=k_0+1}^{d'}Q_j$ of
order $d'-k_0$, where the projections $Q_j$,  $j=1, 2,\dots, d'$, are
uniformly conjugated by $\Phi$ on $\bbR$ and the
corresponding disjoint Bohl segments $[\varkappa'(Q_j),
\varkappa(Q_j)]$ are ordered as follows:
\begin{equation}
\varkappa'(Q_j)\le
\varkappa(Q_j)<\varkappa'(Q_{j+1})\le \varkappa(Q_{j+1}), \quad 1\le
j\le d'-1.
\end{equation}
\end{hypothesis}

\begin{definition}\label{jostsplitn}
Assume $R\in L^1_{\loc}(\bbR)^{d\times d}$ and Hypothesis \ref{mnogo}.
Then $d \times d$ matrix-valued solutions $Y_+^{(j)}$, $j=1,\dots, k_0$, on
$\bbR_+$ and $Y_-^{(j)}$, $j=k_0+1,\dots, d'$, on $\bbR_-$ of the
differential equation
\begin{equation}
Y'(x)=(A(x)+R(x))Y(x), \quad x\in \bbR_\pm,
\end{equation}
are called {\em generalized matrix-valued Jost
solutions} associated with the exponential splitting
$\{Q_j\}_{j=1}^{d'}$ if
\begin{align}\label{lapunn}
& \limsup_{x\to  \infty } \frac{\log\|Y_+^{(j)}(x)-\Phi(x)Q_j
\|_{\bbC^{d\times d}}}{x}<\varkappa_+'(Q_j),\\
&  Y_+^{(j)}(0)=Y_+^{(j)}(0)Q_j,\quad
j=1,2,\dots, k_0,  \label{qq+}
\end{align}
and
\begin{align}
& \liminf_{x\to  -\infty }
\frac{\log\|Y_-^{(j)}(x)-\Phi(x)Q_j \|_{\bbC^{d\times
d}}}{x}>\varkappa_-(Q_j),  \label{lapunn-}  \\
&  Y_-^{(j)}(0)=Y_-^{(j)}(0)Q_j, \quad
 j=k_0+1,\dots, d'.   \label{qq-}
\end{align}
\end{definition}

Assuming Hypothesis \ref{mnogo}, and referring to Lemmas
\ref{stub_dich}--\ref{stab_split}, let
\begin{equation}\label{P++--}
P_+=P_1+P_2+ \dots +P_{k_0}
\end{equation}
be the
exponential splitting of order $k_0$ of the dichotomy projection
$P_+$ for \eqref{perteqn} on $\bbR_+$ such that
\begin{equation}
[\varkappa_+'(P_j),\varkappa_+(P_j)]=[\varkappa_+'(Q_j),\varkappa_+(Q_j)],
\quad j=1,\dots,k_0.
\end{equation}
Similarly, let
\begin{equation}\label{P++--A}
P_-=P_{k_0+1}+P_{k_0+2}+ \dots +P_{d'}
\end{equation}
be the exponential splitting of
order $(d'-k_0)$ of the dichotomy projection $P_-$ for
\eqref{perteqn} on $\bbR_-$ such that
\begin{equation}
[\varkappa_-'(P_j),\varkappa_-(P_j)]=[\varkappa_-'(Q_j),\varkappa_-(Q_j)],
\quad j=k_0+1, \dots, d'.
\end{equation}
Our next result shows that the generalized
matrix-valued Jost solutions exist and are unique up to lower exponential
order terms.

\begin{theorem}\label{+splitn}
Assume Hypothesis \ref{mnogo} and the condition
\begin{equation}\label{falloffn}
\|R\|_{\bbC^{d\times d}}\in L^1(\bbR;e^{\beta |x| }dx)
\end{equation}
with some
\begin{equation}\label{fallooofn}\beta>
\max\Big \{ \max_{1\le j\le
k_0}\{\lambda_+(Q_j)-\varkappa_+'(Q_j)\}, \max_{k_0+1\le j\le
d'}\{\varkappa_-(Q_j)-\lambda_-'(Q_j)\} \Big \}.
\end{equation}
 Then the following conclusion holds: \\
$(i)$ The perturbed equation \eqref{perteqn} has
 generalized matrix-valued Jost  solutions in the sense of Definition
 \ref{jostsplitn}, $Y_+^{(j)}$,
$j=1,2,\dots, k_0$, on $\bbR_+$ and $Y_-^{(j)}$, $j=k_0+1,
\dots, d'$, on $\bbR_-$ associated with the exponential splitting
$\{Q_j\}_{j=1}^{d'}$.

For any set of generalized matrix-valued Jost  solutions in
the sense of Definition \ref{jostsplitn}, the following assertions
hold: \\
$(ii)$ The generalized matrix-valued Jost solutions
$Y_+^{(1)}$ and $Y_-^{(d')}$, associated with the
projections $Q_1$ and $Q_{d'}$, are  uniquely determined and
satisfy 
\begin{equation}\label{INCL0}
\im(Y_+^{(1)}(0))\subseteq\im (P_1), \quad
\im(Y_-^{(d')}(0))\subseteq\im (P_{d'}).
\end{equation}
$(iii)$ The generalized matrix-valued Jost  solutions
$Y_+^{(j)}$, $j=2,\dots, k_0$, and $Y_-^{(j)}$,
$j=k_0+1, \dots, d'-1$, satisfy
\begin{align}\label{INCL1}
\im(Y_+^{(j)}(0))\backslash\{0\}&\subseteq\im\bigg(\sum_{k=1}^{j}P_k\bigg)
\bigg\backslash\im\Big(\sum_{k=1}^{j-1}P_k\bigg),\quad j=2,\dots, k_0,\\
\im(Y_-^{(j)}(0))\backslash\{0\}&\subseteq
\im\bigg(\sum_{k=j}^{d'}P_k\bigg)
\bigg\backslash\im\Big(\sum_{k=j+1}^{d'}P_k\Big),\quad
j=k_0+1,\dots,d'-1.\label{INCL1-}
\end{align}
$(iv)$ If $Y_+^{(j)}$ and $\widetilde Y_+^{(j)}$
 are any two generalized matrix-valued Jost solutions associated with the
projections
$Q_j$ for $j=2,\dots, k_0$ and $Y_-^{(j)}$
and $\widetilde Y_-^{(j)}$ are any two generalized Jost
solutions associated with the projection $Q_{j}$ for
$j=k_0+1,\dots,d'-1$, then
\begin{align}\label{incl+}
\im \big (Y_+^{(j)}(0)-\widetilde Y_+^{(j)}(0)\big )&
\subseteq \im
\bigg (\sum_{k=1}^{j-1}P _k\bigg ),\quad
 j=2,\dots,k_0,   \\
\im \big (Y_-^{(j)}(0)-\widetilde
Y_-^{(j)}(0)\big )& \subseteq  \im \bigg (\sum_{k=j+1}^{d'}P
_k\bigg ),\quad
 j=k_0+1,\dots, d'-1.  \label{incl-}
\end{align}
\end{theorem}
\begin{proof} We will provide a sketch of the proof for the case
of $\bbR_+$ referring for details to the proof of Theorem
\ref{+split}. The case of $\bbR_-$ is treated similarly.

To prove $(i)$ and the uniqueness statement in $(ii)$, we fix
$j\in\{1,\dots,k_0\}$ and denote $\mu_j=\lambda_+(Q_j)+\varepsilon$,
$j=1,\dots,k_0$, for some
$\varepsilon >0$. In addition, for some $\tau \ge 0$, we introduce the
integral operator $F_{\mu_j}^{\tau, j}$  with integral kernel
\begin{equation}\label{FKERJ}
F_{\mu_{j}}^{\tau,
j}(x,x')=e^{-\mu_j(x-x')}\begin{cases}
 -\Phi(x) \bigg [ \sum_{k=j}^{d'}Q_k \bigg ]\Phi(x')^{-1}R(x'),& \tau\le x<
x',\\ \Phi(x)\bigg [\sum_{k=1}^{j-1} Q_k \bigg
]\Phi(x')^{-1}R(x'),& \tau\le x'\le x, \end{cases}
\end{equation}
on the Banach space $C_{\rm b}([\tau, \infty))^{d\times d}$.
Repeating the arguments in the proof of Theorem \ref{+split} under
hypothesis \eqref{fallooofn}, one concludes that the operator $
F_{\mu_{j}}^{\tau, j}$ is a contraction on
$C_{\rm b}([\tau,\infty))^{d\times d}$, provided $\tau$
is large enough and $\varepsilon $ is sufficiently small.
Thus, for each such $\tau$, the Fredholm-type integral equation
\begin{equation}  \label{fred}
 Z^{(j)}(x)=Z^{(j)}_0(x)-(F_{\mu_{j}}^{\tau, j}Z^{(j)})(x), \quad
x\geq \tau, 
\end{equation}
where $Z^{(j)}_0$ is defined by
\begin{equation}
Z^{(j)}_0(x)=e^{-\mu_{j}
x}\Phi(x)Q_j,\quad x\ge \tau,
 \end{equation}
has a unique bounded matrix-valued solution on $[\tau,\infty)$ that can be
obtained by iterations similar to \eqref{it1}. We note that if $j=1$, then
equation \eqref{fred} is a
 Volterra integral equation, and there is no need to pass to a large $\tau$
 to ensure the contraction property: Indeed, the Volterra integral operator
 $F_{\mu_{1}}^{\tau, 1}$ has zero spectral radius on
 $C_{\rm b}([\tau,\infty))^{d\times d}$ for any
$\tau\ge 0$ and one can start the iteration process to obtain the solution
even at  $\tau=0$ (cf.\ the proof of Theorem \ref{+}). In particular, this
shows the uniqueness part of
$(ii)$. If $j=2,\dots,k_0$, then the generalized Jost solution
$Y^{(j)}_{+}$ on $\bbR_+$ is first constructed on the
interval $[\tau, \infty)$ by $ Y^{(j)}_{+}(x)=e^{\mu_j
x}Z^{(j)}(x)$, $x\ge \tau$, and then extended to $[0,\tau)$ by
solving the initial value problem
\begin{equation}
\big ({Y^{(j)}_{+}}\big )'(x)=(A(x)+R(x))Y^{(j)}_{+}(x), \quad x\in
[0,\tau], \quad Y^{(j)}_{+}(\tau )=e^{\mu_j \tau}Z^{(j)}(\tau).
\end{equation}
One verifies as in the proof of Theorem \ref{+split} that
\eqref{lapunn} and \eqref{qq+} hold. (We note the fact that for
$j=2,\dots,k_0$, the generalized matrix-valued Jost solutions depend on
$\tau$). This yields the existence of the generalized Jost solutions
associated with the exponential splitting $Q=\sum_{j=1}^{k_0}Q_j$ of the
dichotomy projection $Q$.

Inclusion \eqref{INCL0} in $(ii)$, and assertions $(iii)$ and $(iv)$
follow from \eqref{lapunn}, \eqref{qq+} and the elementary
properties of Lyapunov exponents of $\bbC^d$-valued solutions of
differential equations listed in Remark \ref{VVSOLS}. Indeed,
turning to the proof of \eqref{INCL0} and \eqref{INCL1}, we first remark
that for $j=1,\dots,k_0$,
\begin{equation}
\limsup_{x\to  \infty } \frac{\log\|Y_+^{(j)}(x)
\|_{\bbC^{d\times d}}}{x}\le\max\bigg\{\limsup_{x\to  \infty }
\frac{\log\|Y_+^{(j)}(x)-\Phi(x)Q_j \|_{\bbC^{d\times d}}}{x},
\lambda_+(Q_j) \bigg\}
\end{equation}
does not exceed $\varkappa_+(Q_j)=\varkappa_+(P_j)$ by
\eqref{lapunn}. Thus, the Lyapunov exponent $\lambda_+(y)$ of
the $\bbC^d$-valued solution $y$ of \eqref{perteqn} with initial
data $y(0)\in\im\big(Y_+^{(j)}(0)\big)$ is strictly smaller than
$\lambda'_+(P_{j+1})$ for $j=1,\dots,k_0-1$ and is negative for
$j=k_0$. Applying \eqref{IMPINT}  to the
perturbed equation \eqref{perteqn}, we obtain the inclusion
$\im\big(Y_+^{(j)}(0)\big)\subseteq\im\big(\sum_{k=1}^{j}P_k\big)$ for
$j=1,\dots,k_0$. To finish the proof of \eqref{INCL1} for
$j=2,\dots,k_0$, let us suppose that there is a nonzero vector
$y_0=Y_+^{(j)}(0)y_1$, with some $y_1\in\bbC^d$, such that
$y_0\in\im\big(\sum_{k=1}^{j-1}P_k\big)$. By \eqref{qq+} we have
$Q_jy_1\neq 0$. Then the Lyapunov exponent of the $\bbC^d$-valued
solution $y$ of the initial value problem \eqref{INPR}
given by $y(x)=Y_+^{(j)}(x)y_1$, $x\ge0$, must satisfy
$\lambda_+(y)\le\varkappa_+(P_{j-1})$ by Remark \ref{VVSOLS}. This
leads to a contradiction proving \eqref{INCL1}. Indeed, adding and
subtracting $y(x)=Y_+^{(j)}(x)y_1$ in $\|\Phi(x)Q_jy_1\|_{\bbC^{d}}$,
we conclude from \eqref{lapunn} that
\begin{align}
& \limsup_{x\to  \infty } \frac{\log(\|\Phi(x)Q_jy_1\|_{\bbC^{d}})}{x}
\no \\ 
& \quad \le\max\bigg\{\limsup_{x\to  \infty }
\frac{\log\|Y_+^{(j)}(x)-\Phi(x)Q_j
\|_{\bbC^{d\times d}}}{x},\varkappa_+(P_{j-1})\bigg\}  \no  \\
& \quad < \varkappa_+'(Q_j).   
\end{align}
But this is impossible since the Lyapunov
exponent of the nonzero $\bbC^d$-valued solution $\Phi \, Q_jy_1$ of
\eqref{unperturbed} on $\bbR_+$ with initial data
$Q_jy_1\in\im (Q_j)$ must belong to the Bohl segment
$[\varkappa'_+(Q_j),\varkappa_+(Q_j)]$ (see again Remark \ref{VVSOLS}).

Finally, to verify assertion $(iv)$ we add and subtract $\Phi(x)Q_j$ in
$Y_+^{(j)}(x)-\wti{Y}_+^{(j)}(x)$ and use \eqref{lapunn} to conclude that
$\lambda_+(y)<\lambda'_+(P_j)$ for every $\bbC^d$-valued solution $y$ of
\eqref{perteqn} with
$y(0)\in\im\big(Y_+^{(j)}(0)-\wti{Y}_+^{(j)}(0)\big)$. This proves
inclusion \eqref{incl+}.
\end{proof}

\begin{remark}\label{CONTRN}
In addition to assertions $(ii)$--$(iv)$ in Theorem \ref{+splitn}, the
generalized matrix-valued Jost solutions $Y^{(j)}_{\pm}$,
constructed in the existence part $(i)$ of Theorem \ref{+splitn} by
means of the Fredholm-type integral equations \eqref{fred} and
their analogs for $\bbR_-$, have the following property: If
$R_n$, $n\in\bbN$, is the sequence of truncated
perturbations as in \eqref{TRUNC}, and $Y^{(j)}_{\pm,n}$
denote the generalized matrix-valued Jost solutions associated
with the perturbation $R_n$ (see the proof of Theorem
\ref{detcomp}), then
\begin{equation}\label{contarg}
  \lim_{n\to  \infty}Y^{(j)}_{\pm,n}(0)=Y^{(j)}_{\pm}(0).
  \end{equation}
To see this one follows the course of the proof of Theorem \ref{+split}. 
First one establishes continuity of the map $R \mapsto
Z^{(j)}$ from $L^1(\bbR_+;e^{\beta x}dx)$ to
$C_{\rm b}([\tau,\infty))^{d\times d}$, defined in a
neighborhood ${\cU}(R)$ of $R$ in the space $L^1(\bbR_+;e^{\beta x}dx)$.
Subsequently, one shows
 that the convergence of the sequence
 $\{R_n\}_{n\in\bbN}$ from this neighborhood to
$R$  as $n\to  \infty$  in the $L^1(\bbR_+;e^{\beta x}dx)$-topology
yields \eqref{contarg}. \hfill$\Diamond$
\end{remark}

 Next, given any set
of generalized matrix-valued Jost solutions, we introduce the
following {\em Evans determinant} (the terminology is related to
the Evans function, which is further discussed in Section \ref{EVFUN}).

\begin{definition}\label{DEFEVANSDET}
For a given set of generalized matrix-valued Jost solutions $Y_+^{(j)}$,
$j=1,\dots,k_0$, on $\bbR_+$ and $Y_-^{(j)}$,
$j=k_0+1\dots,d'$, on $\bbR_-$, the {\em Evans determinant}, $D$,
is defined by
\begin{equation}\label{DEFDD}
D={\det}_{\bbC^d} ({\cY}_++{\cY}_-),\,\text{ where } \, {\cY}_+=
\sum_{j=1}^{k_0} Y_+^{(j)}(0), \quad {\cY}_-
=\sum_{j=k_0+1}^{d'} Y_-^{(j)}(0).\end{equation}
\end{definition}

The following important and purely algebraic result shows that
although the generalized matrix-valued Jost solutions are not
unique, the Evans determinant $D$ is uniquely determined by 
equations \eqref{unperturbed} and \eqref{perteqn}.

\begin{lemma}\label{mix} Assume Hypothesis \ref{mnogo}.
Then the determinant $D$ in \eqref{DEFDD} is independent of the
choice of the generalized matrix-valued Jost solutions
$Y_+^{(j)}$, $j=1,\dots, k_0$, and $Y_-^{(j)}$,
$j=k_0+1\dots, d'$.
\end{lemma}
\begin{proof}
Using the notation \eqref{P++--} and \eqref{P++--A}, let $P_+$ and $P_-$ be
the dichotomy projections and $\{P_j\}_{j=1}^{k_0}$ and
$\{P_j\}_{j=k_0+1}^{d'}$ be the exponential splittings for \eqref{perteqn}
on $\bbR_+$ and
$\bbR_-$, having the same Bohl segments as the
splitting $\{Q_j\}_{j=1}^{d'}$ for \eqref{unperturbed}. Without
loss of generality, we will assume that the subspaces $\im (P_j)$,
$j=1,\dots, d'$, satisfy $ \im (P_1)\dot +\im (P_2)\dot +\cdots \dot +\im
(P_{d'})=\bbC^d$ (otherwise the determinant $D$ is equal  to zero
for any choice of the system of the generalized Jost solutions).
We note that $ \im ({\cY}_\pm)\subseteq \im (P_\pm) $ by assertions
$(ii)$ and $(iii)$ in Theorem \ref{+splitn}, and ${\cY}_+={\cY}_+Q$ and
${\cY}_-={\cY}_-(I_d-Q)$ by \eqref{qq+} and
\eqref{qq-}. Thus, treating the matrices ${\cY}_\pm$ in
\eqref{DEFDD} as operators ${\cY}_+:\im (Q)\to \im (P_+)$ and $
{\cY}_-:\im (I_d-Q)\to \im (P_-)$, we observe that the
block-operator $\cY=\cY_+ + \cY_-$ is diagonal in the following direct sum
decomposition:
\begin{equation}
{\cY}= \begin{pmatrix}
{\cY}_+&0\\
0&{\cY}_-
 \end{pmatrix}:\bbC^d=\im (Q) \dot +\im (I_d-Q)\to \bbC^d=\im (P_+)\dot
+\im (P_-).
\end{equation}
Assertions $(ii)$ and $(iii)$ of  Theorem \ref{+splitn} then also
yield that the matrix of the operator
\begin{align}\label{D1}
\begin{split}
& {\cY}_+:\im (Q)=\im (Q_1)\dot +\cdots \dot +\im (Q_{k_0}) \\
& \hspace*{1cm} \to \im (P_+)=\im
(P_1) \dot + \cdots \dot +\im (P_{k_0})
\end{split}
\end{align}
 is upper-triangular, while the matrix of the operator
 \begin{align}  \label{D2}
\begin{split}
& {\cY}_-:\im (I_d-Q)=\im (Q_{k_0+1})\dot +\cdots \dot +\im (Q_{d'})  \\
& \hspace*{1cm} \to
\im (P_+)=\im (P_{k_0+1}) \dot +\cdots \dot +\im (P_{d'})
\end{split}
\end{align}
 is lower-triangular. Moreover, from assertion $(iv)$ in Theorem
\ref{+splitn}
 it also follows that the diagonal blocks of the operators
 ${\cY}_\pm$  with respect to the decompositions \eqref{D1} and \eqref{D2}
are independent of the choice of the generalized Jost solutions,
 completing the proof.
\end{proof}

\begin{remark}\label{TRIANG}
In fact, the proof of Lemma \ref{mix} only uses properties
\eqref{qq+} and \eqref{qq-}, and properties \eqref{INCL0}--\eqref{INCL1-}
of the matrix-valued solutions of \eqref{perteqn}.
\hfill$\Diamond$
\end{remark}

At this point we are ready to prove the principal result of this paper.

\begin{theorem}\label{posl}
Assume Hypothesis \ref{mnogo} and the condition
\begin{equation}\label{falloffnA}
\|R\|_{\bbC^{d\times d}}\in L^1(\bbR;e^{\beta |x| }dx)
\end{equation}
for some
\begin{equation}\label{fallooofnA}
\beta> \max\Big \{ \max_{1\le j\le
k_0}\{\lambda_+(Q_j)-\varkappa_+'(Q_j)\}, \max_{k_0+1\le j\le
d'}\{\varkappa_-(Q_j)-\lambda_-'(Q_j)\} \Big \}.
\end{equation}
Let $Y_+^{(j)}$, $j=1,2,\dots, k_0$,  and $Y_-^{(j)}$,
$j=k_0+1,2,\dots, d'$,  be any system of matrix-valued
generalized Jost solutions. Let $K$ be the integral operator on
$L^2(\bbR)^d$ whose integral kernel is given by \eqref{KK}.
Then the  2-modified perturbation determinant $\DT(I+K)$ admits
the representation
\begin{equation}\label{NEW7.12}
\DT(I+K)=e^{\Theta} {\det}_{\bbC^d}({\cY}_++{\cY}_-),
\end{equation}
where $\Theta$ is defined in \eqref{theta}, $ {\cY}_\pm$ are
defined in \eqref{DEFDD}, and $D={\det}_{\bbC^d}({\cY}_++{\cY}_-)$ is
the Evans determinant $($cf.\ \eqref{DEFDD}$)$.
\end{theorem}
\begin{proof}
 Let  $Y_{+, F}^{(j)}$, $j=1,\dots, k_0$, and
 $Y_{-, F}^{(j)}$, $j=k_0+1,\dots, d'$,
 be the generalized matrix-valued Jost solutions on $\bbR_+$, respectively,
on  $\bbR_-$, constructed in the existence part $(i)$, $(ii)$ of the
proof of Theorem \ref{+splitn} for a given $\tau\ge 0$ sufficiently large.
Here, the subscript $F$ is added to remind the reader that these solutions
are obtained from solutions of the Fredholm-type integral equations
\eqref{fred} and their analogs for $\bbR_-$. In general we note that the
solutions $Y_{\pm, F}^{(j)}$ depend on $\tau$, and they are {\em not} the
solutions
$Y_{\pm}^{(j)}$ given {\em a priori} in the formulation of
Theorem \ref{posl}. Introduce the truncated
perturbations $R_n$, $n\in\bbN$, by formula \eqref{TRUNC}, and denote by
$Y_{+,n,F}^{(j)}$, $j=1, \dots, k_0$, and
$Y_{-,n,F}^{(j)}$, $j=k_0+1,\dots, d'$, the corresponding
generalized matrix-valued Jost solutions of the truncated
perturbed equation $y'(x)=(A(x)+R_n(x))y(x)$, obtained by solving
the Fredholm-type integral equations \eqref{fred} with
$R$ replaced by $R_n$ (see the proof of
Theorem \ref{+splitn}\,$(i)$). Since $R_n$ is compactly supported,
there exist unique matrix-valued Jost solutions in the sense
of Definition \ref{jost} of the truncated perturbed equation on
$\bbR_+$ and on $\bbR_-$; these will be denoted by
$Y_{+,n}$ and $Y_{-,n}$. Using the solutions
$Y_{\pm,n}$ and the projections $Q_j$, define matrix-valued
solutions of the truncated perturbed equation by
\begin{equation}\label{DEFYV}
\begin{split}Y_{+,n,V}^{(j)}(x)&=
 Y_{+,n}(x)Q_j,\,\,j=1,\dots, k_0,\,\,x\ge0,\\
Y_{-,n,V}^{(j)}(x)&=Y_{-,n}(x)Q_{j},\,\, j=k_0+1,\dots,
d',\,\,x\le0. \end{split}\end{equation} Here, the subscript $V$ is
added to remind the reader that the solutions $Y_{\pm,n}$ are obtained
from the Volterra integral equations \eqref{volttrunc}. Since
the support of $R_n$ is compact, equations
\eqref{volttrunc} imply that the solutions $Y_{+,n,V}^{(j)}$
constitute a  system of {\em generalized} matrix-valued Jost
solutions in the sense of Definition \ref{jostsplitn}. In addition
to \eqref{DEFDD}, introduce the notations
\begin{eqnarray}
{\cY}_{+,F}&=\sum\limits_{j=1}^{k_0} Y_{+, F}^{(j)}(0),
\quad\quad  \quad \;\;{\cY}_{-,F}
& =\sum\limits_{j=k_0+1}^{d'} Y_{-, F}^{(j)}(0),\\
{\cY}_{+,n,F}& = \sum\limits_{j=1}^{k_0}Y_{+,n,F}^{(j)}(0),
\quad\quad  {\cY}_{-,n,F} & =\sum\limits_{j=k_0+1}^{d'}
Y_{-,n,F}^{(j)}(0),\\
{\cY}_{+,n,V}&= \sum\limits_{j=1}^{k_0}Y_{+,n,V}^{(j)}(0),
\quad\quad  {\cY}_{-,n,V} & =\sum\limits_{j=k_0+1}^{d'}
Y_{-,n,V}^{(j)}(0).
\end{eqnarray}
Using \eqref{NEW4.16} and the definition of
$Y_{\pm,n,V}^{(j)}$ in \eqref{DEFYV} we obtain
\begin{equation}
{\cY}_{+,n,V}=
\sum\limits_{j=1}^{k_0}Y_{+,n}(0)Q_j=Y_{+,n}(0)Q=Y_{+,n}(0),
\end{equation}
and similarly, 
\begin{equation}
{\cY}_{-,n,V}=Y_{-,n}(0). 
\end{equation}
Using Remark \ref{look_1} this results in
\begin{equation}\label{f1}
\DT(I+K)=e^{\Theta} \lim_{n\to  \infty}{\det}_{\bbC^d}\big
({\cY}_{+,n,V}+{\cY}_{-,n,V}\big ).
\end{equation}
Since both $Y_{\pm,n,V}^{(j)}$ and
$Y_{\pm,n,F}^{(j)}$ constitute systems of generalized
matrix-valued Jost solutions, Lemma \ref{mix} yields
\begin{equation}\label{f2}
{\det}_{\bbC^d}\big ({\cY}_{+,n,V}+{\cY}_{-,n,V}\big )
={\det}_{\bbC^d}\big
({\cY}_{+,n,F}+{\cY}_{-,n,F}\big ), \quad n\in\bbN.
\end{equation}
By the continuity argument
discussed in Remark \ref{CONTRN} (cf.\ \eqref{contarg}), we get
\begin{equation}\label{f3}
\lim_{n\to  \infty}{\det}_{\bbC^d}\big ({\cY}_{+,n,F}+{\cY}_{-,n,F}\big
) ={\det}_{\bbC^d}\big ({\cY}_{+, F}+{\cY}_{-,F}\big ),
\end{equation}
and next, again by Lemma \ref{mix},
\begin{equation}\label{f4}
{\det}_{\bbC^d}\big ({\cY}_{+, F}+{\cY}_{-,F}\big )={\det}_{\bbC^d}\big
({\cY}_{+}+{\cY}_{-}\big ).
\end{equation}
Combining \eqref{f1}--\eqref{f4} yields
\eqref{NEW7.12},
completing the proof.
\end{proof}

Next, we will give an extension of Lemma \ref{mix} and Theorem
\ref{posl} for the case where the matrix-valued solutions
$\Phi \, Q_j$, $j=1,\dots,d'$, of \eqref{unperturbed} in
Definition \ref{jostsplitn} are replaced by arbitrary solutions of
\eqref{unperturbed} with certain asymptotic properties. We assume
Hypothesis \ref{mnogo} and fix the splitting $\{P_j\}_{j=1}^{d'}$
for \eqref{perteqn} described in \eqref{P++--} and \eqref{P++--A}. Let
$N_j$,
$j=1,\dots,d'$, be any matrices in $\bbC^{d\times d}$ satisfying
the following conditions:
\begin{equation}\label{DEFC}
\begin{split}
(i) & \quad N_j= N_jQ_j=Q_jN_j, \quad j=1,\dots,d'.  \\
(ii) & \quad {\det}_{\bbC^d} (N)\neq
0,\quad\text{where}\quad
N=\sum_{j=1}^{d'}N_j.
\end{split}
\end{equation}
Let $\Phi_j$
for $j=1,\dots,k_0$, respectively, $j=k_0+1,\dots,d'$, denote the
matrix-valued solutions of the unperturbed equation 
\eqref{unperturbed} on $\bbR_+$, respectively, $\bbR_-$,
satisfying the initial conditions $\Phi_j(0)=N_j$, $j=1,\dots,d'$.
By $(i)$ in \eqref{DEFC}, the Lyapunov exponents of the
$\bbC^d$-valued solutions of \eqref{unperturbed} given by the
columns of $\Phi_j$ belong to the Bohl segments
corresponding to the projections $Q_j$, $j=1,\dots,d'$. Next, assume
that $\widetilde{Y}_+^{(j)}$, $j=1,2,\dots, k_0$, on
$\bbR_+$, and $\widetilde{Y}_-^{(j)}$, $j=k_0+1,2,\dots,
d'$, on $\bbR_-$, are any given $d \times d$ matrix-valued
solutions of the differential equation $ Y'(x)=(A(x)+R(x))Y(x)$,
$x\in \bbR_\pm$, satisfying the following properties (cf.\
Definition \ref{jostsplitn}):
\begin{align}\label{lapunnNEW}
& \limsup_{x\to  \infty }
\frac{\log\|\widetilde{Y}_+^{(j)}(x)-\Phi_j(x)
\|_{\bbC^{d\times d}}}{x}<\varkappa_+'(Q_j),\\
\label{qq+NEW} &
\widetilde{Y}_+^{(j)}(0)=\widetilde{Y}_+^{(j)}(0)Q_j,\quad
j=1,2,\dots, k_0, \\
\intertext{and}
\label{lapunn-NEW} & \liminf_{x\to  -\infty }
\frac{\log\|\widetilde{Y}_-^{(j)}(x)-\Phi_j(x) \|_{\bbC^{d\times
d}}}{x}>\varkappa_-(Q_j),\\
\label{qq-NEW} &
\widetilde{Y}_-^{(j)}(0)=\widetilde{Y}_-^{(j)}(0)Q_j, \quad
 j=k_0+1,\dots, d'.
\end{align}
Define the following matrix $\widetilde{N}$ (cf.\ Definition
\ref{DEFEVANSDET}):
\begin{equation}\label{DEFDDNEW}
\widetilde{N}=\widetilde{\cY}_++\widetilde{\cY}_-,\,\text{
where }\, \widetilde{\cY}_+= \sum_{j=1}^{k_0}
\widetilde{Y}_+^{(j)}(0), \quad \widetilde{\cY}_-
=\sum_{j=k_0+1}^{d'} \widetilde{Y}_-^{(j)}(0).
\end{equation}
We remark that if $N_j=Q_j$ in \eqref{DEFC}, then ${\det}_{\bbC^d} (N)=1$
and $\wti{Y}_\pm^{(j)}$ are generalized matrix-valued Jost
solutions $Y_\pm^{(j)}$ in the sense of Definition
\ref{jostsplitn}, and $\wti{N}=\mathcal{Y}_++\mathcal{Y}_-$, in
the notations of \eqref{DEFDD}.

\begin{lemma}\label{mixNEW}
Assume Hypothesis \ref{mnogo}. In addition suppose that
$N_j$ satisfy \eqref{DEFC}, that $\Phi_j$ satisfy $\Phi_j(0)=N_j$, and
that $\tilde{Y}_\pm^{(j)}$ satisfy conditions
\eqref{lapunnNEW}--\eqref{qq-NEW} for $j=1,\dots,d'$. Finally, assume
that $\wti{N}$ is defined as in  \eqref{DEFDDNEW}. Then
the ratio of the determinants ${\det}_{\bbC^d} (\widetilde{N})$ and
${\det}_{\bbC^d} (N)$ is independent of the choice of the solutions
$\Phi_j$ with the initial data satisfying \eqref{DEFC} and the choice of
the solutions $\widetilde{Y}_\pm^{(j)}$, $j=1,\dots,d'$. Thus, one obtains
for the Evans determinant $D$ in Definition \ref{DEFEVANSDET},
\begin{equation}
D=\frac{{\det}_{\bbC^d}(\widetilde{N})}{{\det}_{\bbC^d} (N)}.  \lb{7.46}
\end{equation}
\end{lemma}
\begin{proof}
Using $(ii)$ in \eqref{DEFC}, we denote
${Y}_\pm^{(j)}(x)=\widetilde{Y}_\pm^{(j)}(x)N^{-1}$,
$x\in\bbR_\pm$, for $j=1,\dots,d'$. We
claim that ${Y}_\pm^{(j)}$ are generalized matrix-valued
Jost solutions in the sense of Definition \ref{jostsplitn}.
Indeed,  since $Q_j=N_jN^{-1}$ due to $N_j=Q_jN$, the equalities
\begin{equation}
{Y}_+^{(j)}(x)-\Phi(x)Q_j=
\widetilde{Y}_+^{(j)}(x)N^{-1}-\Phi(x)N_jN^{-1}=
\big(\widetilde{Y}_+^{(j)}(x)-\Phi_j(x)\big)N^{-1}
\end{equation}
show that
condition \eqref{lapunn} follows from \eqref{lapunnNEW}. Also,
since $N^{-1}$ and $Q_j$ commute, condition \eqref{qq+NEW} implies
\eqref{qq+}, proving the claim. By Lemma \ref{mix}, the
Evans determinant $D$ constructed using the solutions
${Y}_\pm^{(j)}$ as described in Definition
\ref{DEFEVANSDET}, is independent of the choice of $N_j$ and
$\widetilde{Y}_\pm^{(j)}$, and hence equals
\begin{align}
\begin{split}
D&={\det}_{\bbC^d} ({\cY}_++{\cY}_-)={\det}_{\bbC^d}\bigg(
\sum_{j=1}^{k_0} \widetilde{Y}_+^{(j)}(0)N^{-1}+
\sum_{j=k_0+1}^{d'}
\widetilde{Y}_-^{(j)}(0)N^{-1}\bigg)  \\
&={\det}_{\bbC^d}(\widetilde{N}N^{-1}),
\end{split}
\end{align}
proving \eqref{7.46}.
\end{proof}

\begin{corollary}\label{RATIO} Under the assumptions imposed in Theorem
\ref{posl} and Lemma \ref{mixNEW}, the 2-modified perturbation
determinant $\DT(I+K)$ admits the representation
\begin{equation}
\DT(I+K)=e^{\Theta} \frac{{\det}_{\bbC^d}(\widetilde{N})}{{\det}_{\bbC^d}
(N)},
\end{equation}
where $\Theta$ is defined in \eqref{theta}, and $N$ and $\widetilde{N}$
are defined in \eqref{DEFC} and \eqref{DEFDDNEW}.
\end{corollary}

\begin{remark}\label{remLT1} Theorem \ref{+splitn} can be viewed
as a further development of the celebrated Levinson theorem (see,
e.g., \cite[Chap.1]{E} and the bibliography cited therein). The
Levinson theorem for the asymptotically diagonal perturbed
equation \eqref{perteqn} on $\mathbb{R}_+$ deals with the
situation when $A(x)=\diag\{\lambda_j(x)\}_{j=1}^d$, $x\in\bbR_+$,
and $\|R\|_{\bbC^{d\times d}}\in L^1(\mathbb{R}_+)$ and asserts the
existence of
$\bbC^d$-valued solutions of \eqref{perteqn} satisfying the
asymptotic relation
\begin{equation}\label{LevR}
y_j(x)\underset{x\to \infty}{=}
(\mathbf{e}_j+o(1))\exp\int_0^x\lambda_j(s)\,ds, \quad j=1,\dots,d,
\end{equation}
where $\mathbf{e}_j$, $j=1,\dots,d$, are the vectors of the
standard basis in $\mathbb{C}^d$. The underlying assumption in the
traditional Levinson theorem (i.e., for \eqref{LevR} to be valid) is that
the following alternative holds: For each pair of integers $j,k
\in\{1,\dots,d\}$,
$j\neq k$,  either
\begin{equation}
\liminf_{x\to \infty}I_{j,k}(0,x)=-\infty \, \text{ and } \, \sup_{x\ge
x'\ge 0}I_{j,k}(x',x)<\infty,
\end{equation}
or
\begin{equation}
\inf_{x\ge x'\ge 0}I_{j,k}(x,x')>-\infty,
\end{equation}
where we denote
$I_{j,k}(x',x)=\int_{x'}^x ds \, \re(\lambda_j (s)-\lambda_k (s))$, $0\leq
x'\leq x$. To compare our results and the traditional Levinson theorem, we
first note that, unlike the assumptions in the Levinson
theorem, we do not assume that $A$ is diagonal. In the
class of diagonal unperturbed equations \eqref{unperturbed} our assumptions
on the unperturbed diagonal system are more special (since we are
interested only in bounded solutions on $\bbR_+$, and assume the
exponential dichotomy on $\bbR$), but in turn we derive more conclusions
in addition to the asymptotic relation \eqref{LevR}, see
assertions $(ii)$--$(iv)$ in Theorem \ref{+splitn}. If, in addition,
we assume that the diagonal unperturbed system has $d$ disjoint
Bohl segments, then the above mentioned alternative holds; thus
our assumptions are stronger than that in the Levinson theorem in
this particular case. However, generally, our hypotheses are more flexible
since we can group the diagonal elements of $A$ related to
the same Bohl segment, thus avoiding conditions on all pairs $j,k$
as in the traditional Levinson theorem. Our assumptions \eqref{falloffn} on
the perturbation are stronger than those in the Levinson theorem, but in
turn (cf.\ \eqref{LevR}), we conclude that the columns of the
generalized matrix-valued Jost solutions approximate the reference
solutions up to terms $o(e^{-\delta x})$ as $x\to \infty$ for a
positive $\delta$ (cf.\ \eqref{lapunn} and \eqref{vano1}). Finally, the
generalized matrix-valued Jost solutions yield formula
\eqref{NEW7.12} for the perturbation determinant, a feature that
is not discussed in the context of the traditional Levinson theorem.
 \hfill$\Diamond$
\end{remark}

\section{Autonomous Perturbed Equations}\label{APERTE}

In this section we treat the case of an autonomous unperturbed
equation \eqref{unperturbed}, that is, we consider the differential
equations
\begin{align}
y'(x) &= Ay(x), \quad x\in\bbR, \label{AUE}  \\
y'(x) &= (A+R(x))y(x), \quad x\in\bbR, \label{APE}  
\end{align}
on $\bbR$, where $A\in\bbC^{d\times d}$ is a constant matrix and $R\in
L^1_{\loc}(\bbR)^{d\times d}$ (cf.\ Examples \ref{CC}, \ref{CC1}),
and \ref{CC2}). Our objective is to show that the conclusions of
Theorems \ref{+splitn} and \ref{posl} hold for
\eqref{AUE} and \eqref{APE} under much weaker assumptions on the
perturbation $R$ than the exponential decay imposed in
\eqref{falloffn} and \eqref{fallooofn}. Specifically, we will merely
assume that $\|R\|_{\bbC^{d\times d}} \in L^1(\bbR)$ when the
eigenvalues of $A$ are semi-simple, or $\|R\|_{\bbC^{d\times d}}$ decays
polynomially when
$A$ has non-diagonal Jordan blocks. Under these weaker
assumptions, the asymptotic behavior of the generalized
matrix-valued Jost solutions will no longer be measured in terms of 
exponential weight factors as indicated in Definition \ref{jostsplitn}.
Instead, we will define the generalized matrix-valued Jost
solutions of \eqref{APE} by means of solutions of a certain mixed
system of Volterra- and Fredholm-type integral equations in the
spirit of Definition \ref{FGjost}.

\begin{hypothesis}\label{HAUTOM}
Assume that $d\ge 2$, $A\in\bbC^{d\times d}$, and $\spec(A)\cap
i\,\bbR=\emptyset$. Let $Q_j$ denote the spectral projection of $A$
corresponding to the spectral subset of the eigenvalues of $A$
having equal real parts denoted by $\varkappa_j$, $j=1,\dots,d'$,
where $2\le d'\le d$. In addition, assume that for some $k_0$,
$1\le k_0\le d'-1$, the inequalities
$\varkappa_1<\dots<\varkappa_{k_0}<0<\varkappa_{k_0+1}<\dots<\varkappa_{d'}$
hold, so that $Q=\sum_{j=1}^{k_0} Q_j$ is the dichotomy projection
for \eqref{AUE} on $\bbR$ and on $\bbR_\pm$.
\end{hypothesis}

Passing to an appropriate coordinate system, we may assume that
the matrix $A$ is in Jordan normal form, and thus for each
$j=1,\dots,d'$ the operator $A|_{\im (Q_j)}$ is represented by a
direct sum of diagonal matrices $\nu I$ and/or matrices $\nu
I+J$, where $\nu\in\spec(A)$, $\re(\nu)=\varkappa_j$, and
$J$ is the matrix of an appropriate size with $1$'s above the main
diagonal and all other entries equal to zero. We introduce
$m_j\in\bbN$ such that $m_j+1$ is equal to the maximal size of the
non-diagonal Jordan blocks of $A|_{\im (Q_j)}$, noting that $m_j\ge
1$ if such blocks exists, and setting $m_j=0$ if all Jordan
blocks are diagonal (i.e., if all eigenvalues $\nu$ of
$A$ with $\re(\nu)=\varkappa_j$ are semi-simple). Clearly, the
matrix exponent $e^{A|_{\im (Q_j)}}$ can be computed explicitly using the
Jordan blocks of $A|_{\im (Q_j)}$ (see, e.g., \cite[Example I.2.5]{EN}),
so that with this notation one has
the estimate
\begin{equation}\label{SGEST}
\|e^{xA|_{\im (Q_j)}}\|\le ce^{\varkappa_jx}(1+|x|)^{m_j},\quad
j=1,\dots,d',\quad x\in\bbR,
\end{equation}
for some constant $c>0$. Finally, we introduce $m=\max_{1\leq
j\leq d'} \{m_j\}$.

Next, we fix $\tau>0$, and consider the following ``mixed'' system of
Volterra- and Fredholm-type integral equations (cf.\ \eqref{FTSYS} with the
weight function $f(x)=(1+|x|)^{m_j}$, $x\in\bbR$, and \eqref{fred} with
$\mu_j=\varkappa_j$):
\begin{align}
& Z_+^{(j)}(x)-(1+|x|)^{-m_j}e^{x(A-\varkappa_j)}Q_j\no \\
&\quad =-\int_x^\infty dx' \, e^{(x-x')(A-\varkappa_j)}
(1+|x|)^{-m_j}\bigg(\sum_{k=j}^{d'}Q_k\bigg)
(1+|x'|)^{m_j}R(x')Z_+^{(j)}(x')  \no\\
&\qquad +\int_\tau^x dx' \, e^{(x-x')(A-\varkappa_j)}
(1+|x|)^{-m_j}\bigg(\sum_{k=1}^{j-1}Q_k\bigg)
(1+|x'|)^{m_j}R(x')Z_+^{(j)}(x'),  \no\\
 &\hspace*{7.7cm}  j=1,\dots,k_0,\quad x\ge\tau,  \label{ZEq+} \\
& Z_-^{(j)}(x)-(1+|x|)^{-m_j}e^{x(A-\varkappa_j)}Q_j  \no \\
&\quad =\int^x_{-\infty} dx' \, e^{(x-x')(A-\varkappa_j)}
(1+|x|)^{-m_j}\bigg(\sum_{k=1}^{j}Q_k\bigg)
(1+|x'|)^{m_j}R(x')Z_-^{(j)}(x')   \no\\
&\qquad -\int^{-\tau}_x dx' \, e^{(x-x')(A-\varkappa_j)}
(1+|x|)^{-m_j}\bigg(\sum_{k=j+1}^{d'}Q_k\bigg)
(1+|x'|)^{m_j}R(x')Z_-^{(j)}(x'),  \no\\
 &\hspace*{6.8cm} j=k_0+1,\dots,d',\quad x\le-\tau. \label{ZEq-}
\end{align}
Here, we set $\sum_{k=1}^{j-1}Q_k=0$ when $j=1$ in
\eqref{ZEq+} and $\sum_{k=j+1}^{d'}Q_k=0$ when $j=d'$ in
\eqref{ZEq-}; thus, the first and the last equations in
\eqref{ZEq+} and \eqref{ZEq-} are of Volterra-type, and the remaining
equations are of Fredholm-type.

\begin{definition}\label{AUTjsp}
Assume Hypothesis \ref{HAUTOM} and $R\in L^1_{\loc}(\bbR)^{d\times d}$.
Then $d\times d$ matrix-valued solutions $Y_+^{(j)}$, $j=1,\dots,
k_0$, on $\bbR_+$ and $Y_-^{(j)}$, $j=k_0+1,\dots, d'$, on $\bbR_-$ of the
differential equation
\begin{equation}
Y'(x)=(A+R(x))Y(x), \quad x\in \bbR_\pm,
\end{equation}
are called {\em generalized matrix-valued Jost
solutions} associated with the exponential splitting
$\{Q_j\}_{j=1}^{d'}$ if
\begin{align}
Y_+^{(j)}(x)&=(1+|x|)^{m_j}e^{\varkappa_jx}Z_+^{(j)}(x),\quad
j=1,\dots,k_0,\quad x\ge\tau, \label{RESCAL1}\\
Y_-^{(j)}(x)&=(1+|x|)^{m_j}e^{\varkappa_jx}Z_-^{(j)}(x),\quad
j=k_0+1,\dots,d',\quad x\le-\tau, \label{RESCAL2}
\end{align}
where $Z_+^{(j)}$, $j=1,\dots,k_0$, are bounded solutions of equation
\eqref{ZEq+} on $[\tau,\infty)$ and $Z_-^{(j)}$, $j=k_0+1,\dots,d'$, are
bounded solutions of equation \eqref{ZEq-} on $(-\infty,\tau]$. The
solutions $Y_+^{(j)}$ and $Y_-^{(j)}$ are extended to $[0,\tau)$ and to
$(-\tau,0]$, respectively, by solving the initial value problem for
equation \eqref{APE} with initial data $Y_\pm^{(j)}(\pm\tau)=
(1+\tau)^{m_j}e^{\pm\varkappa_j\tau}Z_\pm^{(j)}(\pm\tau)$ for the
corresponding values of $j=1,\dots,d'$.
\end{definition}

We emphasize that the generalized matrix-valued Jost solutions in the
sense of Definition \ref{AUTjsp} are not unique (unless $d'=2$,
see Definition \ref{FGjost}) and depend on the choice of $\tau$ in
\eqref{ZEq+} and \eqref{ZEq-}. We will continue to use Definition
\ref{DEFEVANSDET} of the Evans determinant, but the generalized
matrix-valued Jost solutions in this definition will be understood
in the sense of Definition \ref{AUTjsp}. We recall the exponential
splitting $\{P_j\}_{j=1}^{d'}$  in \eqref{P++--} and \eqref{P++--A}, for
which the Bohl and Lyapunov exponents for the perturbed equation 
\eqref{APE} are also equal to $\varkappa_j$, $j=1,\dots,d'$.

\begin{theorem}\label{MAINAUT}
Assume Hypothesis \ref{HAUTOM} and the condition
\begin{equation}\label{POLDEC}
\|R\|_{\bbC^{d\times d}}\in L^1(\bbR; (1+|x|)^{2m}\,dx),
\end{equation}
where $m+1$, $m\ge 1$, is the maximal size of the non-diagonal
Jordan blocks of $A$ $($if such blocks exist$)$, or $m=0$ if all
Jordan blocks of $A$ are diagonal.  Then the following conclusions
hold: \\
Assertions $(i)$--$(iii)$ of Theorem \ref{+splitn} hold, where the
generalized matrix-valued Jost  solutions are understood in the
sense of Definition \ref{AUTjsp}. \\
In addition:  \\
$(iv)$ The generalized matrix-valued Jost  solutions satisfy
 \begin{align}
\begin{split}
e^{-\varkappa_jx}\|Y_+^{(j)}(x)-e^{xA}Q_j\|_{\bbC^{d\times d}} &
\underset{x\to \infty}{=} o(1), \quad j=1,\dots,k_0,  \\
e^{-\varkappa_jx}\|Y_-^{(j)}(x)-e^{xA}Q_j\|_{\bbC^{d\times d}} &
\underset{x\to  -\infty}{=} o(1),  \quad j=k_0+1,\dots,d'.
\end{split}
 \end{align}
$(v)$ The  2-modified perturbation determinant $\DT(I+K)$ admits
the representation
\begin{equation}\label{VeryNEW7.12}
\DT(I+K)=e^{\Theta} {\det}_{\bbC^d}({\cY}_++{\cY}_-),
\end{equation}
where
\begin{equation}\label{VeryNEW7.13}
\Theta=\int_0^\infty dx \, {\tr}_{\bbC^d}[QR(x)]
-\int_{-\infty}^0 dx \, {\tr}_{\bbC^d}[(I_d-Q)R(x)],
\end{equation}
$ {\cY}_\pm$ are defined as in \eqref{DEFDD}, and
$D={\det}_{\bbC^d}({\cY}_++{\cY}_-)$ is the Evans determinant $($cf.\
\eqref{DEFDD}$)$.
\end{theorem}
\begin{proof}
We will sketch the proof for the case of $\bbR_+$ referring for
details to the proofs of Theorem \ref{+splitn}, Lemma \ref{mix},
and Theorem \ref{posl}. The proof for the case of $\bbR_-$ is
similar. We fix $j\in\{1,\dots,k_0\}$. Starting the proof of
assertion $(i)$, let $F_{\varkappa_j}^{\tau,j}$ denote the integral
operator on $C_{\rm b}([\tau,\infty))^{d\times d}$
defined by the right-hand side of \eqref{ZEq+} (cf.\ \eqref{FKERJ}).

We claim that under assumption \eqref{POLDEC} the following
assertions hold:
\begin{enumerate}\item[$(a)$] $F_{\varkappa_j}^{\tau,j}$ is a contraction
 on $C_{\rm b}([\tau,\infty))^{d\times d}$ for $\tau$
sufficiently large.
 \item[$(b)$] There exists a function $g\in C_{\rm b}([\tau,\infty))$
 so that $\lim_{x\to \infty}g(x)=0$ and
\begin{equation}
\|(F_{\varkappa_j}^{\tau,j}Z)(x)\|_{\bbC^{d\times d}}\le
 g(x)(1+|x|)^{-m_j}\|Z\|_{C_{\rm b}([\tau,\infty))^{d\times d}}
\end{equation}
 for all $x\in[\tau,\infty)$ and $Z\in C_{\rm b}([\tau,\infty))^{d\times
d}$.
 \end{enumerate}
To prove the claims $(a)$ and $(b)$ we use \eqref{SGEST} and bound the norm
of the first integral in \eqref{ZEq+} by the expression
\begin{align}
& c\|
Z_+^{(j)}\|_{C_{\rm b}([\tau,\infty))^{d\times
d}}(1+|x|)^{-m_j}\label{ERZ1} \no \\
&\qquad\times
\int_x^\infty dx' \, \sum_{k=j}^{d'}
e^{(\varkappa_k-\varkappa_j)(x-x')}(1+|x-x'|)^{m_k}
(1+|x'|)^{m_j}\|R(x')\|_{\bbC^{d\times d}} \no \\
& \quad \le c\|
Z_+^{(j)}\|_{C_{\rm b}([\tau,\infty))^{d\times
d}}\int_x^\infty dx' \, (1+|x'|)^{2m_j} \|R(x')\|_{\bbC^{d\times d}},
\end{align}
since, due to $x\le x'$, for $k\ge j+1$ in the sum above, we have
$\varkappa_k-\varkappa_j>0$ and
$\exp((\varkappa_k-\varkappa_j)(x-x'))(1+|x-x'|)^{m_k}\le c$,
while for $k=j$ we have $(1+|x-x'|)^{m_j}\le c(1+|x'|)^{m_j}$.
Similarly, the norm of the second integral in \eqref{ZEq+} is
dominated by the expression
\begin{align}
&c\| Z_+^{(j)}\|_{C_{\rm b}([\tau,\infty))^{d\times
d}}(1+|x|)^{-m_j} \no  \\
&\qquad\times \int_\tau^x dx' \, \sum_{k=1}^{j-1}
e^{(\varkappa_k-\varkappa_j)(x-x')}(1+|x-x'|)^{m_k} (1+|x'|)^{m_j}
\|R(x')\|_{\bbC^{d\times d}}  \no  \\
& \quad \le c\|
Z_+^{(j)}\|_{C_{\rm b}([\tau,\infty))^{d\times
d}}  \no \\
& \qquad \times \int^x_\tau dx' \, \sum_{k=1}^{j-1}
e^{(\varkappa_k-\varkappa_j)(x-x')/2}
(1+|x'|)^{2m_j}\|R(x')\|_{\bbC^{d\times d}},  \label{ERZ2}
\end{align}
since $\exp((\varkappa_k-\varkappa_j)(x-x')/2) [1+|x-x'|]^{m_k}\le
c$ due to $\varkappa_k-\varkappa_j<0$ and $x-x'\ge 0$. Now
condition \eqref{POLDEC} and Lemma \ref{nulek} yield claims $(a)$ and
$(b)$.

Using claim $(a)$, assertion $(i)$ follows. Also, the
uniqueness part in assertion $(ii)$ holds since \eqref{ZEq+} is a
Volterra integral equation for $j=1$. Assertions \eqref{INCL0}, and $(iii)$
in Theorem \ref{+splitn}, and assertion $(iv)$ in Theorem
\ref{MAINAUT} follow from claim $(b)$ similar to the proof
of Theorem \ref{+splitn}. For instance, to prove $(iii)$, we note
first that claim $(b)$ implies the estimate
\begin{equation}\label{ESTGX}
\|Y_+^{(j)}(x)-e^{xA}Q_j\|_{\bbC^{d\times d}}\le cg(x)e^{\varkappa_jx},\quad
x\ge\tau,\quad\text{where}\quad\lim_{x\to \infty} g(x)=0.
\end{equation}
If $y$ is a nonzero $\bbC^d$-valued solution of \eqref{APE} with
$y(0)\in\im(Y_+^{(j)}(0))$ then, using \eqref{ESTGX}, its Lyapunov
exponent satisfies the inequality
\begin{align}
&\lambda_+(y)\le\max\bigg\{\limsup_{x\to \infty}
\frac{\log\|Y_+^{(j)}(x)-e^{xA}Q_j\|_{\bbC^{d\times
d}}}{x},\varkappa_j\bigg\} <\varkappa_{j+1}=\lambda'_+(P_{j+1}), \no \\
& \hspace*{8cm}  j=2,\dots,k_0-1,   \\
& \lambda_+(y)<0, \quad j=k_0.
\end{align}
Applying Remark \ref{VVSOLS} to \eqref{APE} we conclude that
$y(0)\in\im\big(\sum_{k=1}^jP_k\big)$ for
$j=1,\dots,k_0$. However, for $j=2,\dots,k_0$, the assumption
$y(0)=y_0\in\im\big(\sum_{k=1}^{j-1}P_k\big)$ leads to a contradiction
similar to the proof of Theorem \ref{+splitn}\,$(iii)$. Indeed,
using the same notation as in that proof, we have
\begin{equation}
\|e^{xA}Q_jy_1\|_{\bbC^{d}}\le c\|e^{xA}Q_j-Y_+^{(j)}(x)\|_{\bbC^{d\times
d}}+\|y(x)\|_{\bbC^{d}}
\le cg(x)e^{\varkappa_jx}+ce^{\varkappa_{j-1}x},
\end{equation}
where we used \eqref{ESTGX}, the assumption
$y(0)\in\im\big(\sum_{k=1}^{j-1}P_k\big)$, and Remark \ref{VVSOLS} applied
to \eqref{APE}. It follows that
$e^{-\varkappa_jx}\|e^{xA}Q_jy_1\|_{\bbC^{d}}\to  0$ as
$x\to \infty$, in contradiction with the explicit formula for
$e^{xA|_{\im (Q_j)}}$ in terms of the Jordan blocks (cf.\
\cite[Example I.2.5]{EN}). This proves assertion
$(iii)$. Assertion $(iv)$ follows from \eqref{ESTGX}.

It remains to prove assertion $(v)$. For this purpose we consider the
truncated perturbed equation, $y'(x)=(A+R_n(x))y(x)$, $x\in\bbR$,
with $R_n$ as in \eqref{TRUNC}. For a sufficiently large
$\tau>0$, we define the operator $F^{\tau,j}_{\varkappa_j,n}$ on
$C_{\rm b}([\tau,\infty))^{d\times d}$ by the right-hand side of
equation \eqref{ZEq+} with $R$ replaced by $R_n$. Since
$R_n\underset{n\to \infty}{\to} R$ in
$L^1(\bbR; (1+|x|)^{2m}\,dx)^{d\times d}$, estimates similar to
\eqref{ERZ1} and \eqref{ERZ2} show that
$F^{\tau,j}_{\varkappa_j,n}\underset{n\to \infty}{\rightarrow}
F^{\tau,j}_{\varkappa_j}$ in operator norm on the Banach space
$C_{\rm b}([\tau,\infty))^{d\times d}$.
It follows that
\begin{equation}\label{YFTOY}
Y^{(j)}_{\pm,n,F}(0)\underset{n\to \infty}{\to}
Y^{(j)}_{\pm,F}(0),
\end{equation}
where $Y^{(j)}_{\pm,n,F}$ are the generalized matrix-valued Jost
solutions corresponding to the truncated perturbed equation, and
$Y^{(j)}_{\pm}=Y^{(j)}_{\pm,F}$ are the generalized
matrix-valued Jost solutions of \eqref{APE} constructed in part
$(i)$ of the current proof. The rest of the proof is similar to the
arguments in the proof of Theorem \ref{posl}. Indeed, using the
notations introduced in that proof, to establish formula
\eqref{VeryNEW7.12}, we need to show that
${\det}_2(I+K)=e^{\Theta}
{\det}_{\bbC^d}(\mathcal{Y}_{+,F}+\mathcal{Y}_{-,F})$. The latter equality
follows from \eqref{f1}, \eqref{f2}, and
\eqref{f3}. As before, \eqref{f1} follows from Remark
\ref{look_1}, and we already know from \eqref{YFTOY} that
\eqref{f3} holds. Thus, it remains to prove \eqref{f2}. For this 
we will take advantage of the fact that the support of
$R_n$ is compact (and thus we may drop the improper
integral in \eqref{ZEq+} for sufficiently large $x>0$). Indeed, we may
assume that $n>\tau$. It follows from \eqref{volttrunc} and \eqref{DEFYV}
that $Y^{(j)}_{+,n,V}(x)=e^{xA}Q_j$ for $x>n$. On the other hand,
for $x>n$ we also have (cf.\ \eqref{ZEq+})
\begin{equation}
Y^{(j)}_{+,n,F}(x)=e^{xA}Q_j+\int_\tau^n dx' \,
e^{(x-x')A}\bigg(\sum_{k=1}^{j-1}Q_k\bigg)R(x')Y^{(j)}_{+,n,F}(x').
\end{equation}
Thus, using \eqref{SGEST} for $x>n$ and a sufficiently
small $\varepsilon>0$, we obtain the estimate
\begin{align}
& \|Y^{(j)}_{+,n,F}(x)-Y^{(j)}_{+,n,V}(x)\|_{\bbC^{d\times d}}\le
c\int_\tau^n dx' \,
\sum_{k=1}^{j-1}e^{\varkappa_k(x-x')}(1+|x-x'|)^{m_k}  \no \\
&\hspace*{4.7cm} \times\|R(x')\|_{\bbC^{d\times d}}
\|Y^{(j)}_{+,n,F}(x')\|_{\bbC^{d\times d}}  \no \\
 & \quad \le c \int_\tau^n dx' \,
\sum_{k=1}^{j-1}e^{(\varkappa_k-\varepsilon)(x-x')} \max_{x\ge
x'\in\bbR}\Big[e^{(\varkappa_k+\varepsilon)(x-x')}(1+|x-x'|)^{m_k}\Big]
\no  \\
&\qquad \times
\|R(x')\|_{\bbC^{d\times d}}\|Y^{(j)}_{+,n,F}(x')\|_{\bbC^{d\times
d}}  \no \\
& \quad \le ce^{(\varkappa_j-\varepsilon)x} \, \text{ with }\,
c=c(n,\tau),
\end{align}
since $\varkappa_k<\varkappa_j$ for $k=1,\dots,j-1$. Using Remark
\ref{VVSOLS} we therefore obtain the inclusion
$\im\big(Y^{(j)}_{+,n,F}(0)-Y^{(j)}_{+,n,V}(0)\big)
\subseteq\im\big(\sum_{k=1}^{j-1}P_k\big)$, $j=2,\dots,k_0$. Now
\eqref{f2} follows as in the proof of Lemma \ref{mix} since
$Y^{(j)}_{+,n,F}$ and $Y^{(j)}_{+,n,V}$ possess all
properties $(ii)$--$(iv)$ in Theorem \ref{+splitn} (see Remark
\ref{TRIANG}). This concludes the proof of assertion $(v)$ in Theorem
\ref{MAINAUT}.
\end{proof}

\begin{remark}
We emphasize that although the generalized matrix-valued Jost  solutions
in the sense of Definition \ref{AUTjsp} are not uniquely
determined and depend on the choice of $\tau$ in
\eqref{ZEq+} and \eqref{ZEq-}, formula \eqref{VeryNEW7.12} shows that
the Evans determinant $D$ is independent of the choice of the
solutions and hence independent of $\tau$.  \hfill$\Diamond$
\end{remark}

\begin{remark}\label{remLT2}
We note that
Theorem \ref{MAINAUT}\,$(iv)$ also follows from
the Levinson theorem (cf. Remark \ref{remLT1}) for asymptotically
constant coefficients under the additional assumption that the
spectrum of $A$ is simple. In this context we refer to
\cite[Theorem 1.8.1]{E} and also to \cite[Sect.\ 1.10]{E} for additional
results on the asymptotic behavior of solutions of \eqref{perteqn} with $A$
of Jordan block-type.
\hfill$\Diamond$
\end{remark}

\section{The Evans Function}\label{EVFUN}

In this section we relate the Evans determinant $D$ introduced in
Definition \ref{DEFEVANSDET} to the {\em Evans function}, $E$.
First, we recall one of many equivalent definitions of the Evans
function available in the literature (see, e.g.,
\cite{AGJ,BD,JK}), namely, the definition using exponential
dichotomies, see \cite[Definition 4.1]{S}. Consider a family of
differential equations
\begin{equation}
y'(x)=B(z,x)y(x), \quad x\in\bbR,
\end{equation}
parameterized by a complex spectral parameter $z\in\Omega$, where
$\Omega\subseteq\bbC$  is open and simply connected. It is assumed that the
locally integrable $\bbC^{d\times d}$-valued functions $B(z,\cdot)$ depend
on $z\in\Omega$ analytically. Since in this paper we are not
concerned\footnote{See, however, Remark \ref{EDz} below.} with
function theoretic issues related to the Evans function such as
its analytic continuations, etc., we will fix a
value $z_0\in\Omega$ in what follows, and will give the definition of
$E=E(z_0)$ for this $z_0\in\Omega$. Accordingly, we will suppress the
$z$-dependence in the notations for the differential equations above, and
consider just one differential equation $y'(x)=B(x)y(x)$.

\begin{hypothesis}\label{HYPeqnB}
Assume that the differential equation $y'(x)=B(x)y(x)$, $x\in\bbR$, has an
exponentially bounded propagator on $\bbR$ and, in addition, an
exponential dichotomy $P_+$ on $\bbR_+$ and an exponential
dichotomy $P_-$ on $\bbR_-$. Moreover, assume that the ranks of
the projections $P_+$ and $P_-$ are equal, and denote the common
value by $n$ so that
\begin{equation}\label{mEQ} n=\dim(\im (P_+))=\dim(\im (P_-))
\, \text{ and }\, d-n=\dim(\ker (P_+))=\dim(\ker (P_-)).
\end{equation}
\end{hypothesis}

\begin{definition}\label{DEFevfB} Assume Hypothesis \ref{HYPeqnB}.
Choose any vector basis $\{{\rm e}_1,\dots, {\rm e}_n\}$ of the
subspace $\im (P_+)$, and any vector basis $\{{\rm e}_{n+1},\dots,
{\rm e}_d\}$ of the subspace $\ker (P_-)$. The {\it Evans function},
$E$, is then defined as the determinant of the $d\times d$ matrix
with columns ${\rm e}_k$, $k=1,\dots,d$.
\end{definition}

Clearly, $E=0$ if and only if the subspaces $\im (P_+)$ and $\ker
(P_-)$ have a nonzero intersection. $\bbC^d$-valued solutions
$y^{(k)}_+$ on $\bbR_+$ with initial data
$y^{(k)}_+(0)={\rm e}_k$, $k=1,\dots,n$, form a basis in the
subspace of all bounded solutions of the differential
equation $y'(x)=B(x)y(x)$, $x\in\bbR_+$ . Similarly, $\bbC^d$-valued
solutions $y^{(k)}_-$ on $\bbR_-$ with initial data
$y^{(k)}_-(0)={\rm e}_k$, $k=n+1,\dots,d$, form a basis in the
subspace of all solutions bounded on $\bbR_-$. Thus, $E=0$ if and
only if the differential equation has a $\bbC^d$-valued solution
bounded on $\bbR$; see the Introduction. We emphasize, that the Evans
function $E$ just defined is not unique, and depends on the choice of the
vectors ${\rm e}_k$, $k=1,\dots,d$. Moreover, this definition does not
assume any perturbation structure in the differential equation,
and, in particular, can be used to define the Evans functions for
both \eqref{unperturbed} and \eqref{perteqn}.

Returning to the principal theme of the current paper, we will set
$B(x)=A(x)+R(x)$, $x\in\bbR$, with
the matrix-valued functions $A$ and $R$ as in \eqref{unperturbed} and 
\eqref{perteqn}. Next, we intend to show that, in fact, $D=E$,
where $D$ is the Evans determinant in Definition \ref{DEFEVANSDET}, and
$E$ is the Evans function defined using a special choice of the
vectors ${\rm e}_k$, $k=1,\dots,d$, given by the initial data of the
generalized
matrix-valued Jost solutions of \eqref{perteqn}.

For this, we assume in accordance with Hypothesis \ref{diffeq}, that the
unperturbed $d\times d$ matrix differential equation
$y'(x)=A(x)y(x)$, $x\in\bbR$, with $A \in
L^1_{\text{loc}}(\bbR)^{d\times d}$, has an exponentially bounded
propagator and an exponential dichotomy $Q$ on $\bbR$. Assume, in
addition (see Hypothesis \ref{mnogo}), that the unperturbed
equation \eqref{unperturbed} has an exponential splitting of order $d'$,
$2\le d'\le d$, so that $\bbC^d=\im (Q_1)\dot +\cdots\dot +\im (Q_{d'})$,
and that the Bohl segments corresponding to the projections $Q_j$ are
disjoint. We also have $Q=\sum_{j=1}^{k_0}Q_j$ and
$I-Q=\sum_{j=k_0+1}^{d'}Q_j$ for some $1<k_0\le d'$. Suppose that
a $d\times d$ matrix-valued function $R$ is such that
$\|R\|_{\bbC^{d\times d}}\in L^1(\bbR)$. Under these assumptions, due to
Lemma \ref{stub_dich} and its analog for $\bbR_-$, we conclude that the
perturbed equation $y'(x)=(A(x)+R(x))y(x)$, $x\in\bbR_\pm$, has an
exponential dichotomy $P_+$ on $\bbR_+$ and an exponential
dichotomy $P_-$ on $\bbR_-$ such that $\dim(\im (Q))=\dim(\im
(P_+))$ and $\dim(\ker (Q))=\dim(\ker (P_-))$. Thus, \eqref{mEQ}
holds if we let
\begin{equation}\label{DEFn}
n=\dim(\im (Q))=\sum_{j=1}^{k_0}\dim(\im (Q_j)). 
\end{equation}
The same conclusions hold provided the unperturbed equation
\eqref{unperturbed} is autonomous, and we assume Hypothesis \ref{HAUTOM}.
We summarize our assumptions as follows.

\begin{hypothesis}\label{ALLHYP}
For the unperturbed equation \eqref{unperturbed}, respectively,
\eqref{AUE} assume Hypotheses \ref{diffeq} and \ref{mnogo},
respectively, Hypothesis \ref{HAUTOM}. Also, assume that
$\|R\|_{\bbC^{d\times d}}\in L^1(\bbR)$. In addition, assume that
the perturbed equation \eqref{perteqn}, respectively, 
\eqref{APE} has a system of generalized matrix-valued Jost
solutions $Y_+^{(j)}$, $j=1,\dots, k_0$, on $\bbR_+$ and
$Y_-^{(j)}$, $j=k_0+1,\dots, d'$, on $\bbR_-$ in the sense of
Definition \ref{jostsplitn}, respectively, in the sense of
Definition \ref{AUTjsp}.
\end{hypothesis}

Our final result then reads as follows.

\begin{theorem}\label{EVANSFREDHOLM}
 Assume Hypotheses \ref{ALLHYP}. Define the
matrices $\mathcal{Y}_\pm$ and the Evans determinant $D$ as
indicated in Definition \ref{DEFEVANSDET}. In addition, with $n$ as in 
\eqref{DEFn}, let ${\rm e}_k$, $k=1,\dots, n$, denote the nonzero
columns of the matrix $\mathcal{Y}_+$ and ${\rm e}_k$,
$k=n+1,\dots, d$, denote the nonzero columns of the matrix
$\mathcal{Y}_-$. Letting $B(x)=A(x)+R(x)$, $x\in\bbR$, and using
the vectors ${\rm e}_k$, $k=1,\dots, d$, define the Evans function
$E$ for the perturbed equation \eqref{perteqn} as indicated in Definition
\ref{DEFevfB}. Then 
\begin{equation}
E=D. 
\end{equation}
Moreover, if, in addition,  the perturbation
$R$ satisfies the assumptions in Theorem \ref{posl},
respectively, Theorem \ref{MAINAUT},  then the $2$-modified
perturbation determinant $\DT(I+K)$ admits the representation
\begin{equation}\label{NEW9.3}
\DT(I+K)=e^{\Theta} E,
\end{equation}
where $\Theta$ is defined in \eqref{theta}, respectively, in
\eqref{VeryNEW7.13}.
\end{theorem}
\begin{proof}
We claim that $n$ nonzero columns $y_k(x)$ of the $d\times d$
matrix ${\cY}_+(x)=\sum_{j=1}^{k_0}Y_+^{(j)}(x)$, $x\ge 0$,
provide a basis in the space of all bounded $\bbC^d$-valued
solutions of the perturbed equation \eqref{perteqn} on $\bbR_+$, while
$d-n$ nonzero columns $y_k(x)$ of the matrix
${\cY}_-(x)=\sum_{j=k_0+1}^{d'}Y_-^{(j)}(x)$, $x\le  0$, yield a
basis in the space of all bounded $\bbC^d$-valued solutions of the
perturbed equation $y'(x)=(A(x)+R(x))y(x)$ on $\bbR_-$. Indeed, the claim
follows from \eqref{lapunn} and \eqref{lapunn-} (or, for the
autonomous unperturbed equation  \eqref{unperturbed}, from assertion $(iv)$
in Theorem
\ref{MAINAUT})\, and the fact that the nonzero columns of the
matrix $\Phi(x)Q$, respectively, $\Phi(x)(I_d-Q)$, provide a basis
in the space of all bounded solutions of the unperturbed equation
\eqref{unperturbed}  on $\bbR_+$, respectively, $\bbR_-$.

As a result, letting ${\rm e}_k=y_k(0)$, $k=1,\dots,d$, we
selected a basis ${\rm e}_1,\dots, {\rm e}_n$ in the subspace $\im
(P_+)$, and a basis ${\rm e}_{n+1},\dots, {\rm e}_d$ in the
subspace $\ker (P_-)$ needed for the definition of the Evans function $E$
as above. In particular, this proves that the Evans determinant
$D$ equals the Evans function $E$ defined with this choice of bases. The
result in  Theorem \ref{posl} (or, for the autonomous unperturbed equation
\eqref{unperturbed}, Theorem
\ref{MAINAUT}\,$(v)$) now states that, under additional assumptions on $R$
(see \eqref{falloffn} and \eqref{fallooofn} or \eqref{POLDEC}), the
$2$-modified perturbation determinant $\DT(I+K)$ admits the
representation \eqref{NEW9.3}.
\end{proof}

\begin{remark}\label{EDz} Assume that the matrix-valued
functions $A(z,\cdot)$ and
$R(z,\cdot)$ in \eqref{perteqn} analytically depend
on a spectral parameter $z\in\Omega\subseteq\bbC$. Moreover,
assume that the equation $y'(x)=A(x,z)y(x)$, $x\in\bbR$,
has an exponential dichotomy on $\bbR$ for each $z\in\Omega$.
Then the dichotomy projection $Q=Q(z)$ for \eqref{unperturbed}
on $\bbR$ is an analytic $\bbC^{d\times d}$-valued function on $\Omega$
as well. This follows, for instance, from the discussion in Remark
\ref{invertNEW} (alternatively, see \cite{SS2} or \cite[p.\ 995]{S}). Thus,
$\Theta(z)$ in \eqref{theta} is an analytic function for $z\in\Omega$. 
Moreover, the corresponding function $K(z)$ is also analytic as a
function with values in the space of Hilbert--Schmidt operators
equipped with the $\cB(L^2(\bbR)^d)$-norm. 
Hence, the 2-modified perturbation determinant
$\DT(I+K(z))$ is also analytic for $z\in\Omega$. As a final result, it
follows from formula \eqref{NEW7.12} in Theorem \ref{posl} (or from
formula \eqref{VeryNEW7.12} in Theorem \ref{MAINAUT}) that the
Evans determinant $D(z)$ also depends on $z\in\Omega$ analytically. In
view of the discussion in this section we just proved that the
Evans function $E(z)$ is an analytic function of the spectral
parameter $z\in\Omega$ (cf.\ also \cite{AGJ} for an entirely different
proof of this result).
\hfill$\Diamond$
\end{remark}

\begin{remark} The assertions in Lemma \ref{mixNEW} and Corollary \ref{RATIO}
(which are more general than Lemma \ref{mix} and Theorem
\ref{posl}) can be used in the context of this section as
follows. Let ${\rm e}_i^{(j)}$, $i=1,\dots, \dim(\im (Q_j))$, denote
an arbitrary basis in $\im (Q_j)$, $j=1,\dots,d'$, and let $N$
denote the $d\times d$ matrix with the columns ${\rm e}_i^{(j)}$.
For brevity, we enumerate the columns as ${\rm e}_k$,
$k=1,\dots,d$. We note that $\{{\rm e}_k\}_{k=1}^d$ is a basis in
$\bbC^d$. Passing from the standard basis in $\bbC^d$ (i.e.,
$(1,0,\dots,0), \dots, (0,\dots,0,1)$) to this new basis,
we observe that for each $j=1,\dots,d'$, the matrix
$\big(q^{(j)}_{\ell k}\big)_{1\leq \ell,k\leq d}$ corresponding to the
operator $Q_j$ in the basis $\{{\rm e}_k\}_{k=1}^d$ is a diagonal
matrix with all entries being equal to zero except the entries
$q_{\ell\ell}^{(j)}=1$ for $\ell=n_{j-1}+1,\dots,n_{j}$. Here we
denote $n_0=0$ and $n_j=\sum_{k=1}^{j}\dim(\im (Q_k))$ for
$j=1,\dots,d'$. Next, we introduce the $\bbC^d$-valued solutions of the
unperturbed equation \eqref{unperturbed} by $y_k(x)=\Phi(x){\rm e}_k$,
$x\in\bbR$, $k=1,\dots, d$. Let $\tilde y_k(x)$, $x\in\bbR$, $k=1,\dots,
d$, be (not necessarily unique) solutions of the perturbed equation
\eqref{perteqn} that satisfy the asymptotic conditions
\begin{align}
\begin{split}
\widetilde y_k(x)&\underset{x\to  \infty}{=}
y_k(x)+o\big(e^{\varkappa_+'(Q)x}\big),\quad k=1, \dots ,n,  \\
\widetilde y_k(x) &\underset{x\to  -\infty}{=}
 =y_k(x)+o\big(e^{\varkappa_-(I_d-Q)x}\big),\quad k=n+1, \dots, d,
\end{split}
\end{align}
where $n=\dim(\im (Q))$. Now consider two Evans
functions, $E_{\{\widetilde y_k\}}$ and $E_{\{ y_k\}}$ for the
perturbed and unperturbed equations, \eqref{perteqn} and
\eqref{unperturbed}, respectively. Here,
$E_{\{ y_k\}}={\det}_{\bbC^d} (N)$, and $E_{\{\widetilde
y_k\}}={\det}_{\bbC^d}(\widetilde{N})$ with the matrix $\widetilde{N}$
defined in \eqref{DEFDDNEW} and constructed using the matrices
$N_j=NQ_j=Q_jN$,
$j=1,\dots,d'$. The representation
\begin{equation}
\DT(I+K)=e^{\Theta}\frac{E_{\{\widetilde y_k\}}}{E_{\{ y_k\}}}
\end{equation}
now follows from Corollary
\ref{RATIO}.\hfill$\Diamond$
\end{remark}

\appendix
\section{Operators with Semi-Separable Kernels}\label{AP1}

In this appendix  we recall some of the results derived in \cite{GM}. Let
$d_1,d_2\in\bbN$, $d=d_1+d_2$, and let $f_j$ and
$g_j$, $j=1,2$, be given matrix-valued functions on $\bbR$
satisfying
\begin{equation}
f_j\in L^2(\bbR)^{d\times d_j}, \quad  g_j \in
L^2(\bbR)^{d_j\times d}, \quad j=1,2. \lb{2.1}
\end{equation}
On the space $L^2(\bbR)^d$ we consider the Hilbert--Schmidt integral
operator defined by
\begin{equation}
(Ku)(x)=\int\limits_\bbR dx' \, K(x,x')u(x'), \quad u\in
L^2(\bbR)^{d},
\end{equation}
where $K(\cdot,\cdot)$ is a semi-separable
$d\times d$ matrix-valued integral kernel defined by
\begin{equation}
K(x,x')=\begin{cases} f_1(x)g_1(x'), & x\ge x', \\
f_2(x)g_2(x'), & x<x', \end{cases}, \quad x, x' \in \bbR.  \lb{2.3}
\end{equation}
In addition, we introduce the $d\times d$ matrix-valued integral
kernel
\begin{equation}
H(x,x')=f_1(x)g_1(x')-f_2(x)g_2(x'), \lb{2.6}
\end{equation}
and the corresponding Volterra integral equations
\begin{align}
\hat f_1(x )&=f_1(x)- \int_x^\infty dx' \, H(x,x')\hat
f_1(x' ), \lb{2.35} \\
\hat f_2(x )&=f_2(x)+ \int_{-\infty}^x dx' \, H(x,x')\hat f_2(x' ).
\lb{2.36}
\end{align}
We note that equations \eqref{2.35} and \eqref{2.36} have a unique pair of
solutions satisfying $\hat f_j \in L^2(\bbR)^{d\times d_j}$, $j=1,2$.
Finally, we introduce the $d\times d$ matrix-valued function $B$ by
\begin{equation}
B(x)= \begin{pmatrix} g_1(x)f_1(x) & g_1(x)f_2(x) \\
-g_2(x)f_1(x) & -g_2(x)f_2(x)  \end{pmatrix},\quad x\in\bbR.
\lb{2.16}
\end{equation}

\begin{theorem} [\cite{GM}] \label{t2.6} 
Assume \eqref{2.1}. Then:  \\
$(i)$ The first-order $d\times d$ matrix differential equation
\begin{equation}
U'(x)= B(x)U(x),  \quad x\in\bbR,
\end{equation}
permits an explicit particular solution given by the formula
\begin{equation}
U(x)=
 \begin{pmatrix} I_{d_1}-\int_x^\infty dx' \, g_1(x')\hat
f_1(x' )
& \int_{-\infty}^x dx' \, g_1(x')\hat f_2(x' ) \\
\int_x^\infty dx' \, g_2(x')\hat f_1(x' ) & I_{d_2}-
\int_{-\infty}^x dx' \, g_2(x')\hat f_2(x')  \end{pmatrix}, \quad x\in\bbR.
\lb{2.37}
\end{equation}

$(ii)$ The modified Fredholm determinant ${\det}_2 (I- K)$ has the
following representation:
\begin{align}
\begin{split}
{\det}_2 (I- K)&= {\det}_{\bbC^{n}}(U(-\infty)) \exp\bigg(\int_\bbR dx \,
 {\tr}_{\bbC^m}(f_1(x)g_1(x))\bigg) \\
&  ={\det}_{\bbC^{n}}(U(\infty)) \exp\bigg( \int_\bbR dx \,
{\tr}_{\bbC^m}(f_2(x)g_2(x))\bigg).
\end{split}
\end{align}
\end{theorem}

By Liouville's formula (cf., e.g., \cite[Theorem IV.1.2]{H}) one infers for
any $x,x_0\in\bbR$,
\begin{equation}
{\det}_{\bbC^d}(U(x))={\det}_{\bbC^d}(U(x_0)) \exp\bigg( \int_{x_0}^x dx' \,
 {\tr}_{\bbC^d}(B(x'))\bigg),
\end{equation}
and thus, for any $x_0\in\bbR$, one has from Theorem \ref{t2.6}\,$(ii)$,
\begin{align}
 {\det}_2 (I - K) &= {\det}_{\bbC^{d}}(U(\infty)) \exp\bigg(
\int_\bbR dx \,
 {\tr}_{\bbC^d}(f_2(x)g_2(x))\bigg) \no \\
& ={\det}_{\bbC^d}(U(x_0)) \exp\bigg(\int_{x_0}^\infty dx \,
{\tr}_{\bbC^d}(g_1(x)f_1(x)-g_2(x)f_2(x))\bigg) \no \\
 &\quad \times \exp\bigg( \int_\bbR dx \,
 {\tr}_{\bbC^d}(f_2(x)g_2(x))\bigg).
\end{align}
This results in the representation
\begin{align}
{\det}_2 (I - K)&={\det}_{\bbC^d}(U(x_0))
\exp\bigg(\int_{-\infty}^{x_0} dx\,
{\tr}_{\bbC^d}(f_2(x)g_2(x)) \no \\
& \hspace*{3.4cm}  + \int_{x_0}^\infty dx\,
{\tr}_{\bbC^d}(f_1(x)g_1(x)) \bigg)  \lb{3.30} \\
& ={\det}_{\bbC^d} (U(x_0))\exp\bigg({\tr}_{\bbC^d}
 \bigg (\int_{-\infty}^{x_0} dx\,
 K(x-0,x)  \no\\
&\hspace*{4.3cm} +\int_{x_0}^\infty dx\, K(x+0,x) \bigg )\bigg).
\end{align}
We note that if $K$ is a trace class operator with a continuous integral
kernel, then  ${\tr}_{\bbC^d}(B(x))=0$, and applying \cite[p.\
1086--87]{DS88} the  Fredholm determinant can be represented in the simpler
form
\begin{equation}
\det (I-K)={\det}_{\bbC^d} (U(x_0))={\det}_{\bbC^d}(U(-\infty))
={\det}_{\bbC^d} (U(\infty)).
\end{equation}

\section{The Proof of Lemmas \ref{stub_dich} and
\ref{stab_split}}\label{AppendB}

\begin{proof}[Proof of Lemma \ref{stub_dich}.] Using rescaling (see
Remark \ref{RESC}), we will assume without loss of generality that
$\varkappa_+(Q)=-\varkappa'_+(I_d-Q)$. We choose $\varkappa$ so that
$0<\varkappa<-\varkappa_+(Q)$. We claim that: \\
$(i)$ $\dim N=\dim (\im (Q))$, \\
and \\
$(ii)$ for a constant $C(\varkappa)>0$ solutions $y$ of
\eqref{perteqn} satisfy:
\begin{align}
\|y(x)\|_{\bbC^{d}}&\le C(\varkappa)e^{-\varkappa(x-x')}
\|y(x')\|_{\bbC^{d}},  \quad
x\ge x'\ge 0,\,\text{ if }\, y(0)\in N,\label{DICH1} \\
\|y(x)\|_{\bbC^{d}}&\le C(\varkappa)e^{\varkappa(x-x')}
\|y(x')\|_{\bbC^{d}},  \quad
x'\ge x\ge 0,\,\text{ if }\, y(0)\in \ker (Q). \label{DICH2}
\end{align}
As soon as these claims are proved, the assertions in the lemma
follow. Indeed, since $\bbC^d=N\oplus\ker (Q)$ by $(i)$ and
\eqref{DICH1}--\eqref{DICH2}, the exponential boundedness of the
propagator of \eqref{unperturbed} implies that $P$ is an exponential
dichotomy for \eqref{perteqn} by the proof of Lemma IV.3.2 in
\cite{DK}. (We note that in the proof of this lemma only the
exponential boundedness of the propagator has been used.) Since
$-\varkappa$ can be chosen arbitrary close to $\varkappa_+(Q)$,
it follows from \eqref{DICH1} and the definition of
$\varkappa_+(P)$ (cf.\ \eqref{DEFVAR}), that
$\varkappa_+(P)\le\varkappa_+(Q)$. Considering \eqref{unperturbed}
as an $L^1$-perturbation of \eqref{perteqn}, we obtain the
opposite inequality; the second equality in \eqref{bohl12} is
proved similarly using \eqref{DICH2}.

Turning to the proof of the claim, let $\tau>0$, and denote by
$F^\tau$ the integral operator on
$C_{\rm b}([\tau,\infty))^d$ with the
integral kernel $F(x,x')R(x')$ for $x,x'\ge\tau$, where $F(x,x')$
is defined in \eqref{KERNF}. By a direct estimate, inequalities
\eqref{bouddich} and condition $\|R\|_{\bbC^{d\times d}}\in
L^1(\bbR_+)$ imply that $F^\tau$ is a contraction on
$C_{\rm b}([\tau,\infty))^d$
for $\tau$ sufficiently large. Thus, the formula
\begin{equation}\label{CORRbound}
y(x)=\Phi(x)q+(F^\tau y)(x),\quad q\in\im (Q),\quad
x\in[\tau,\infty),
\end{equation}
gives a one-to-one correspondence between the 
solutions $y$, $y(0)\in N$, of \eqref{perteqn} bounded on $\bbR_+$, and the
solutions $\Phi q$, $q\in\ker (Q)$, of \eqref{unperturbed} bounded on
$\bbR_+$. In particular,
$(i)$ in the claim above holds. Also, due to \eqref{CORRbound}, the
bounded solutions $y$ of \eqref{perteqn} admit the representation (cf.\
\cite[p.\ 181]{DK})
\begin{equation} \label{REP}
y(x)=\Phi(x)Q\Phi(x')^{-1}y(x')+(F^{x'}y)(x'),
\quad x\ge x'\ge\tau.
\end{equation}
Taking $\varepsilon>0$ such that
$\alpha=-\varkappa_+(Q)-\varepsilon>\varkappa$, and applying the
exponential dichotomy estimates in the right-hand side of \eqref{REP}, we
arrive at the integral inequality
\begin{equation}
\|y(x)\|_{\bbC^{d}}\le
a(\varepsilon)e^{-\alpha(x-x')}\|y(x')\|_{\bbC^{d}}
+c(\varepsilon)\int_{x'}^\infty ds \,
e^{-\alpha|x-s|}\|R(s)\|_{\bbC^{d\times d}}
\|y(s)\|_{\bbC^{d}}.
\end{equation}
Thus, the function
$u(x)=e^{\varkappa x}\|y(x)\|_{\bbC^{d}}$, $x\ge x'\ge\tau$, 
satisfies the inequality
\begin{equation}\label{INEQu} u(x)\le
a(\varepsilon)e^{-\alpha(x-x')+\varkappa
x}\|y(x')\|_{\bbC^{d}}+c(\varepsilon)\int_{x'}^\infty ds \,
e^{-\alpha|x-s|+\varkappa(x-s)}\|R(s)\|_{\bbC^{d\times
d}}u(s).
\end{equation}
Claim \eqref{DICH1} can now be proved using
arguments in
\cite[Sect.\ III.2]{DK}. Indeed, on the space
$C_{\rm b}([\tau,\infty))$ we consider the 
integral operator $T$ with integral kernel
\begin{equation}
T(x,x')=c(\varepsilon)e^{-\alpha|x-x'|+\varkappa(x-x')}
\|R(x')\|_{\bbC^{d\times d}}, \quad x\geq x'\geq \tau.
\end{equation}
Using equation (2.15) in \cite[Sect.\ III.2]{DK} one then shows that
$\|T\|$ is bounded by the expression
\begin{equation}
c(\varepsilon)\big[1+(1-e^{-(\alpha-\varkappa)})^{-1}+
(1-e^{-(\alpha+\varkappa)})^{-1}\big]
\sup_{x\in[\tau,\infty)}\int_x^{x+1} ds \, \|R(s)\|_{\bbC^{d\times
d}}.
\end{equation}
Thus, the operator $T$ is a contraction for
$\tau=\tau(\varkappa,\varepsilon)>0$ suficiently large. Passing to
equality in the inequality
\eqref{INEQu}, solving the resulting integral equation for $u$ by applying
$(I-T)^{-1}$, and using \cite[Lemma III.2.1]{DK}, we derive
\eqref{DICH1} from the estimate
\begin{equation}
u(x)\le
a(\varepsilon)(1-\|T\|)^{-1}e^{\varkappa
x'}\|y(x')\|_{\bbC^{d}}, \quad  x\ge x'\ge\tau.
\end{equation}

Turning to the proof of \eqref{DICH2}, we note that any solution
$y$ of \eqref{perteqn} on $\bbR_+$ with $y(0)\in\ker (Q)$
solves the equation
\begin{align}
\begin{split}
y(x)&=\Phi(x)(I_d-Q)\Phi(x')^{-1}y(x')+
\int_0^x ds \, \Phi(x)Q\Phi(s)^{-1}R(s)y(s)\\
& \quad -\int_x^{x'} ds \, \Phi(x)(I_d-Q)\Phi(s)^{-1}R(s)y(s), \quad
0\le x\le x',
\end{split}
\end{align}
which yields the integral inequality
\begin{align}
\|y(x)\|_{\bbC^{d}} & \le a(\varepsilon)e^{\alpha(x-x')}
\|y(x')\|_{\bbC^{d}}
+c(\varepsilon)\int_0^{x'} ds \, e^{\alpha(x-s)}\|R(s)\|_{\bbC^{d\times d}}
\|y(s)\|_{\bbC^{d}} \no \\
&\quad  +
c(\varepsilon)\int_{x}^{x'} ds \, e^{-\alpha(x-s)}\|R(s)\|_{\bbC^{d\times
d}} \|y(s)\|_{\bbC^{d}} \no \\
& = a(\varepsilon)e^{\alpha
(x-x')}\|y(x')\|_{\bbC^{d\times d}}
+c(\varepsilon)\int_0^{x'} ds \, e^{-\alpha|x-s|}\|R(s)\|_{\bbC^{d\times d}}
\|y(s)\|_{\bbC^{d}} \no  \\
& \hspace*{8cm} 0\le x\le x'.
\end{align}
Arguments similar to the proof of \eqref{DICH1} conclude the proof
of Lemma \ref{stub_dich}.
\end{proof}

\begin{proof}[Proof of Lemma \ref{stab_split}.] We will use the
rescaling from Remark \ref{RESC} as follows: For each $j=1,\dots,d'-1$
we fix
$\mu_j\in (\varkappa_+(Q_j),\varkappa'_+(Q_{j+1}))$, and for $j=d'$ we fix
$\mu_{d'}>\varkappa_+(Q_d)$. Then, for each $j=1,\dots, d'$,
the rescaled unperturbed equation \eqref{RESCeqn} with
$\mu=\mu_j$ has an exponential dichotomy on $\bbR_+$ with 
dichotomy projection $\mathcal{Q}_j=\sum_{k=1}^jQ_k$. By Lemma
\ref{stub_dich}, the rescaled perturbed equation,
$y'(x)=(A(x)+R(x)-\mu_j I_d)y(x)$, $x\in\bbR$, has an exponential dichotomy
$\mathcal{P}_j$, where $\mathcal{P}_j$ is the projection parallel
to $\ker(\mathcal{Q}_j)$ on the subspace $N_j$ consisting of the
values $y(0)$ of the bounded solutions of the
rescaled perturbed equation on $\bbR_+$ . We note that
$\mathcal{P}_{d'}=I_d$ and set $\mathcal{P}_0=0$. Letting
$P_j=\mathcal{P}_j-\mathcal{P}_{j-1}$, $j=1,\dots,d'$, finishes
the proof of Lemma \ref{stab_split}.
\end{proof}

\noindent {\bf Acknowledgments.}
Fritz Gesztesy and Yuri Latushkin gratefully acknowledge a
research leave for the academic year 2005/06 granted by the Research
Council and the Office of Research of the University of Missouri--Columbia.
Moreover, Yuri Latushkin gratefully acknowledges support by the Research
Board of the University of Missouri.


\end{document}